\documentclass[a4,12pt,reqno]{amsart}
\usepackage{amsmath}
\usepackage{amsthm}
\usepackage{graphicx}
\usepackage{amsfonts}
\usepackage{bm}
\usepackage{bbm}
\usepackage{textcomp}
\usepackage{tikz}
\usetikzlibrary{er,positioning}
\usepackage{amssymb}
\usepackage{enumerate}
\usepackage[margin=1.1in]{geometry}
\usepackage{xcolor}
\usepackage{mathrsfs}

\numberwithin{equation}{section}

\newtheorem{theorem}{Theorem}[section]
\newtheorem{lemma}[theorem]{Lemma}
\newtheorem{proposition}[theorem]{Proposition}
\newtheorem{corollary}[theorem]{Corollary}

\theoremstyle{definition}

\newtheorem{definition}[theorem]{Definition}
\newtheorem{remark}[theorem]{Remark}

\def\d{\,{\mathrm d}}
\def\E{{\mathbb E}}

\def\R{{\mathbb R}}
\def\SS{{\mathcal S}}
\def\N{{\mathbb N}}
\def\PP{{\mathbb P}}
\def\FF{{\mathbb F}}
\def\P{{\mathcal P}}

\def\M{{\mathcal M}}
\def\H{{\mathcal H}}

\def\L{{\mathcal L}}
\def\G{{\mathcal G}}

\def\A{{\mathcal A}}

\def\F{{\mathcal F}}
\def\GG{{\mathbb G}}
\def\HH{{\mathbb H}}

\def\L{{\mathcal L}}
\def\gen{ L}

\def\stabtimes{\times_{\mathrm{s}}}

\newcommand{\indep}{\perp \!\!\! \perp}

\begin{document}
	
\title[Superposition and mimicking theorems]{Superposition and mimicking theorems for conditional McKean-Vlasov equations}
\author{Daniel Lacker, Mykhaylo Shkolnikov, and Jiacheng Zhang}
\address{IEOR Department, Columbia University, New York, NY 10027}
\email{daniel.lacker@columbia.edu}
\address{ORFE Department, Bendheim Center for Finance, and Program in Applied \& Computational Mathematics, Princeton University, Princeton, NJ 08544.}
\email{mshkolni@gmail.com}
\address{ORFE Department, Princeton University, Princeton, NJ 08544.}
\email{jiacheng@princeton.edu}
\footnotetext[1]{M.~Shkolnikov is partially supported by the NSF grant DMS-1811723 and a Princeton SEAS innovation research grant.}

\begin{abstract}
We consider conditional McKean-Vlasov stochastic differential equations (SDEs), such as the ones arising in the large-system limit of mean field games and particle systems with mean field interactions when  common noise is present. The conditional time-marginals of the solutions to these SDEs satisfy non-linear stochastic partial differential equations (SPDEs) of the second order, whereas the laws of the conditional time-marginals follow Fokker-Planck equations on the space of probability measures. We prove two superposition principles: The first establishes that any solution of the SPDE can be lifted to a solution of the conditional McKean-Vlasov SDE, and the second guarantees that any solution of the Fokker-Planck equation on the space of probability measures can be lifted to a solution of the SPDE. We use these results to obtain a mimicking theorem which shows that the 
conditional time-marginals of an It\^o process
can be emulated by those of a solution to a conditional McKean-Vlasov SDE with Markovian coefficients. This yields, in particular, a tool for converting open-loop controls into Markovian ones in the context of controlled McKean-Vlasov dynamics.
\end{abstract}

\maketitle

\tableofcontents


\section{Introduction}\label{sec_intro}

We are concerned with stochastic differential equations (SDEs) with random coefficients of the form
\begin{equation}\label{SDE_gen}
	\d\bm{X}_t = \bm{b}(t,\omega,\bm{X}_t)\d t + \bm{\sigma}(t,\omega,\bm{X}_t)\d\bm{W}_t + \bm{\gamma}(t,\omega,\bm{X}_t)\d\bm{B}_t,\quad t\in[0,T].
\end{equation}
This SDE is posed on a filtered probability space $(\Omega,\FF\!=\!(\F_t)_{t\in[0,T]},\PP)$ supporting two independent $d$-dimensional $\FF$-Brownian motions $\bm B$ and $\bm W$, as well as a subfiltration $\GG=\!(\G_t)_{t\in[0,T]}\subset \FF$, with respect to which $\bm{B}$ is adapted and $\bm{W}$ is independent.
Additionally, 
the coefficients $(\bm{b},\bm{\sigma},\bm{\gamma}):\,[0,T]\times\Omega\times\R^d\to\R^d\times\R^{d\times d}\times\R^{d\times d}$ are measurable with respect to the product of the $\GG$-progressive $\sigma$-algebra on $[0,T]\times\Omega$ and the Borel $\sigma$-algebra on $\R^d$. We assume also that
\begin{align}
\E\bigg[\int_{0}^{T}  |\bm b(t,\omega,\bm{X}_t)\big| +  |(\bm\sigma\bm\sigma^\top + \bm\gamma\bm\gamma^\top)(t,\omega,\bm{X}_t)| \d t\bigg]<\infty. \label{integrability-SDEtoSPDE}
\end{align}

\begin{remark} \label{re:MVcase}
Fitting into the framework of \eqref{SDE_gen} is the important special case of conditional McKean-Vlasov SDEs, in which the dependence on $\omega$ is through the conditional law of $\bm{X}_t$ given $\G_t$, denoted $\L(\bm{X}_t\,|\,\G_t)$:
\begin{equation}\label{McV_special}
	\d\bm{X}_t=\bm{b}\big(t,\L(\bm {\bm{X}}_t|\G_t),\bm{X}_t\big)\!\d t + \bm{\sigma}\big(t,\L(\bm {\bm{X}}_t|\G_t),\bm{X}_t\big)\!\d\bm{W}_t + \bm{\gamma}\big(t,\L(\bm {\bm{X}}_t|\G_t),\bm{X}_t\big)\!\d\bm{B}_t,\; t\in[0,T].
\end{equation}
Indeed, if $\bm{X}$ solves \eqref{McV_special}, then by freezing the nonlinear term $\L(\bm {\bm{X}}_t\,|\,\G_t)(\omega)$ we see that $\bm{X}$ also solves an equation of the form \eqref{SDE_gen}.
\end{remark}

The SDE \eqref{McV_special} describes the state process of a representative player in the large-population limit of mean field games and control problems with common noise (see \cite[Section 2.1]{CD2}, and also \cite{GLLpre,GLL,CD,BFY,CFS,LW,Ah,BZ,CDL,La, ARY, HSS,KT}), as well as the dynamics of a particle in the large-scale limit of a particle system with mean field interactions and common noise (see, e.g., \cite{Va,DV,Ko,kurtz1999particle,CF}).
Often, as in \cite{kurtz1999particle}, the filtration $\GG$ is the one generated by $(\bm{B}_t)_{t \in [0,T]}$, but the general case corresponds to a notion of \emph{weak solution}; see \cite{HSS} and \cite[Section 8]{LaSv} for further discussion.

\medskip

\subsection{Superposition from SPDE to SDE}
Our first results iron out the general connection between the SDE \eqref{SDE_gen} and the associated SPDE which should govern the evolution of the conditional measure flow $\mu_t=\L(\bm{X}_t\,|\,\G_t)$:
\begin{equation}\label{SPDE}
\d\langle \mu_t, \varphi\rangle = \big\langle \mu_t, \gen^{\bm b,\bm a}_{t,\omega}\varphi\big\rangle \d t + \big\langle \mu_t,(\nabla\varphi)^\top \bm \gamma(t,\omega,\cdot)\big\rangle \d\bm B_t, \quad t \in [0,T], \quad\varphi \in C_c^\infty(\R^d),
\end{equation}
where $\bm{a} := \bm{\sigma}\bm{\sigma}^\top + \bm{\gamma}\bm{\gamma}^\top$, and the operator $\gen^{\bm b,\bm a}_{t,\omega}$ acts on a smooth function $\varphi$ of compact support by
\begin{equation*}
	\gen^{\bm b,\bm a}_{t,\omega}\varphi:= \bm b(t,\omega,\cdot)\cdot\nabla \varphi +\frac12\,\bm a(t,\omega,\cdot):\nabla^2\varphi.
\end{equation*}
Here $\nabla$ and $\nabla^2=\nabla\nabla^\top$ denote the gradient and Hessian, and $s\cdot r=s^\top r$ and $S : R=\mathrm{Tr}[SR^\top]$ denote the usual inner products in $\R^d$ and $\R^{d \times d}$.

\medskip

We state first a form of a fairly well known proposition which asserts that the conditional measure flow of a solution of the SDE \eqref{SDE_gen} is a solution of the SPDE \eqref{SPDE}. In the following, we use the symbol $\indep$ to stand for ``independent of" or ``is independent of" depending on the grammatical context. In addition, for $\sigma$-algebras $\A_1$, $\A_2$, and $\A_3$, we write $\A_1\indep\A_2\,|\,\A_3$ to mean that $\A_1$ and $\A_2$ are conditionally independent given $\A_3$. For two $\sigma$-algebras $\A_1$ and $\A_2$ we write  $\A_1 \vee \A_2$ for $\sigma(\A_1 \cup \A_2)$ as usual. If $Z=(Z_t)_{t \in [0,T]}$ is a stochastic process, we write $\FF^Z=(\F^Z_t)_{t \in [0,T]}$ for the filtration it generates.

\begin{proposition} \label{pr:SDEtoSPDE}
Suppose $(\Omega,\FF=(\F_t)_{t\in[0,T]},\PP)$ is a filtered probability space supporting independent $d$-dimensional $\FF$-Brownian motions $\bm{B}$ and $\bm{W}$, as well as a subfiltration $\GG =(\G_t)_{t \in [0,T]}$ with respect to which $\bm{B}$ is adapted and $\bm{W}$ is independent. Suppose the triple $(\bm{b},\bm{\sigma},\bm{\gamma}):\,[0,T]\times\Omega\times\R^d\to\R^d\times\R^{d\times d}\times\R^{d\times d}$ is measurable with respect to the product of the $\GG$-progressive $\sigma$-algebra on $[0,T]\times \Omega$ and the Borel $\sigma$-algebra on $\R^d$. Let $\bm{X}=(\bm{X}_t)_{t\in [0,T]} $ be a continuous $\FF$-adapted $\R^d$-valued process satisfying \eqref{SDE_gen} and \eqref{integrability-SDEtoSPDE}, as well as $\F^{\bm{X}}_t \indep \F^{\bm{W}}_T \vee \G_T\,|\,\F^{\bm{W}}_t \vee \G_t$ for each $t \in [0,T]$.
Then the process $(\mu_t := \L(\bm{X}_t \, | \, \G_T))_{t \in [0,T]}$ admits a continuous version, and the following hold:
\begin{enumerate}[(1)]
\item $\mu_t=\L(\bm{X}_t\,|\,\G_t)$ a.s.\ for each $t \in [0,T]$.
\item The SPDE \eqref{SPDE} holds.
\end{enumerate}
\end{proposition}

Proposition \ref{pr:SDEtoSPDE} is well known in various forms, see, e.g., \cite{DV,kurtz1999particle,HSS}. It is fairly straightforward to prove by applying It\^o's formula to $\varphi(\bm{X}_t)$, taking the conditional expectations with respect to $\G_T$, and then interchanging the order of stochastic integration and conditional expectation. We give the proof in Appendix \ref{sec_app_SDEtoSPDE}, because our statement appears to be somewhat more general than its relatives in prior literature, and because there is more subtlety than one might expect in the aforementioned interchange. The latter is one purpose of the conditional independence or \emph{compatibility} assumption, $\F^{\bm{X}}_t \indep \F^{\bm{W}}_T \vee \G_T\,|\,\F^{\bm{W}}_t \vee \G_t$ for each $t \in [0,T]$.
This compatibility condition takes several equivalent forms, compiled in Lemma \ref{le:compatibility}, and is natural from the perspective of prior work on weak solutions of SDEs with random coefficients \cite{jacod1981weak,kurtz2014weak}.

\medskip

Our first main result is the following superposition theorem, which shows that, conversely, each solution of the SPDE \eqref{SPDE} gives rise to a solution of the SDE \eqref{SDE_gen}.

\begin{theorem}\label{Thm_SPDEtoSDE_general}
Suppose $(\Omega,\mathbb G=(\G_t)_{t\in[0,T]},\PP)$ is a filtered probability space supporting a $d$-dimensional $\GG$-Brownian motion $\bm{B}$. Suppose the triple $(\bm{b},\bm{a},\bm{\gamma}):\,[0,T]\times\Omega\times\R^d\to\R^d\times\R^{d\times d}\times\R^{d\times d}$ is measurable with respect to the product of the $\mathbb G$-progressive $\sigma$-algebra on $[0,T]\times \Omega$ and the Borel $\sigma$-algebra on $\R^d$,  $(\mu_t)_{t\in [0,T]}$ is a continuous $\GG$-adapted probability measure-valued process satisfying \eqref{SPDE}, and $\G_T$ is countably generated.
Suppose finally that
\begin{align} 
\E\bigg[\int_{0}^{T} \big\|\bm b(t,\omega,\cdot)\big\|^p_{L_p(\mu_t)}+\big\|\bm a(t,\omega,\cdot)\big\|^p_{L_p(\mu_t)}\,\d t\bigg]<\infty \label{integ}
\end{align}
for some $p>1$,
and $\bm a-\bm\gamma\bm\gamma^\top$ is symmetric and positive semidefinite with symmetric square root $\bm\sigma$. Then there exists an extension $(\widetilde{\Omega},\FF=(\F_t)_{t\in[0,T]},\widetilde{\PP})$ of the probability space $(\Omega,\GG,\PP)$ supporting a $d$-dimensional $\FF$-Brownian motion $\bm{W}$ independent of $\G_T$ and a continuous $\FF$-adapted $d$-dimensional process $\bm{X}$ such that:
\begin{enumerate}[(1)]
\item The SDE \eqref{SDE_gen} holds.
\item $\mu_t=\L(\bm{X}_t\,|\,\G_T)=\L(\bm{X}_t\,|\,\G_t)$ a.s., for each $t \in [0,T]$.
\item $\F^{\bm{X}}_t \indep\F^{\bm{W}}_T \vee \G_T\,|\,\F^{\bm{W}}_t \vee \G_t$, for each $t\in [0,T]$.
\item $\bm{B}$ is an $\FF$-Brownian motion.
\end{enumerate}
\end{theorem}

For classical SDEs with bounded coefficients and the corresponding Fokker-Planck equations, the seminal paper \cite{figalli2008existence} develops a general ``superposition'' theory that associates to each solution of a Fokker-Planck equation a weak solution of the SDE. The results of \cite{figalli2008existence} have been  extended to a wide range of SDEs with unbounded coefficients in \cite{trevisan} (see also  \cite{BRS,RXZ} and the comparison with the previous work \cite{LL} in \cite[Remark 3.4]{trevisan}). The precursors to this line of research can be traced back to the analogous results in \cite{DL,Am} for ordinary differential equations. Moreover, superposition principles for controlled SDEs and two classes of (particularly singular) McKean-Vlasov SDEs can be found in \cite{KuSt} and \cite{BR,BR2}, respectively. Theorem \ref{Thm_SPDEtoSDE_general} gives the first superposition principle for SDEs with \textit{random} coefficients, in which the SPDE \eqref{SPDE} plays the role of the Fokker-Planck equation, and covers, in particular, the conditional McKean-Vlasov SDEs of the special form \eqref{McV_special}.
Our proof of Theorem \ref{Thm_SPDEtoSDE_general} employs, in addition to many of the approximation techniques developed in \cite{figalli2008existence,trevisan}, novel tools in the context of superposition principles to address the randomness of the coefficients: the stable convergence topology, the measurable maximum theorem, and a new strong existence theorem for SDEs with random coefficients. 

\subsection{Superposition from Fokker-Planck equation  on $\P(\R^d)$ to SPDE}

Our next results pertain to the SPDE corresponding to the McKean-Vlasov SDE \eqref{McV_special}: 
\begin{equation}\label{SPDE_special}
	\d\langle \mu_t, \varphi\rangle = \big\langle \mu_t, \gen_{t,\mu_t}^{\bm b,\bm a}\varphi\big\rangle \d t + \big\langle \mu_t,(\nabla\varphi)^\top \bm \gamma(t,\mu_t,\cdot)\big\rangle \d\bm B_t,\quad t \in [0,T], \quad\varphi \in C_c^\infty(\R^d),
\end{equation}
where
\begin{align}
\gen_{t,m}^{\bm b,\bm a}\varphi:=\bm b(t,m,\cdot)\cdot\nabla \varphi +\frac12\,\bm a(t,m,\cdot):\nabla^2\varphi, \quad t\in [0,T], \quad m \in \P(\R^d) \label{def:gen-MV}
\end{align}
and $\bm a = \bm\sigma{\bm\sigma}^\top\!+\bm\gamma\bm\gamma^\top$ as before. We denote by $P_t$ the law $\L(\mu_t)$ of $\mu_t$, for $t\in[0,T]$. Throughout the paper, for a Polish space $E$ we write $\P(E)$ for the space of Borel probability measures on $E$ equipped with the topology of weak convergence and the corresponding Borel $\sigma$-algebra. Then $(\mu_t)_{t \in [0,T]}$ is a continuous $\P(\R^d)$-valued process, and $P=(P_t)_{t \in [0,T]}$ belongs to $C([0,T];\P(\P(\R^d)))$, the space of continuous functions from $[0,T]$ to $\P(\P(\R^d))$.

\medskip

First, let us describe the easier direction, which lies in deriving an equation for $P$ from the equation for $\mu$. For $m \in \P(\R^d)$, $k \in \N$, and a vector of test functions $\bm\varphi=(\varphi_1,\ldots,\varphi_k)\in(C^\infty_c(\R^d))^k$, we write $\langle m,\bm\varphi \rangle := \big(\langle m,\varphi_1\rangle,\ldots,\langle m,\varphi_k\rangle\big) \in \R^k$.
If we assume
\begin{align} 
\E\bigg[\int_{0}^{T} \big\|\bm b(t,\mu_t,\cdot)\big\|_{L_1(\mu_t)}+\big\|\bm a(t,\mu_t,\cdot)\big\|_{L_1(\mu_t)}+\big\|\bm \gamma(t,\mu_t,\cdot)\big\|_{L_1(\mu_t)}^2\,\d t\bigg]<\infty, \label{integrability-SPDEtoPDE}
\end{align}
then, for $f \in C_c^\infty(\R^k)$, we may apply It\^o's formula to $f\big(\langle \mu_t,\bm\varphi\rangle \big)$ and take the expectation to find that $(P_t)_{t\in[0,T]}$ satisfies a Fokker-Planck equation on $\P(\R^d)$: 
\begin{align}
\begin{split}
&\int_{\P(\R^d)} f\big(\langle  m,\bm\varphi\rangle \big) (P_t-P_0)(\!\d m)  = \int_0^t\int_{\P(\R^d)}\left[\sum_{i=1}^k\partial_if\big(\langle  m,\bm\varphi\rangle \big) \,\big\langle  m,\gen^{\bm b,\bm a}_{s, m}\varphi_i\big\rangle \right. \\
& \;\;\;\left. + \frac12\sum_{i,j=1}^k\partial_{ij}f\big(\langle  m,\bm\varphi\rangle \big) \, \big\langle  m,(\nabla \varphi_i)^\top \bm\gamma(s, m,\cdot)\big\rangle\cdot\big\langle  m, (\nabla \varphi_j)^\top \bm\gamma(s, m,\cdot)\big\rangle \right] P_s(\!\d m)\d s, \\
&\qquad\qquad\qquad\qquad\qquad\qquad\;\;\; 
t\in[0,T],\quad k\in\N,\quad f\in C^2_c(\R^k),\quad \bm\varphi\in (C^\infty_c(\R^d))^k.
\end{split} \label{FPE}
\end{align}
Indeed, \eqref{integrability-SPDEtoPDE} ensures that the local martingale in the It\^o expansion of $f\big(\langle \mu_t,\bm\varphi\rangle \big)$ is a true martingale, and the expectations are well-defined.
This proves the following:

\begin{proposition} \label{pr:SPDEtoPDE}
Let $(\bm b,\bm a,\bm \gamma):\,[0,T]\times \P(\R^d)\times\R^d\to \R^d\times\R^{d\times d}\times \R^{d\times d}$ be measurable. Suppose $(\Omega,\mathbb G=(\G_t)_{t\in[0,T]},\PP)$ is a filtered probability space supporting a $d$-dimensional $\GG$-Brownian motion $\bm B$ and a continuous $\GG$-adapted $\P(\R^d)$-valued process $(\mu_t)_{t\in [0,T]}$  satisfying the SPDE \eqref{SPDE_special} and the integrability condition \eqref{integrability-SPDEtoPDE}.
Let $P_t:=\L(\mu_t)$ for $t \in [0,T]$. Then $(P_t)_{t \in [0,T]}$ satisfies the Fokker-Planck equation \eqref{FPE}.
\end{proposition}

Our second main result is the following superposition theorem, which shows that, conversely, each solution of the Fokker-Planck equation \eqref{FPE} in the space $C([0,T];\P(\P(\R^d)))$ stems from a solution of the SPDE \eqref{SPDE_special}. 

\begin{theorem}\label{Thm_PDEtoSPDE}
Suppose that $P\in C([0,T];\P(\P(\R^d)))$ and measurable $(\bm b,\bm a,\bm \gamma):\,[0,T]\times \P(\R^d)\times\R^d\to \R^d\times\R^{d\times d}\times \R^{d\times d}$ are such that \eqref{FPE} holds and 
\begin{align}
\int_{0}^{T}\int_{\P(\R^d)} \big\|\bm b(t, m,\cdot)\big\|^p_{L_p( m)}
		+\big\|\bm a(t, m,\cdot)\big\|^p_{L_p( m)}
		+\big\|\bm \gamma(t, m,\cdot)\big\|^{2p}_{L_p( m)}\, P_t(\!\d m)\d t<\infty \label{integ2}
\end{align}
for some $p > 1$. Then there exists a filtered probability space $(\Omega,\GG,\PP)$, with $\G_T$ countably generated, supporting a $d$-dimensional $\GG$-Brownian motion $\bm B$ and a continuous $\GG$-adapted $\P(\R^d)$-valued process $(\mu_t)_{t\in[0,T]}$ solving \eqref{SPDE_special} with $\L(\mu_t)=P_t$ for each $t\in[0,T]$.
\end{theorem}

\smallskip

Superposition principles in infinite-dimensional spaces, specifically in $\R^\infty$ and general metric measure spaces, have been established recently in \cite{trevisan_diss} (see also \cite{AT,ST} for deterministic counterparts). In Theorem \ref{Thm_PDEtoSPDE}, the infinite-dimensional space in consideration is $\P(\R^d)$, and even if one were to replace $\R^d$ by a compact subset, restrict the attention to the subspace of probability measures with finite second moments and equip the latter with the $2$-Wasserstein distance, the associated metric measure spaces (see \cite{Stu,vRS}) do not seem to admit a $\Gamma$-calculus (cf.~\cite[Remark 5.6]{Gi}) as required for the superposition principle in general metric measure spaces \cite[Theorem 7.3]{trevisan_diss}. Instead, our proof of Theorem \ref{Thm_PDEtoSPDE}  relies on the specifics of the space $\P(\R^d)$ and the Fokker-Planck equation \eqref{FPE} to define suitable test functions on $\P(\R^d)$, allowing to deal with solutions in $C([0,T];\P(\P(\R^d)))$ and to reduce Theorem \ref{Thm_PDEtoSPDE}  to the superposition principle  in $\R^\infty$ of \cite[Theorem 7.1]{trevisan_diss}.

\medskip

Combining Theorems \ref{Thm_PDEtoSPDE} and \ref{Thm_SPDEtoSDE_general} shows that to any continuous solution $P=(P_t)_{t \in [0,T]}$ of the Fokker-Planck equation \eqref{FPE} one can associate continuous processes $\mu=(\mu_t)_{t \in [0,T]}$ and $\bm{X}=(\bm{X}_t)_{t \in [0,T]}$, with values in $\P(\R^d)$ and $\R^d$, respectively, and solving the SPDE \eqref{SPDE_special} and the McKean-Vlasov SDE \eqref{McV_special}, respectively, with the time-marginal relations $P_t = \L(\mu_t)$ and $\mu_t=\L(\bm{X}_t\,|\,\G_t)$. In a sense, $P$ and $\mu$ represent successive marginalizations of the process $\bm{X}$. For another perspective on this hierarchy, one may view $P$, $\L(\mu)$, and $\L(\L(\bm{X}\,|\,\G_T))$ as elements of $C([0,T];\P(\P(\R^d)))$, $\P(C([0,T];\P(\R^d)))$, and $\P(\P(C([0,T];\R^d)))$, respectively. The connections between $P$, $\mu$, and $\bm{X}$ are summarized in the following diagram:
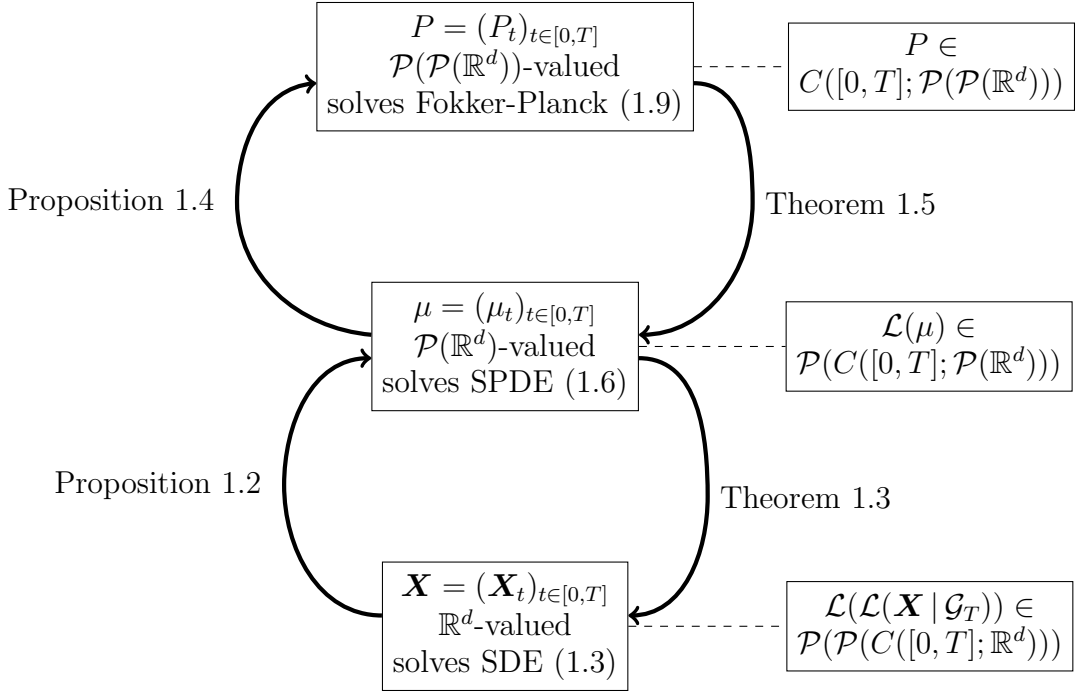
\begin{tikzpicture}[auto,node distance=2cm]
  \node[align=center,rectangle,draw=black] (level1) {$\bm{X}=(\bm{X}_t)_{t \in [0,T]}$ \\ $\R^d$-valued \\ solves SDE \eqref{McV_special} };
  \node[align=center,rectangle,draw=black] (level2) [above = of level1] {$\mu=(\mu_t)_{t \in [0,T]}$ \\  $\P(\R^d)$-valued \\ solves SPDE \eqref{SPDE_special} };
  \node[align=center,rectangle,draw=black] (level3) [above = of level2] {$P=(P_t)_{t \in [0,T]}$ \\ $\P(\P(\R^d))$-valued \\ solves Fokker-Planck \eqref{FPE} };
  \draw[black,ultra thick,<-]
  	(level2.185) to[out=180,in=90] ++(-33pt,-48pt) node[left,align=center] { Proposition \ref{pr:SDEtoSPDE} } to[out=270,in=180] (level1.175);
  \draw[black,ultra thick,<-]
  	(level3.185) to[out=185,in=90] ++(-30pt,-45pt) node[left,align=center] { Proposition \ref{pr:SPDEtoPDE}  } to[out=270,in=175] (level2.175);
  \draw[black,ultra thick,<-]
  	(level1.5) to[out=0,in=270] ++(+30pt,+45pt) node[right,align=center] {Theorem \ref{Thm_SPDEtoSDE_general} } to[out=90,in=0] (level2.-5);
  \draw[black,ultra thick,<-]
  	(level2.5) to[out=0,in=270] ++(+43pt,+50pt) node[right,align=center] {Theorem \ref{Thm_PDEtoSPDE}} to[out=90,in=0] (level3.-5);
  \node[align=center,rectangle,draw=black] (level1note) [right = 2.1cm of level1] {$\L(\L(\bm X\,|\,\G_T)) \in $ \\ $\P(\P(C([0,T];\R^d)))$};
  \node[align=center,rectangle,draw=black] (level2note) [right =1.95cm of level2] {$\L(\mu) \in $ \\ $\P(C([0,T];\P(\R^d)))$};
  \node[align=center,rectangle,draw=black] (level3note) [right =1.25cm of level3] {$P \in $ \\ $C([0,T];\P(\P(\R^d)))$};
  \draw[black,dashed] (level1) -- (level1note);
  \draw[black,dashed] (level2) -- (level2note);
  \draw[black,dashed] (level3) -- (level3note);
\end{tikzpicture}

\subsection{Mimicking theorem}
In addition to being of interest in their own right, Theorems \ref{Thm_SPDEtoSDE_general} and \ref{Thm_PDEtoSPDE} can be used to obtain the following ``mimicking'' theorem.

\begin{corollary}\label{thm_mim}
Suppose $(\Omega,\FF,\PP)$ is a filtered probability space supporting $d$-dimensional $\FF$-Brownian motions $\bm{B}$ and $\bm{W}$, as well as a subfiltration $\GG$ with respect to which $\bm{B}$ is adapted and $\bm{W}$ is independent. Let $(\bm{b},\bm{\sigma})$ be an $\FF$-progressive process with values in $\R^d\times\R^{d\times d}$, and $\widehat{\bm \gamma}:\,[0,T] \times \P(\R^d) \times \R^d \to \R^{d \times d}$ be measurable.
Suppose $(\bm{X}_t,\mu_t)_{t\in [0,T]}$ is a continuous $\FF$-adapted process with values in $\R^d \times \P(\R^d)$ and satisfying: 
\begin{enumerate}[(a)]
\item $\E\big[\int_{0}^{T} |\bm b_t|^p + |\bm \sigma_t\bm \sigma_t^\top + \widehat{\bm \gamma}\widehat{\bm \gamma}^\top(t,\mu_t,\bm X_t)|^p\, \d t\big]<\infty$ for some $p>1$. 
\item $\F_t \indep \F^{\bm{W}}_T \vee \G_T\,|\,\F^{\bm{W}}_t \vee \G_t$ for each $t \in [0,T]$.
\item The following SDE holds:
\begin{align}
\begin{split}
\d\bm {\bm{X}}_t &=\bm b_t \d t+\bm\sigma_t\d\bm W_t+ \widehat{\bm \gamma}(t,\mu_t,\bm X_t)\d\bm B_t, \\
\mu_t &=\L(\bm {\bm{X}}_t\,|\,\G_T), \quad t\in[0,T].
\end{split} \label{SDE-mimicking}
\end{align}
\end{enumerate}
Then $\mu_t=\L(\bm{X}_t\,|\,\G_t)$ a.s., for each $t \in [0,T]$, and there exist measurable functions $(\widehat{\bm b},\widehat{\bm \sigma}):\,[0,T]\times\P(\R^d)\times\R^d\to \R^d\times\R^{d\times d}$ such that
\begin{align*}
\widehat{\bm b}(t,\mu_t,\bm{X}_t)= \E[\bm b_t \,|\,\mu_t,\,\bm{X}_t], \ \ \widehat{\bm\sigma}\widehat{\bm\sigma}^\top (t,\mu_t,\bm{X}_t)= \E[\bm\sigma_t\bm\sigma_t^\top\,|\,\mu_t,\,\bm{X}_t]\;\;\text{a.s., for a.e.~} t \in [0,T].
\end{align*}
Finally, for any such functions $(\widehat{\bm b},\widehat{\bm \sigma})$, there exists a filtered probability space $(\widehat\Omega,\widehat{\mathbb F},\widehat{\PP})$ supporting a subfiltration $\widehat{\mathbb G}\subset\widehat{\mathbb F}$, two $d$-dimensional $\widehat\FF$-Brownian motions $\widehat{\bm W}$ and $\widehat{\bm B}$ with $\widehat{\bm B}$ adapted to $\widehat\GG$ and $\widehat{\bm{W}}$ independent of $\widehat\G_T$, and a continuous $\widehat\FF$-adapted $\R^d \times \P(\R^d)$-valued process $(\widehat{\bm X}_t,\widehat\mu_t)_{t \in [0,T]}$ such that:
\begin{enumerate}[(1)]
\item The McKean-Vlasov SDE holds:
\begin{align*}
\d\widehat{\bm X}_t &= \widehat{\bm b}(t,\widehat{\mu}_t,\widehat{\bm X}_t)\d t+\widehat{\bm\sigma}(t,\widehat{\mu}_t,\widehat{\bm X}_t)\d\widehat{\bm W}_t+\widehat{\bm\gamma}(t,\widehat{\mu}_t,\widehat{\bm X}_t)\d\widehat{\bm B}_t, \\
\widehat{\mu}_t &=\L\big(\widehat{\bm X}_t \,|\, \widehat{\G}_T\big)=\L\big(\widehat{\bm X}_t \,|\, \widehat{\G}_t\big)\;\;\text{a.s.}, \quad t\in[0,T].
\end{align*}
\item $(\widehat{\bm X}_t,\widehat \mu_t)\overset{d}{=}(\bm {\bm{X}}_t,\mu_t)$, for each $t\in[0,T]$.
\item $\F^{\widehat{\bm X}}_t \indep \F^{\widehat{\bm W}}_T \vee \G_T\,|\,\F^{\widehat{\bm W}}_t \vee \G_t$, for each $t \in [0,T]$.
\end{enumerate}
\end{corollary} 

\smallskip

The proof is given in Section \ref{sec_application}. The idea is that an It\^o process of the form \eqref{SDE-mimicking} can be ``mimicked," in terms of the time-marginals $(\L(\bm X_t,\mu_t))_{t \in [0,T]}$, by a process solving a \emph{Markovian} conditional McKean-Vlasov equation. We allow extra randomness in only the $(\bm b,\bm\sigma)$ coefficients of \eqref{SDE-mimicking}, whereas the $\widehat{\bm \gamma}$ term is assumed to already be in Markovian form; one can generalize this easily to allow a process $\bm\gamma$ satisfying $\bm\gamma_t \indep \G_t\,|\,(\bm X_t,\mu_t)$ for all $t \in [0,T]$.

\medskip

Theorems on mimicking aspects (such as the time-marginal distributions) of an It\^o process by those of a diffusion process go back to the seminal papers \cite{Kry,Gyo} (see also \cite{Kle} for an independent development). A very general result on mimicking aspects of It\^o processes, including the distributions of their marginals, running maxima and running integrals, is given in \cite[Theorem 3.6]{BS} (see also \cite[Theorem 1.13]{FP} for a mimicking theorem in another degenerate setting). The analogue of this question in the context of optimal control of Markov processes has been studied in \cite{KuSt}. Related results have appeared in \cite{Ke,MY}, where Markov (sub)martingales with given fixed-time distributions have been constructed. Corollary \ref{thm_mim} differs from the previous mimicking theorems in that both the fixed-time distributions and the fixed-time conditional distributions of an It\^o process are being mimicked. 

\medskip

The idea of using a superposition principle to prove a mimicking theorem seems to be new, and to illustrate the idea we record here a simple proof of the classical mimicking theorem \cite{Gyo,BS} using the superposition principle of Trevisan \cite{trevisan} (or Figalli \cite{figalli2008existence} if the coefficients are bounded). Indeed, suppose a filtered probability space $(\Omega,\FF,\PP)$ supports an $\FF$-Brownian motion $\bm{W}$ and an $\FF$-progressive process $(\bm{b}_t,\bm{\sigma}_t)_{t \in [0,T]}$ of suitable dimension satisfying $\E\big[\int_0^T|\bm{b}_t| + |\bm{\sigma}_t\bm{\sigma}_t^\top|\d t\big] < \infty$. Consider a $d$-dimensional It\^o process $\bm{X}$ with
\begin{align*}
\d\bm{X}_t = \bm{b}_t \d t + \bm{\sigma}_t \d\bm{W}_t.
\end{align*}
Apply It\^o's formula to a test function $\varphi \in C_c^\infty(\R^d)$ and take expectations to find
\begin{align*}
\E[\varphi(\bm{X}_t)] &= \E[\varphi(\bm{X}_0)] + \E\bigg[\int_0^t \bm{b}_s \cdot \nabla \varphi(\bm{X}_s) + \tfrac12\bm{\sigma}_s\bm{\sigma}_s^\top : \nabla^2 \varphi(\bm{X}_s)\d s\bigg].
\end{align*}
If we define  $\widehat{\bm{b}}(t,x) = \E[\bm b_t \, | \,\bm{X}_t=x]$ and $\widehat{\bm{\sigma}}(t,x) = \E[\bm{\sigma}_t\bm{\sigma}_t^\top \, | \,\bm{X}_t=x]^{1/2}$, then Fubini's theorem and the tower property in the above equation yield
\begin{align*}
\E[\varphi(\bm{X}_t)] &= \E[\varphi(\bm{X}_0)] + \E\bigg[\int_0^t \widehat{\bm{b}}(s,\bm{X}_s) \cdot \nabla \varphi(\bm{X}_s) + \tfrac12 \widehat{\bm{\sigma}}\widehat{\bm{\sigma}}^\top(s,\bm{X}_s) : \nabla^2 \varphi(\bm{X}_s)\d s\bigg].
\end{align*}
This shows that the marginal flow $(\mu_t=\L(\bm{X}_t))_{t \in [0,T]}$ is a solution of the Fokker-Planck equation
\begin{align}
\langle \mu_t-\mu_0,\varphi\rangle = \int_0^t \big\langle\mu_s, \widehat{\bm{b}}(s,\cdot) \cdot \nabla \varphi + \tfrac12  \widehat{\bm{\sigma}}\widehat{\bm{\sigma}}^\top(s,\cdot) : \nabla^2\varphi\big\rangle \d s, \ \ t \in [0,T], \ \ \varphi \in C^\infty_c(\R^d), \label{intro:FPE2}
\end{align}
which is associated with the SDE
\begin{align}
\d\bm{\widehat{X}}_t &= \widehat{\bm{b}}(t,\bm{\widehat{X}}_t)\d t + \widehat{\bm{\sigma}}(t,\bm{\widehat{X}}_t)\d\bm{\widehat{\bm W}}_t. \label{intro:mimickedSDE}
\end{align}
According to the superposition principle of Trevisan \cite[Theorem 2.5]{trevisan}, since $(\mu_t)_{t \in [0,T]}$ solves the Fokker-Planck equation, there must exist a weak solution of the SDE \eqref{intro:mimickedSDE} which shares the same marginal flow, i.e.,  $\L(\bm{\widehat{X}}_t)=\mu_t=\L(\bm{X}_t)$ for all $t \in [0,T]$. This recovers precisely the classical mimicking theorem in the form of \cite[Corollary 3.7]{BS}. Our method of proving Corollary \ref{thm_mim} is completely analogous, using a combination of Theorems \ref{Thm_PDEtoSPDE} and \ref{Thm_SPDEtoSDE_general} in place of \cite[Theorem 2.5]{trevisan}, and the Fokker-Planck equation \eqref{FPE} on $\P(\R^d)$  in place of the Fokker-Planck equation \eqref{intro:FPE2} on $\R^d$.

\subsection{Mean field games and control}

Our specific setup in Corollary \ref{thm_mim} is motivated in part by questions from the theory of mean field games and controlled McKean-Vlasov dynamics, also known as mean field control.
In particular, we show in Section \ref{MFG_appl} how to use Corollary \ref{thm_mim} to convert an optimal control into a Markovian one, in the sense that the control is a function of $(t,\bm{X}_t,\mu_t)$ only, in the context of controlled McKean-Vlasov dynamics with common noise.
In the case that $(\mu_t)_{t \in [0,T]}$ is non-random, this has been done using the classical mimicking theorem, but the case of stochastic $(\mu_t)_{t \in [0,T]}$, as in any model with common noise, seems to be out of reach of prior methods.
See Section \ref{MFG_appl} for more details.

\subsection{Further remarks on the Fokker-Planck equation on $\P(\R^d)$} \label{se:FPE}

The Fokker-Planck equation \eqref{FPE} is stated in terms of test functions $F :\P(\R^d) \to \R$ of the form
\begin{align}
F(m) = f\big(\langle m, \varphi_1\rangle,\ldots,\langle m, \varphi_k\rangle\big), \label{test functions-simple}
\end{align}
where $k \in \N$, $f \in C^\infty_c(\R^k)$, and $\varphi_1,\ldots,\varphi_k \in C^\infty_c(\R^d)$. Alternatively, we can write the equation in terms of a larger class of test functions, specified through differentiation for real-valued functions on $\P(\R^d)$, known in some recent literature as the \emph{L-derivative}. See \cite[Chapter I.5]{CD2} for a careful and thorough development of this notion of derivative. The test functions of interest for us are the following: 

\begin{definition} \label{def:Lions}
Define $C^2_b(\P(\R^d))$ as the set of bounded continuous functions $F : \P(\R^d) \to \R$ such that:
\begin{itemize}
\item There exists a bounded continuous function $\partial_m F : \P(\R^d) \times \R^d \to \R$ satisfying
\begin{align*}
\lim_{h \downarrow 0}\frac{F(m + h(m'-m)) - F(m)}{h} &= \int_{\R^d}\partial_m F(m,v)\,(m'-m)(\!\d v),
\end{align*}
for all $m,m' \in \P(\R^d)$. Given that such a function $\partial_m F$ exists, it is unique provided we impose the additional requirement that
\begin{align*}
\int_{\R^d}\partial_m F(m,v) \, m(\!\d v) = 0, \quad m \in \P(\R^d).
\end{align*}
\item The function $v \mapsto \partial_m F(m,v)$ is continuously differentiable, and its gradient, denoted $D_m F(m,v)$, is uniformly bounded in $(m,v)$.
\item For every fixed $v' \in \R^d$, each component of the $\R^d$-valued function $m \mapsto D_mF(m,v')$ satisfies the properties described in the first two bullet points, resulting in some $D_m^2F(m,v',v) \in \R^{d \times d}$. The function $D_m^2F$ is bounded and continuous.
\item For $m \in \P(\R^d)$, we write $D_vD_mF(m,v)$ for the Jacobian of the function $v \mapsto D_mF(m,v)$ and assume it to be continuous and bounded in $(m,v)$.
\end{itemize}
\end{definition}

An equivalent form of the Fokker-Planck equation \eqref{FPE}, under the integrability assumption \eqref{integ2}, is then 
\begin{align}
\langle P_t-P_0,F\rangle = \int_0^t\langle P_s, \M_s F\rangle \d s, \quad t \in [0,T], \quad F \in C^2_b(\P(\R^d)), \label{FPE-Lionsderivative}
\end{align}
where we define, for $t \in [0,T]$ and $m \in \P(\R^d)$,
\begin{align*}
\M_t F(m) := \ &\int_{\R^d} \Big[D_m F(m,v) \cdot \bm{b}(t,m,v) + \frac12 D_vD_mF(m,v) : \bm{a}(t,m,v) \Big]  m(\!\d v) \\
	&+ \frac12\int_{\R^d}\int_{R^d}D_m^2F(m,v,v') : [\bm{\gamma}(t,m,v)\bm{\gamma}^\top(t,m,v')]\,m(\!\d v)\,m(\!\d v').
\end{align*}
One direction of this equivalence is easy: If $F$ is of the form \eqref{test functions-simple}, then straightforward calculus shows that $F \in C^2_b(\P(\R^d))$ and
\begin{align*}
D_mF(m,v) &=  \sum_{i=1}^k \partial_i f\big(\langle m, \varphi_1\rangle,\ldots,\langle m, \varphi_k\rangle\big)\nabla \varphi_i(v), \\
D_m^2F(m,v,v') &=  \sum_{i,j=1}^k \partial_{ij} f\big(\langle m, \varphi_1\rangle,\ldots,\langle m, \varphi_k\rangle\big)\nabla \varphi_i(v)\nabla \varphi_j(v')^\top, \\
D_vD_mF(m,v) &=  \sum_{i=1}^k \partial_i f\big(\langle m, \varphi_1\rangle,\ldots,\langle m, \varphi_k\rangle\big)\nabla^2 \varphi_i(v).
\end{align*}
It is then easy to check that for these test functions the equation \eqref{FPE-Lionsderivative} rewrites exactly as \eqref{FPE}. On the other hand, suppose we start from some $P=(P_t)_{t\in [0,T]} \in C([0,T];\P(\P(\R^d)))$ satisfying \eqref{FPE}. Applying Theorems \ref{Thm_PDEtoSPDE} and \ref{Thm_SPDEtoSDE_general}, we find processes $\mu$ and $\bm X$ following \eqref{SPDE_special} and \eqref{McV_special}, respectively, along with $\mu_t=\L({\bm{X}}_t\,|\,\G_t)$ and $P_t=\L(\mu_t)$.
For $F \in C^2_b(\P(\R^d))$, we may then apply a recent form of It\^o's formula for conditional measure flows of It\^o processes, see \cite[Theorem II.4.14]{CD2}. This yields 
\begin{align*}
F(\mu_t) &= F(\mu_0) +  \int_0^t\M_s F(\mu_s) \d s + R_t,
\end{align*}
for a certain mean-zero martingale $R$ which can be identified explicitly, but we have no need to do so. Taking expectations and using $P_t=\L(\mu_t)$ gives the Fokker-Planck equation in the form of \eqref{FPE-Lionsderivative}.

\medskip

\subsection{Existence and uniqueness}
One obvious implication of our superposition principles, Theorems \ref{Thm_PDEtoSPDE} and \ref{Thm_SPDEtoSDE_general}, is that the following are equivalent:
\begin{itemize}
\item There exists a solution of the Fokker-Planck equation \eqref{FPE}, or its equivalent form \eqref{FPE-Lionsderivative}.
\item There exists a solution of the SPDE \eqref{SPDE_special}.
\item There exists a solution of the conditional McKean-Vlasov SDE \eqref{McV_special}.
\end{itemize}
Uniqueness is somewhat more subtle. Using Theorems \ref{Thm_PDEtoSPDE} and \ref{Thm_SPDEtoSDE_general}, we deduce that uniqueness for the conditional McKean-Vlasov SDE \eqref{McV_special}\footnote{Here and more generally for SDEs with random coefficients such as \eqref{SDE_gen}, the appropriate notion of uniqueness in law is in the sense of the \emph{very good solution measures} of \cite{jacod1981weak}, i.e., the uniqueness of the measure induced on $\Omega \times C([0,T];\R^d)$ by appending the solution process $\bm X$.} implies uniqueness for the SPDE \eqref{SPDE_special}, and uniqueness for the SPDE \eqref{SPDE_special} implies uniqueness for the Fokker-Planck equation \eqref{FPE}. 
The reverse implications are likely true but would require much more machinery to prove; one must show, for example, that the law of a solution $(\bm X_t)_{t \in [0,T]}$ of the SDE \eqref{McV_special} is uniquely determined by its \emph{marginal} conditional laws $\L(\bm X_t\,|\,\G_t)$. In the unconditional case studied by Figalli \cite{figalli2008existence} and Trevisan \cite{trevisan}, this relies on a result of Stroock and Varadhan \cite[Theorem 6.2.3]{stroock-varadhan}, which shows a martingale problem is determined by its marginals in a certain carefully specified sense.

\medskip

Under various assumptions (e.g., Lipschitz coefficients), existence and uniqueness results are known for McKean-Vlasov SDEs like \eqref{McV_special} (see \cite{kurtz1999particle,HSS} or \cite[Section II.2.1.3]{CD2}) and for the SPDE \eqref{SPDE_special} (see \cite{kurtz1999particle,DV}).  From such results we may, rather remarkably, deduce corresponding existence and uniqueness results for the Fokker-Planck equation on $\P(\R^d)$ given by \eqref{FPE} or \eqref{FPE-Lionsderivative}, which appear to be the first of their kind.

\medskip

Moreover, our Theorem \ref{Thm_PDEtoSPDE} allows us to recover and extend prior results on well-posedness for stochastic Fokker-Planck equations. For instance, if we assume the coefficients $(\bm b,\bm \sigma,\bm \gamma)$ from \eqref{SDE_gen} are Lipschitz in $\bm x$ uniformly in $(t,\omega)$, and $\E\big[\int_0^T |\bm b(t,\omega,0)|^2 + |\bm a(t,\omega,0)|\,\mathrm{d}t\big] < \infty$, then standard arguments yield existence and uniqueness for the SDE \eqref{SDE_gen}. Using Theorem \ref{Thm_PDEtoSPDE} we immediately deduce existence and uniqueness for the stochastic Fokker-Planck equation \eqref{SPDE}. This quickly recovers and extends results of \cite{kurtz1999particle,coghi2019stochastic} on linear SPDEs of Fokker-Planck type. Likewise, the conditional McKean-Vlasov equation \eqref{McV_special} is well-posed under Lipschitz assumptions, using a Wasserstein metric for the measure argument; see \cite[Propositions II.2.8 and II.2.11]{CD2}, \cite[Theorem 3.3]{HSS}, \cite[Theorem 3.1]{kurtz1999particle}, or \cite[Theorem 3.3]{coghi2019stochastic}. 
From well-posedness of \eqref{McV_special} and our superposition principle, Theorem \ref{Thm_PDEtoSPDE}, we immediately deduce well-posedness for the corresponding nonlinear SPDE \eqref{SPDE_special}, which again recovers results of \cite{kurtz1999particle,coghi2019stochastic}.

\medskip

We also stress that our results apply to many cases of singular or local interactions, thanks to the minimal regularity required. For instance, suppose the SPDE \eqref{SPDE_special} is \emph{local} in the sense that $\bm b(t,m,x) = \widetilde{\bm b}(t,m(x),x)$ for some $\widetilde{\bm b} : [0,T] \times \R_+ \times \R^d \to \R$, whenever $m \in \P(\R^d)$ admits a density with respect to Lebesgue measure (denoted $m(x)$), and suppose $\bm a$ and $\bm \gamma$ take similar forms. Indeed, to apply Theorem \ref{Thm_PDEtoSPDE}, we may extend the domain of the coefficients in an arbitrary (measurable, bounded) fashion to include all of $\P(\R^d)$, not just those measures which admit a Lebesgue density. Hence, if there exists a solution of the (local, in this case) SPDE \eqref{SPDE_special} with $\mu_t(\omega)$ absolutely continuous with respect to Lebesgue measure for a.e.\ $(t,\omega)$ and satisfying the requisite integrability conditions, then there exists a solution of the corresponding (density-dependent) McKean-Vlasov SDE \eqref{McV_special}. Similarly, if there is a unique solution of the (density-dependent) SDE \eqref{McV_special}, we again deduce existence and uniqueness for the (local) SPDE.

\medskip

The Fokker-Planck equation on $\P(\R^d)$ given in \eqref{FPE-Lionsderivative} might be compared to certain \emph{backward} PDEs on Wasserstein space studied in recent literature on McKean-Vlasov equations and mean field games. The equation \eqref{FPE-Lionsderivative} is a Fokker-Planck or Kolmogorov forward equation associated with the process $(\mu_t)_{t \in [0,T]}$ of \eqref{SPDE_special}, and we are not aware of any prior studies of this equation in the literature. On the other hand, there have been a number of recent studies of the Kolmogorov \emph{backward} equation (i.e., Feynman-Kac formulas) associated with $(\bm{X}_t,\mu_t)_{t \in [0,T]}$ from \eqref{McV_special}, \eqref{SPDE_special}. That is, the functions $(t,x,m) \mapsto \E[F(\bm{X}_T,\mu_T) \, | \, \bm{X}_t=x,\mu_t=m]$ should satisfy a certain linear PDE on $[0,T] \times \R^d \times \P(\R^d)$, provided $F$ is sufficiently nice \cite{CrMc,RaFr1}, and regularity estimates on the solution are useful for quantitative propagation of chaos arguments \cite{RaFr2,CST}. Nonlinear analogues of these backward PDEs appear in the guise of the \emph{master equation} from the theory of mean field games \cite{CCD,CDLL} and in mean field control problems \cite{PhWe,CD2}.

\begin{remark}
It is worth noting that we always interpret the SPDE \eqref{SPDE} \emph{in the weak PDE sense}, i.e., in the sense of distributions.
In general, we are also interpreting this SPDE \emph{in the weak probabilistic sense}, with $\mu$ not necessarily adapted to the filtration generated by the driving Brownian motion $\bm{B}$. If $\mu$ is adapted to $\FF^{\bm{B}}$ (i.e., if $\GG=\FF^{\bm{B}}$), then we might say we have a (probabilistic) strong solution of the SPDE \eqref{SPDE}. We will not dwell on this point, because whether we have a strong or weak solution is irrelevant to our purposes.
\end{remark}

\subsection{Outline of the paper}

The rest of the paper is structured as follows.
We begin with a short discussion of compatibility conditions in Section \ref{se:compatibility}.
Section \ref{sec_pre} proves Theorem \ref{Thm_SPDEtoSDE_general} first under additional smoothness restrictions on the coefficients. In this regime, the claimed superposition follows from the well-posedness (which we establish) of the SDE \eqref{SDE_gen} and SPDE \eqref{SPDE}. Before turning to the general case, Section \ref{se:stableconvergence} then develops some preliminary results on tightness and stable convergence which will aid in our successive approximations in the subsequent Section \ref{sec_proof_SPDEtoSDE}. Section \ref{sec_proof_SPDEtoSDE} contains the main line of the proof of Theorem \ref{Thm_SPDEtoSDE_general}, which follows a sequence of reductions to the aforementioned smooth case.
Next, Section \ref{sec_proof_PDEtoSPDE} is devoted to the (surprisingly short) proof of Theorem \ref{Thm_PDEtoSPDE}.
Our mimicking result, Corollary \ref{thm_mim}, is derived in Section \ref{sec_application} from Theorems \ref{Thm_SPDEtoSDE_general} and \ref{Thm_PDEtoSPDE}.
The final Section \ref{MFG_appl} details applications to controlled McKean-Vlasov dynamics.
Lastly, three short appendix sections give less central proofs omitted from the body of the paper. 

Notably, Section \ref{sec_proof_PDEtoSPDE}, Section \ref{sec_application}, and Section \ref{MFG_appl} are independent of each other and also of Sections \ref{se:compatibility}--\ref{sec_proof_SPDEtoSDE}. That is, after this introduction, one could read either Sections \ref{se:compatibility}--\ref{sec_proof_SPDEtoSDE}, Section \ref{sec_proof_PDEtoSPDE}, Section \ref{sec_application}, or Section \ref{MFG_appl}, without loss of continuity.


\section{Compatibility preliminaries} \label{se:compatibility}

This short section collects a few elementary implications related to the recurring condition that $\F^{\bm{X}}_t \indep \F^{\bm{W}}_T \vee \G_T\,|\,\F^{\bm{W}}_t \vee \G_t$ for all $t \in [0,T]$.
We continue here with our notational conventions for $\sigma$-algebras and filtrations.
Given a stochastic process $\bm{X}=(\bm{X}_t)_{t \in [0,T]}$,  we write $\FF^{\bm{X}}=(\F^{\bm{X}}_t)_{t \in [0,T]}$ for the filtration it generates. For two filtrations $\GG =(\G_t)_{t \in [0,T]}$ and $\HH=(\H_t)_{t \in [0,T]}$, we let $\GG \vee \HH := (\G_t \vee \H_t)_{t \in [0,T]}$.
We should stress that, as usual, a process $\bm{W}$ is said to be a \emph{Brownian motion with respect to a filtration $\HH=(\H_t)_{t \in [0,T]}$} if $\bm{W}$ is a Brownian motion which is adapted to $\HH$ and also has independent increments with respect to $\HH$, meaning $\bm{W}_t-\bm{W}_s \indep \H_s$ for $0 \le s < t \le T$.
The following is proven in Appendix \ref{se:misc}.

\begin{lemma} \label{le:compatibility}
Suppose $(\Omega,\FF=(\F_t)_{t \in [0,T]},\PP)$ is a filtered probability space supporting a $d$-dimensional $\FF$-Brownian motion $\bm{W}$ as well as two subfiltrations $\GG =(\G_t)_{t \in [0,T]}$ and $\HH=(\H_t)_{t \in [0,T]}$.
Assume $\bm{W}$ is independent of $\G_T$.
Then the following four statements are equivalent:
\begin{enumerate}[(1)]
\item $\H_t \indep \F^{\bm{W}}_T \vee \G_T \,|\, \F^{\bm{W}}_t \vee \G_t$, for each $t \in [0,T]$. 
\item The following three conditions hold:
\begin{enumerate}[(2a)]
\item $\H_t \indep \H_0 \vee \F^{\bm{W}}_T \vee \G_T \,|\, \H_0 \vee \F^{\bm{W}}_t \vee \G_t$, for each $t \in [0,T]$. 
\item $\H_t \indep \G_T \,|\, \G_t$, for each $t \in [0,T]$.
\item $\bm{W} \indep \H_0 \vee \G_T$.
\end{enumerate}
\item The three $\sigma$-algebras $\H_t$, $\G_T$, and $\F^{\bm{W}}_T$ are conditionally independent given $\F^{\bm{W}}_t \vee \G_t$, for each $t \in [0,T]$.
\item The following two conditions hold:
\begin{enumerate}[(4a)]
\item $\bm{W}$ is a Brownian motion with respect to the filtration $(\G_T \vee \F^{\bm{W}}_t \vee \H_t)_{t \in [0,T]}$. 
\item $\H_t \vee \F^{\bm{W}}_t \indep \G_T \,|\, \G_t$ for all $t \in [0,T]$.
\end{enumerate}
\item $\H_t \vee \F^{\bm W}_T \indep \G_T  \,|\,  \G_t$, for each $t \in [0,T]$.
\end{enumerate}
\end{lemma}

We will most often apply Lemma \ref{le:compatibility} with $\HH=\FF^{\bm{X}}$ for a given process $\bm{X}$. Hence, in Proposition \ref{pr:SDEtoSPDE}, the assumption $\F^{\bm{X}}_t \indep \F^{\bm{W}}_T \vee \G_T \,|\, \F^{\bm{W}}_t \vee \G_t$, for each $t \in [0,T]$, implies that $\F^{\bm{X}}_t \indep \G_T \,|\,\G_t$, for each $t \in [0,T]$, which in turn yields $\L(\bm{X}_t\,|\,\G_T)=\L(\bm{X}_t\,|\,\G_t)$ a.s.

\medskip

Despite possibly appearing tangential at first glance, carefully formulated compatibility conditions are essential for dealing with weak solutions for SDEs \cite{jacod1981weak}, McKean-Vlasov equations \cite{HSS,LaSv}, control problems \cite{DjPoTa,DjPoTa2}, and mean field games \cite{CDL,CD2}.
For this paper, we favor the most concise formulation $\F^{\bm X}_t \indep \F^{\bm W}_T \vee \G_T \,|\, \F^{\bm W}_t \vee \G_t$ for all $t \in [0,T]$.


\section{Superposition from SPDE to SDE: smooth case}\label{sec_pre}

In this section, we establish Theorem \ref{Thm_SPDEtoSDE_general} in the case where the coefficients ${\bm b},\,{\bm a},\,{\bm \gamma}$ are sufficiently smooth in the ${\bm x}$-variable. In this case we are able to construct a \textit{strong} solution of the SDE \eqref{SDE_gen}, which is not feasible in general. 

\medskip

We use the following notation for norms throughout the paper. We let $|\bm x|$ denote the Euclidean norm of a vector $\bm x \in \R^d$ and $|\bm A|$ the Frobenius norm of a matrix $\bm A$. For a real-valued function $f : \R^d \to \R$, we set
\begin{align*}
&\|f\|_{C_b(\R^d)} = \sup_{\bm x \in \R^d}|f(\bm x)|, \quad \|f\|_{C^1_b(\R^d)} = \|f\|_{C_b(\R^d)} + \sup_{\bm x \in \R^d}|\nabla f(\bm x)|, \\
&\|f\|_{C^2_b(\R^d)} = \|f\|_{C^1_b(\R^d)} + \sup_{\bm x \in \R^d}|\nabla^2 f(\bm x)|.
\end{align*}
Set $\|f\|_{C^1_b(\R^d)} = \infty$ if $f$ is not continuously differentiable, and similarly for $\|f\|_{C^2_b(\R^d)}$.
We extend these norms to vector-valued functions in the natural way; for instance, for $f=(f_1,\ldots,f_k) : \R^d \to \R^k$, we take $\|f\|_{C_b(\R^d)}^2 := \|f_1\|_{C_b(\R^d)}^2 + \cdots +\|f_k\|_{C_b(\R^d)}^2$. Lastly, for a Lipschitz function $f:\,\R^d \to \R^k$, we denote by $\|f\|_{\mathrm{Lip}(\R^d)}$ its Lipschitz constant, i.e., the smallest constant $c \ge 0$ such that $|f(\bm x)-f(\bm y)| \le c|\bm x-\bm y|$ for all $\bm x,\bm y \in \R^d$. If $f$ is not Lipschitz, then $\|f\|_{\mathrm{Lip}(\R^d)} := \infty$.

\begin{proposition}[Theorem \ref{Thm_SPDEtoSDE_general}, smooth case] \label{Thm_SPDEtoSDE_unismooth}
In the setting of Theorem \ref{Thm_SPDEtoSDE_general}, suppose that, for a.e.~$(t,\omega)\in[0,T]\times\Omega$, the functions $\bm b(t,\omega,\cdot),\,\bm a(t,\omega,\cdot),\, \bm \gamma(t,\omega,\cdot)$ are bounded and twice continuously differentiable with bounded derivatives of the first and second order. Then there exists a positive constant $D_{p} <\infty$ such that the following holds: If
\begin{equation}\label{condition_uniformsmooth} 
\max_{i=1,\ldots,N}\int_{\frac{(i-1)T}{N}}^{\frac{iT}{N}} \! \big\|\bm b(t,\omega,\cdot)\big\|^{p}_{C_b^1(\R^d)}
\!+\!\big\|\bm a(t,\omega,\cdot)\big\|^{p}_{C_b^2(\R^d)}
\!+\!\big\|\bm \gamma(t,\omega,\cdot)\big\|^{2p}_{C_b^2(\R^d)}\d t\le D_{p}\;\;\text{a.s.}
\end{equation}
for some $N \in \N$, and $\mu_0$ admits a density $\rho\in L_2(\R^d)$ with $\E[\|\rho\|_2^2]<\infty$, then the conclusion of Theorem \ref{Thm_SPDEtoSDE_general} applies. 
\end{proposition}

\noindent\textbf{Proof.}
Let $\|\bm\sigma(t,\omega,\cdot)\|_{\mathrm{Lip}(\R^d)}$ denote the Lipschitz constant of $x \mapsto \bm\sigma(t,\omega,x)$.
Since $\bm\sigma^2=\bm a -\bm\gamma\bm\gamma^\top$, we have by \cite[Theorem 5.2.3]{stroock-varadhan}:
\begin{equation*}
\begin{split}
\int_0^T\big\|\bm\sigma(t,\omega,\cdot)\big\|^{2 p }_{\mathrm{Lip}(\R^d)}\d t
& \leq\int_0^T(2d^2)^{ p }\,\big\|\bm\sigma(t,\omega,\cdot)^2\big\|^{ p }_{C_b^2(\R^d)}\d t \\
&\leq (2d^2)^{ p }\,\int_0^T  \Big(\big\|\bm a(t,\omega,\cdot)\big\|_{C_b^2(\R^d)}
+ 2\big\|\bm \gamma(t,\omega,\cdot)\big\|^{2}_{C_b^2(\R^d)}\Big)^{ p }\d t \\
&\leq(8d^2)^{ p }\,\int_0^T \big\|\bm a(t,\omega,\cdot)\big\|^{ p }_{C_b^2(\R^d)} + \big\|\bm \gamma(t,\omega,\cdot)\big\|^{2 p }_{C_b^2(\R^d)}\d t \\
&\leq(8d^2)^{ p } ND_{ p }\;\;\text{a.s.}
\end{split}
\end{equation*}
Our construction starts with an extension of the probability space to support a random vector $\bm X_0$ such that $\L(\bm X_0 | \G_T)= \mu_0$, as well as a $d$-dimensional Brownian motion $\bm W$ independent of $\G_T \vee \sigma(\bm{X}_0)$. We endow this extended probability space $\widetilde{\Omega}$ with the filtration $\FF=(\F_t)_{t \in [0,T]}$ defined by $\F_t=\F^{\bm W}_t \vee \G_t\vee \sigma(\bm{X}_0)$ for $t\in[0,T]$. Then, we face the SDE \eqref{SDE_gen} with random coefficients which are Lipschitz in the $\bm x$-variable, with the Lipschitz constants satisfying an integrability condition in $t$ uniformly in the randomness. Under these circumstances, the existence of a strong solution (adapted to $\FF$) to the SDE \eqref{SDE_gen} can be shown along the lines of \cite[proof of Theorem 5.1.1]{stroock-varadhan}, and we postpone the proof to Appendix \ref{sec_appendix} for the sake of continuity of exposition. 

\medskip

To summarize, we have extended the probability space to $\big(\widetilde{\Omega},\FF,\widetilde{\PP}\big)$, where  $\bm B$ is an $\FF$-Brownian motion which is adapted to $\GG$. Moreover, the extended probability space supports a $\G_T$-independent $\FF$-Brownian motion $\bm W$ and a continuous $\FF$-adapted $d$-dimensional process $\bm X$ such that
\begin{equation*}
\d\bm{X}_t = \bm{b}(t,\omega,\bm{X}_t)\d t + \bm{\sigma}(t,\omega,\bm{X}_t)\d\bm{W}_t+\bm{\gamma}(t,\omega,\bm{X}_t)\d\bm{B}_t.
\end{equation*}
Since $\bm{B}$ is an $\FF$-Brownian motion, assertion (4) of Theorem \ref{Thm_SPDEtoSDE_general} holds, and clearly so does assertion (1). To check the compatibility assertion (3) of Theorem \ref{Thm_SPDEtoSDE_general}, it suffices to check property (2) of Lemma \ref{le:compatibility} with $\H_t:=\sigma(\bm{X}_0)$ for all $t \in [0,T]$. Property (2a) holds trivially because $\H_t=\H_0$ for all $t \in [0,T]$. Property (2c) is immediate by construction. Noting that $\L(\bm X_0 | \G_T)= \mu_0$ implies $\L(\bm X_0 | \G_T)= \L(\bm X_0 | \G_0)$, we deduce that property (2b) holds, and finally Lemma \ref{le:compatibility} ensures that $\F^{\bm X}_t \indep \F_T^{\bm{W}} \vee \G_T \,|\, \F^{\bm{W}}_t \vee \G_t$ for all $t \in [0,T]$.

\medskip

It remains to prove that $\mu_t=\L(\bm{X}_t\,|\,\G_T)=\L(\bm{X}_t\,|\,\G_t)$ a.s., for each $t \in [0,T]$. To this end, define $\nu_t=\L(\bm X_t|\G_T)$ for $t\in[0,T]$, and note that $\nu_0=\mu_0$. By Proposition \ref{pr:SDEtoSPDE}, $\nu$ satisfies $\nu_t=\L(\bm X_t|\G_t)$ for each $t$ and the SPDE
\begin{equation*}
\d\langle \nu_t, \varphi\rangle = \big\langle \nu_t, \gen^{\bm b,\bm a}_{t,\omega}\varphi\big\rangle \d t + \big\langle \nu_t,(\nabla\varphi)^\top \bm \gamma(t,\omega,\cdot)\big\rangle \d\bm B_t, \quad t \in [0,T], \quad\varphi \in C_c^\infty(\R^d).
\end{equation*}
Recall that $\mu$ satisfies the exact same linear SPDE \eqref{SPDE}. Thus, the proof is complete up to the following proposition, which establishes uniqueness for this SPDE. \qed
	
\begin{proposition}[Uniqueness for linear SPDEs with coefficients smooth in $\bm x$]\label{prop_uniquesmooth}
Suppose $\nu^{(1)},\,\nu^{(2)}$ are two probability measure-valued processes such that
\begin{equation}\label{uni_SPDE}
\d\langle \nu^{(i)}_t, \varphi\rangle = \big\langle \nu^{(i)}_t, \gen^{\bm b,\bm a}_{t,\omega}\varphi\big\rangle \d t + \big\langle \nu^{(i)}_t,(\nabla\varphi)^\top \bm \gamma(t,\omega,\cdot)\big\rangle \d\bm B_t,\; t \in [0,T], \; \varphi \in C_c^\infty(\R^d),
\end{equation}
for $i=1,2$, where the  coefficients satisfy
\begin{eqnarray} 
&& \int_0^T \big\|\bm b(t,\omega,\cdot)\big\|_{C_b^1(\R^d)}+\big\|\bm a(t,\omega,\cdot)\big\|_{C_b^2(\R^d)}+\big\|\bm \gamma(t,\omega,\cdot)\big\|_{C_b^2(\R^d)}^2\d t<\infty\;\;\text{a.s.,}\qquad\quad \\
&& \E\bigg[ \int_0^T \big\|\bm b(t,\omega,\cdot)\big\|_{C_b(\R^d)}+\big\|\bm a(t,\omega,\cdot)\big\|_{C_b(\R^d)}\d t\bigg]<\infty,\qquad\quad \label{exp_bnd}
\end{eqnarray}
and $\bm a - \bm\gamma\bm\gamma^\top$ is positive semidefinite and symmetric with symmetric square root $\bm\sigma$. If $\nu^{(1)}_0=\nu^{(2)}_0$ admits a density $\rho\in L_2(\R^d)$ with $\E[\|\rho\|_2^2]<\infty$, then $\nu^{(1)}_t=\nu^{(2)}_t$ for all $t\in[0,T]$ a.s.
\end{proposition}

\noindent\textbf{Proof.} The proof is similar to that of \cite[Theorem 3.4]{kurtz1999particle}, but requires the use of an additional Lyapunov function. We start by introducing the Gaussian densities $G_\delta(\bm x)=(2\pi\delta)^{-d/2} e^{-|\bm x|^2/(2\delta)}$, $\delta>0$ and the (random) functions
\begin{equation*}
Z_\delta(t)(\bm x):=(G_\delta*\nu_t)(\bm x)=\int_{\R^d}G_\delta(\bm x-\bm y)\,\nu_t(\d \bm y),\quad\delta>0,\quad t\in[0,T], 
\end{equation*}
for an arbitrary $\nu\in\{\nu^{(1)},\,\nu^{(2)},\,\nu^{(1)}-\nu^{(2)}\}$. Then, $Z_\delta(t)\in L_2(\R^d)$. Using the definition of $Z_\delta(t)$, Fubini's theorem, \eqref{uni_SPDE} and integration by parts we find further:
		\begin{equation*}
			\begin{aligned}
				\big\langle Z_\delta(t),\varphi\big\rangle=\big\langle Z_\delta(0),\varphi\big\rangle+\int_0^t\Big\langle \frac{1}{2}\sum_{i,j=1}^d\partial_{ij}\big(G_\delta\ast(a_s^{ij}\nu_s)\big)
				-\sum_{i=1}^d\partial_i\big(G_\delta*(b^i_s\nu_s)\big),\varphi\Big\rangle\d s
				\\
				-\int_0^t\sum_{j=1}^d\Big\langle\sum_{i=1}^d\partial_i\big(G_\delta*(\gamma^{ij}_s\nu_s)\big),\varphi\Big\rangle \d B^j_s,
			\end{aligned}
		\end{equation*}
with the shorthand notations $b_t^i,\,a_t^{ij},\,\gamma_t^{ij}$ for $b_i(t,\omega,\cdot),\,a_{ij}(t,\omega,\cdot),\,\gamma_{ij}(t,\omega,\cdot)$, respectively.
Next, letting $\|f\|_{W^{1,\infty}(\R^d)} := \|f\|_{C_b(\R^d)} + \|f\|_{\mathrm{Lip}(\R^d)}$ for bounded Lipschitz functions $f$ on $\R^d$, we define
\begin{equation}\label{def_Kt}
\begin{aligned}
& K_t(\omega)=K^*\Big(\|\bm b(t,\omega,\cdot)\|_{C_b^1(\R^d)}\!+\!\|\bm a(t,\omega,\cdot)\|_{C_b^2(\R^d)}\!+\!\|\bm \gamma(t,\omega,\cdot)\|^2_{C_b^2(\R^d)}\!+\!
\|\bm \sigma(t,\omega,\cdot)\|^2_{W^{1,\infty}(\R^d)}\Big), \\
&Y_t(\omega) =\exp\bigg(2\int_0^t K_s(\omega)\d s\bigg),
\end{aligned}
\end{equation}
for a constant $K^*<\infty$ to be chosen later. The estimates at the beginning of the proof of Proposition \ref{Thm_SPDEtoSDE_unismooth} show that $Y_T<\infty$ a.s. Thus, It\^o's formula yields:
\begin{equation*}
\begin{split}
\frac{\big\langle Z_\delta(t),\varphi\big\rangle^2}{Y_t}=\;&\frac{\big\langle Z_\delta(0),\varphi\big\rangle^2}{Y_0}+\int_0^t\frac{1}{Y_s}\sum_{j=1}^d\bigg\langle\sum_{i=1}^d\partial_i\big(G_\delta*(\gamma^{ij}_s\nu_s)\big),\varphi\bigg\rangle^2 \d s \\
&+\int_0^t\frac{2\langle Z_\delta(s),\varphi\rangle}{Y_s}\,\Big\langle \frac{1}{2} \sum_{i,j=1}^d\partial_{ij}\big(G_\delta*(a_s^{ij}\nu_s)\big)-\sum_{i=1}^d\partial_i\big(G_\delta*(b^i_s\nu_s)\big),\varphi\Big\rangle\d s \\
&-\!\int_0^t\sum_{j=1}^d\frac{2\langle Z_\delta(s),\varphi\big\rangle}{Y_s}\,\Big\langle\sum_{i=1}^d\partial_i\big(G_\delta*(\gamma^{ij}_s\nu_s)\big),\varphi\Big\rangle \!\d B^j_s
-\!\int_0^t\frac{2K_s\langle Z_\delta(s),\varphi\rangle^2}{Y_s} \d s,
\end{split}
\end{equation*}
where the functions inside the inner products belong to $L_2(\R^d)$ for a.e.~$s$. We take the expectation (thereby removing the martingale terms thanks to \eqref{exp_bnd}), sum over $\varphi$ in an orthonormal basis of $L_2(\R^d)$, and use $a^{ij}_s=\sum_{k=1}^d(\gamma^{ik}_s\gamma^{jk}_s+\sigma^{ik}_s\sigma^{jk}_s)$ to get
\begin{equation}\label{equa_expez}
\begin{split}
\E\bigg[\frac{\|Z_\delta(t)\|_2^2}{Y_t}\bigg]=&\;\E\bigg[\frac{\|Z_\delta(0)\|_2^2}{Y_0}\bigg]
+\int_0^t\E\bigg[\frac{1}{Y_s}\sum_{k=1}^d\Big\langle Z_\delta(s),\sum_{i,j=1}^d\partial_{ij}\big(G_\delta*(\gamma_s^{ik}\gamma_s^{jk}\nu_s)\big)\Big\rangle\bigg]\d s \\
&\;+\int_0^t\E\bigg[\frac{1}{Y_s}\sum_{k=1}^d\Big\langle Z_\delta(s),\sum_{i,j=1}^d\partial_{ij}\big(G_\delta*(\sigma_s^{ik}\sigma_s^{jk}\nu_s)\big)\Big\rangle\bigg]\d s \\
&\;-\int_0^t\E\bigg[\frac{2}{Y_s} \sum_{i=1}^d \Big\langle Z_\delta(s),\partial_i\big(G_\delta*(b^i_s\nu_s)\big)\Big\rangle\bigg]\d s \\
&\;+\int_0^t\E\bigg[\frac{1}{Y_s}\sum_{j=1}^d\bigg\|\sum_{i=1}^d\partial_i\big(G_\delta*(\gamma^{ij}_s\nu_s)\big)\bigg\|_2^2\bigg]\d s-\int_0^t\E\bigg[\frac{2K_s\|Z_\delta(s)\|_2^2}{Y_s}\bigg]\d s.
\end{split}
\end{equation}
We abbreviate the latter equation as 	
\begin{equation*}
\frac{\d}{\d t}\,\E\bigg[\frac{\|Z_\delta(t)\|_2^2}{Y_t}\bigg]=\E\bigg[\frac{\big\|G_\delta*|\nu_t|\big\|_2^2}{Y_t}\, b^Z_t\bigg]-\E\bigg[\frac{2K_t\|Z_\delta(t)\|_2^2}{Y_t}\bigg]
\end{equation*}
upon combining the first four integrals together (notice that $\|G_\delta*|\nu_t|\|_2=0$ implies that all inner products involved are $0$). By \cite[(3.8) in Lemma 3.3, (3.7) in Lemma 3.2]{kurtz1999particle} (taking $H$ to be a singleton therein) and writing $C_d$ for constants depending only on $d$ we have: 
\begin{equation}\label{inequ}
\begin{split}
&\,\Big\langle Z_\delta(s),\sum_{i,j=1}^d\partial_{ij}\big(G_\delta*(\gamma_s^{ik}\gamma_s^{jk}\nu_s)\big)\Big\rangle+\bigg\|\sum_{i=1}^d\partial_i\big(G_\delta*(\gamma^{ik}_s\nu_s)\big)\bigg\|_2^2 \\
&\qquad \leq C_d\big\|\bm\gamma(s,\omega,\cdot)\big\|^2_{C_b^2(\R^d)}\big\|G_\delta*|\nu_s|\big\|_2^2, \\
&\,\Big\langle Z_\delta(s),\sum_{i,j=1}^d\partial_{ij}\big(G_\delta*(\sigma_s^{ik}\sigma_s^{jk}\nu_s)\big)\Big\rangle\leq C_d\big\|\bm\sigma(s,\omega,\cdot)\big\|^2_{W^{1,\infty}(\R^d)}\big\|G_\delta*|\nu_s|\big\|_2^2, \\
&\,\Big|\Big\langle Z_\delta(s),\partial_i\big(G_\delta*(b^i_s\nu_s)\big)\Big\rangle\Big|
\leq C_d\big\|\bm b(s,\omega,\cdot)\big\|_{C_b^1(\R^d)}\big\|G_\delta*|\nu_s|\big\|^2_2.
\end{split}
\end{equation}
According to these inequalities we can choose the constant $K^*$ in \eqref{def_Kt} such that $b^Z_t\leq K_t$ a.s. Moreover, it holds:
\begin{equation}\label{equa_Zdelta}
\begin{aligned}
\frac{\d}{\d t}\,\E\bigg[\frac{\|Z_\delta(t)\|_2^2}{Y_t}\bigg]=\E\bigg[\frac{\big\|G_\delta*|\nu_t|\big\|_2^2}{Y_t}\big(b^Z_t-2K_t\big)\bigg]+\E\bigg[\frac{2K_t}{Y_t}\Big(\big\|G_\delta*|\nu_t|\big\|_2^2-\|Z_\delta(t)\|_2^2\Big)\bigg].
\end{aligned}
\end{equation}
Applying \eqref{equa_Zdelta} for the positive measures $\nu=\nu^{(1)}$ and $\nu=\nu^{(2)}$ we deduce that  $\frac{\d}{\d t}\,\E\big[\frac{\|Z^{(i)}_\delta(t)\|_2^2}{Y_t}\big]\leq0$, where $Z^{(i)}_\delta(t):=G_\delta*\nu_t^{(i)}$, for  $i=1,2$. Thus,
\begin{equation*}
\E\bigg[\frac{\|Z^{(i)}_\delta(t)\|_2^2}{Y_t}\bigg]\leq\E\big[\|Z^{(i)}_\delta(0)\|_2^2\big]\leq\E\big[\|\rho\|_2^2\big],\quad i=1,2.
\end{equation*} 
Letting $\varphi$ range over an orthonormal basis of $L_2(\R^d)$ consisting of continuous bounded functions and using Fatou's lemma we find:
\begin{equation*}
\E\bigg[\!\sum_{\varphi}\frac{\langle\nu_t^{(i)},\varphi\rangle^2}{Y_t}\bigg]
\!=\!\E\bigg[\!\sum_{\varphi}\frac{\lim_{\delta\to0}\langle Z_\delta^{(i)}(t),\varphi\rangle^2}{Y_t}\bigg]
\leq\liminf_{\delta\to0}\E\big[\|Z^{(i)}_\delta(0)\|_2^2\big]\leq\E\big[\|\rho\|_2^2\big],\;i\!=\!1,2.
\end{equation*}
In particular, for any $i$ and $t$, $\sum_{\varphi}{\langle\nu_t^{(i)},\varphi\rangle^2}<\infty$ a.s., which means that $\nu_t^{(i)}$ has a density $\rho^{(i)}_t\in L_2(\R^d)$ a.s. In addition, $\|Z_\delta^{(i)}(t)\|_2\leq\|\rho_t^{(i)}\|_2$ a.s. Plugging back into \eqref{equa_Zdelta} and relying on the dominated convergence theorem, we have for $\E\big[\frac{\|\rho^{(i)}_t\|_2^2}{Y_t}\big]$:
		\begin{equation*}
			\begin{aligned}
				0\leq\E\bigg[\frac{\|\rho^{(i)}_t\|_2^2}{Y_t}\bigg]&=\lim_{\delta\to0}\,\E\bigg[\frac{\|Z^{(i)}_\delta(t)\|_2^2}{Y_t}\bigg]
				\\
				&=\lim_{\delta\to0}\,\bigg(\E\bigg[\int_0^t\frac{\big\|Z^{(i)}_\delta(s)\big\|_2^2}{Y_s}\big(b^Z_s-2K_s\big)\d s\bigg]+\E\big[\|Z^{(i)}_\delta(0)\|_2^2\big]\bigg)
				\\
				&\leq\limsup_{\delta\to0}\,\bigg(-\E\bigg[\int_0^t\frac{\big\|Z^{(i)}_\delta(s)\big\|_2^2}{Y_s}\,K_s\d s\bigg]+\E\big[\|\rho\|_2^2\big]\bigg)
				\\
				&\leq \E\big[\|\rho\|_2^2\big]-\E\bigg[\int_0^t\frac{\big\|\rho^{(i)}_s\big\|_2^2}{Y_s}\,K_s\d s\bigg],
			\end{aligned}
		\end{equation*}
where the second to last inequality is due to $b_s^Z\le K_s$ and $\|Z^{(i)}_\delta(0)\|_2\leq\|\rho\|_2$; the last inequality follows from Fatou's lemma. All in all,
		\begin{equation}\label{equa_dct}
			\E\bigg[\int_0^t\frac{\big\|\rho^{(i)}_s\big\|_2^2}{Y_s}K_s\d s\bigg]\leq\E\big[\|\rho\|_2^2\big].
		\end{equation}
Now, we take $\nu=\nu^{(1)}-\nu^{(2)}$ in \eqref{equa_Zdelta} and note that, since  $\rho_t^{(1)},\rho_t^{(2)}\in L_2(\R^d)$, also $\rho_t^{(1)}-\rho_t^{(2)}\in L_2(\R^d)$. Thus, both  $\|Z_\delta(t)\|_2^2$ and $\|G_\delta*|\nu^{(1)}_t-\nu^{(2)}_t|\big\|_2^2$ are bounded by $2\|\rho^{(1)}_t\|_2^2+2\|\rho^{(2)}_t\|_2^2$ and converge to $\|\rho_t^{(1)}-\rho_t^{(2)}\|_2^2$ as $\delta\to0$. Passing to the $\delta\to0$ limit in \eqref{equa_Zdelta} and using the dominated convergence theorem with the help of \eqref{equa_dct} we get:
\begin{equation*}
\dfrac{\d}{\d t}\,\E\bigg[\frac{\|\rho_t^{(1)}-\rho_t^{(2)}\|_2^2}{Y_t}\bigg]\leq 0.
\end{equation*}  
Recalling that $\rho_0^{(1)}=\rho_0^{(2)}=\rho$ and $Y_t \in (0,\infty)$ a.s., we deduce that $\rho_t^{(1)}=\rho_t^{(2)}$ a.s.\ for each $t$. By the continuity of $t\mapsto\langle \nu^{(i)}_t,\varphi\rangle$ for all $\varphi\in C^\infty_c(\R^d)$ and $i=1,2$ (cf.~\eqref{uni_SPDE}), it must hold that $\nu_t^{(1)}=\nu_t^{(2)}$ for all $t\in[0,T]$ a.s. \qed
	

\section{Preliminaries on stable convergence} \label{se:stableconvergence}

The proof of Theorem \ref{Thm_SPDEtoSDE_general} follows an approximation-tightness-limit scheme similar to the one in \cite[proof of Theorem 2.6]{figalli2008existence} and \cite[proof of Theorem 2.5]{trevisan}.
In our case, however, we must keep track of the given probability space $(\Omega,\FF,\PP)$ throughout, and this complicates matters. To do this, we work with the notion of \emph{stable convergence}, which is quite natural in the context of SDEs with random coefficients \cite{jacod1981weak}.
We summarize here the minimal definitions and results we need for our purposes; see \cite{jacod1981type} or \cite[Section 8.10(xi)]{bogachev2007measure2} for further background and references, or \cite[Section 3-a]{jacod1981weak} for a summary.

\medskip

Consider for this paragraph a measurable space $(S,\SS)$ and a Polish space $E$ equipped with its Borel $\sigma$-algebra, and endow $S \times E$ with the product $\sigma$-algebra. We write $\P(S \stabtimes E)$ for the set of probability measures on $S \times E$ equipped with the \emph{stable topology}, which is the coarsest topology such that the maps $\P(S \stabtimes E) \ni m \mapsto  \langle m, h\rangle$ are continuous, where $h : S \times E \to \R$ is bounded and measurable with $h(s,\cdot)$ continuous for each $s \in S$.  With the stable topology, $\P(S \stabtimes E)$ is not metrizable in general. But, if $\SS$ is countably generated, then every compact set $K \subset \P(S \stabtimes E)$ is metrizable \cite[Proposition 8.10.64]{bogachev2007measure2}. (In fact, this is precisely why we assume in Theorem \ref{Thm_SPDEtoSDE_general} that $\G_T$ is countably generated.) For any probability measure $\lambda$ on $(S,\SS)$ and any set $K \subset \P(E)$ which is compact with respect to the topology of weak convergence, the set
\begin{align*}
\{m \in \P(S \stabtimes E) :\,m(\cdot \times E) = \lambda, \ m(S \times \cdot) \in K \}
\end{align*}
is compact in the stable topology (and thus also metrizable), see \cite[Theorem 2.8]{jacod1981type}. In particular, a set $K \subset \P(S \stabtimes E)$ is sequentially precompact with respect to stable convergence if all elements share the same $S$-marginal and their $E$-marginals are tight.

\medskip

We write $\P(S \stabtimes E)$, instead of simply $\P(S \times E)$, to emphasize the different roles played by $S$ and $E$ in the stable topology. For example, in the following lemma, the \emph{stable topology} on $\P(\Omega \times C([0,T];\R^d) \stabtimes C([0,T];\R^d))$ is the one generated by the functionals $P \mapsto \langle P,h\rangle$, where $h : \Omega \times C([0,T];\R^d) \times C([0,T];\R^d) \to \R$ is bounded, measurable and continuous in its final argument.

\medskip 

Throughout this section and the next, we take the filtered probability space $(\Omega,\GG,\PP)$ described in Theorem \ref{Thm_SPDEtoSDE_general} as given, and we make frequent use of the canonical space 
\begin{align}
\widetilde{\Omega} := \Omega \times C([0,T];\R^d) \times C([0,T];\R^d). \label{canonicalspace}
\end{align}
Any function $g$ on $\Omega$ extends naturally to $\widetilde{\Omega}$, as does any $\sigma$-algebra on $\Omega$, and we use the same notation in each case. The processes $\mu$ and $\bm{B}$ defined on $\Omega$ thus also live on $\widetilde{\Omega}$. Further, we define the canonical processes $\bm{W}=(\bm{W}_t)_{t \in [0,T]}$ and $\bm{X}=(\bm{X}_t)_{t \in [0,T]}$ on $\widetilde{\Omega}$:
\[
\bm{W}_t(\omega,{\bm w},{\bm x}) := {\bm w}_t, \quad \bm{X}_t(\omega,{\bm w},{\bm x}) :={\bm x}_t.
\]
As usual, the filtrations generated by these canonical processes are denoted $\FF^{\bm{W}}$ and $\FF^{\bm{X}}$, respectively.
We will work with various probability measures $\widetilde\PP$ on $\widetilde\Omega$. We write $\E_{\widetilde{\PP}}$ for the expectation with respect to a given $\widetilde\PP$, as well as $\L_{\widetilde\PP}(Z)$ for the law of a random variable $Z$ defined on $\widetilde\Omega$. An expectation symbol $\E$ without a subscript will always mean the expectation under $\PP$ of a random variable on $(\Omega,\G_T,\PP)$.

\begin{definition} \label{def:canonical}
Let $\P(\widetilde{\Omega};\PP)$ denote the set of probability measures $\widetilde\PP$ on $\widetilde{\Omega}$ such that:
\begin{enumerate}
\item The $\Omega$-marginal of $\widetilde\PP$ equals $\PP$.
\item $\bm{W}$ and $\bm{B}$ are independent $\GG \vee \FF^{\bm W} \vee \FF^{\bm X}$-Brownian motions under $\widetilde\PP$.
\item $\F^{\bm X}_t \indep \F^{\bm W}_T \vee \G_T \,|\,\F^{\bm W}_t \vee \G_t$, for each $t \in [0,T]$.
\item $\bm W$ is independent of $\G_T$.
\end{enumerate}
We always equip $\P(\widetilde{\Omega};\PP)$ with the stable topology it inherits as a subset of $\P(\Omega \times C([0,T];\R^d) \stabtimes C([0,T];\R^d))$.
\end{definition}

The first lemma will ensure that the conditions (2) and (3) of Definition \ref{def:canonical} are preserved by our approximation arguments in Section \ref{sec_proof_SPDEtoSDE}.

\begin{lemma} \label{lem:closed}
As a subset of $\P(\Omega \times C([0,T];\R^d) \stabtimes C([0,T];\R^d))$, the set $\P(\widetilde{\Omega};\PP)$  is closed in the stable topology.
\end{lemma}

\noindent\textbf{Proof.} This is immediate once one notices that conditions (2) and (3) of Definition \ref{def:canonical} may be recast as requiring that $\widetilde\PP \circ (\bm{W},\bm{B})^{-1}$ is the product of two Wiener measures and
\begin{align*}
0 &= \E_{\widetilde\PP}\left[f_t(\omega,\bm{W},\bm{X})\big(g_{t+}(\bm{W},\bm{B}) - \E_{\widetilde\PP}[g_{t+}(\bm{W},\bm{B})]\big) \right], \\
0 &= \E_{\widetilde\PP}\left[f_t(\omega,\bm{W},\bm{X}) \big(\psi(\omega,\bm{W}) - \E_{\widetilde\PP}[\psi(\omega,\bm{W})\,|\,\F^{\bm{W}}_t \vee \G_t](\omega,\bm{W})\big) \right],
\end{align*}
for each $t \in [0,T]$ and bounded functions $f_t$, $g_{t+}$, and $\psi$ which are measurable with respect to $\G_t \vee \F^{\bm W}_t \vee \F^{\bm X}_t$, $\sigma((\bm{W}_s-\bm{W}_t,\bm{B}_s-\bm{B}_t) : s \in [t,T])$, and $\G_T \vee \F^{\bm{W}}_T$, respectively, and with $x \mapsto f_t(\omega,{\bm w},{\bm x})$ continuous for each $(\omega,{\bm w})$. These are clearly closed constraints under the stable topology of $\P(\Omega \times C([0,T];\R^d) \stabtimes C([0,T];\R^d))$.
Indeed, the conditional expectation term poses no difficulty because, in light of conditions (1) and (4) of Definition \ref{def:canonical}, it does not in fact depend on the choice of $\widetilde\PP \in \P(\widetilde\Omega;\PP)$. \qed

\medskip

Our next lemma shows that the conditional measures $\L_{\widetilde\PP}(\bm{X}_t\,|\,\G_T)$ are sufficiently well-behaved under stable convergence:

\begin{lemma} \label{le:stable-condlaw}
Suppose $(\widetilde\PP^n)_{n\in\N}$ is a sequence in $\P(\widetilde\Omega;\PP)$ which converges in the stable topology to some $\widetilde\PP$. Let $\mu^n_t = \L_{\widetilde\PP^n}(\bm{X}_t\,|\,\G_T)$ for each $n$, which may be viewed as a random variable on $(\Omega,\G_T,\PP)$. If $\mu^n_t \to \mu_t$ weakly a.s.~for all $t \in [0,T]$, then $\mu_t=\L_{\widetilde\PP}(\bm{X}_t\,|\,\G_T)$ a.s.~for all $t \in [0,T]$.
\end{lemma}

\noindent\textbf{Proof.} For $t \in [0,T]$, a bounded $\G_T$-measurable random variable $Z : \Omega \to \R$, and a bounded continuous function $f : \R^d \to \R$, we have:
\begin{align*}
\qquad\quad\;\;\E\big[ Z \langle\mu_t,f\rangle\big] &= \lim_{n\to\infty}\E\big[ Z \langle\mu^n_t,f\rangle\big] = \lim_{n\to\infty}\E_{\widetilde\PP^n}\big[ Z f(\bm{X}_t)\big] = \E_{\widetilde\PP}\big[ Z f(\bm{X}_t)\big].\qquad\quad\;\;\qed
\end{align*}

The following proposition is a precompactness criterion for the stable topology, which we will use repeatedly.

\begin{proposition}\label{Prop_tightness}
For each $n\in\N$, suppose $(\bm b^n,\bm \sigma^n,\bm \gamma^n):\,[0,T]\times\Omega\times\R^d\to\R^d\times \R^{d\times d}\times \R^{d\times d}$ is measurable with respect to the product of the $\GG$-progressive $\sigma$-algebra on $[0,T]\times\Omega$ and the Borel $\sigma$-algebra on $\R^d$. Let $\widetilde\PP^n\in \P(\widetilde{\Omega};\PP)$ be such that
\begin{equation*}
\d\bm{X}_t = \bm{b}^{n}(t,\omega,\bm{X}_t)\d t + \bm{\sigma}^{n}(t,\omega,\bm{X}_t)\d\bm{W}_t + \bm{\gamma}^{n}(t,\omega,\bm{X}_t)\d\bm{B}_t,\quad t\in[0,T] \quad \widetilde\PP^n\text{-a.s.},
\end{equation*}
and let $\bm a^n:=\bm\sigma^n({\bm\sigma^n})^\top\!+\bm\gamma^n({\bm\gamma^n})^\top$.
If the family $\{\widetilde\PP^n \circ \bm X_0^{-1} : n \in \N\} \subset \P(\R^d)$ is tight, and if for some $p > 1$ we have	
\begin{align}
C:=\sup_{n\in\N}\,\E_{\widetilde\PP^n}\bigg[\int_0^T\big|\bm b^n(t,\omega,\bm X_t)\big|^{p}+\big|\bm a^n(t,\omega,\bm X_t)\big|^{p}\d t\bigg]<\infty, \label{tightness_cond}
\end{align}
then the sequence $(\widetilde\PP^n)_{n\in\N}$ admits a convergent subsequence in $\P(\widetilde{\Omega};\PP)$.
\end{proposition}

\noindent\textbf{Proof.} As discussed above, the compact sets in the stable topology are metrizable since $\G_T$ is countably generated by assumption. Every element of $\P(\widetilde{\Omega};\PP)$ is a probability measure on $\P(\Omega \times C([0,T];\R^d) \stabtimes C([0,T];\R^d))$ with the $\Omega \times C([0,T];\R^d)$-marginal equal to the product of $\PP$ and the Wiener measure. Hence, all we need to show is that $\{\widetilde\PP^n \circ \bm{X}^{-1} : n \in \N\} \subset \P(C([0,T];\R^d))$ is tight.
 To this end, for any $\delta > 0$ and stopping times $0\leq\tau_1\leq\tau_2\leq T$ with $\tau_2-\tau_1\leq\delta$ a.s., we estimate $\E_{\widetilde\PP^n}\big[|\bm X_{\tau_2}-\bm X_{\tau_1}|^{p}\big]$ with the help of the Burkholder-Davis-Gundy inequality by
\begin{equation*}
\begin{split}
&\,2^{p-1}\,\E_{\widetilde\PP^n}\Bigg[\bigg|\int_{\tau_1}^{\tau_2}\bm b^{n}(u,\omega,\bm X_u)\d u\bigg|^{p}+\bigg|\int_{\tau_1}^{\tau_2} \bm \gamma^{n}(u,\omega,\bm X_u)\d \bm  B_u+\int_{\tau_1}^{\tau_2} \bm\sigma^{n}(u,\omega,\bm X_u)\d\bm W_u\big)\bigg|^{p}\Bigg] \\
&\qquad \leq{C}_{p,d}\,\E_{\widetilde\PP^n}\Bigg[\bigg(\int_{\tau_1}^{\tau_2}\big|\bm b^{n}(u,\omega,\bm X_u)\big|\d u\bigg)^{p}+\bigg(\int_{\tau_1}^{\tau_2}\big|\bm a^{n}(u,\omega,\bm X_u)\big|\d u\bigg)^{\frac{p}{2}}\Bigg] \\
&\qquad \leq {C}_{p,d}\,\Bigg(\E_{\widetilde\PP^n}\bigg[\delta^{p-1}\int_{\tau_1}^{\tau_2}\big|\bm b^n(u,\omega,\bm X_u)\big|^p\d u\bigg]+\E_{\widetilde\PP^n}\bigg[\delta^{p-1}\int_{\tau_1}^{\tau_2} \big|\bm a^{n}(u,\omega,\bm X_u)\big|^p\d u\bigg]^{\frac{1}{2}}\Bigg) \\
&\qquad \leq {C}_{p,d}\big(\delta^{p-1}C+\sqrt{\delta^{p-1}C}\big),
\end{split}
\end{equation*}
where $C_{p,d} <\infty$ is a constant depending only on $p$ and $d$. Therefore, 
\begin{equation*}
\lim_{\delta\to 0}\,\limsup_{n\to\infty}\,\sup_{0\leq\tau_1\leq\tau_2\leq T:\,\tau_2-\tau_1<\delta}\,\E_{\PP^n}\Big[\big|\bm X_{\tau_2}-\bm X_{\tau_1}\big|^{p}\Big] =0.
\end{equation*}
In view of the Aldous tightness criterion (see \cite[Lemma 16.12, Theorems 16.11, 16.10, 16.5]{kallenberg2006foundations}), the sequence $(\widetilde\PP^n \circ \bm{X}^{-1})_{n \in \N}$ in $\P(C([0,T];\R^d))$ is indeed tight. \qed

\medskip

For the passage to the limit, we repeatedly use the following proposition. 

\begin{proposition}\label{prop_mtgconverge}
Let $(\bm{b},\bm{a},\bm{\gamma})$ be as in Theorem \ref{Thm_SPDEtoSDE_general}, and assume that $\bm b(t,\omega,\cdot)$, $\bm a(t,\omega,\cdot)$, and $\bm \gamma(t,\omega,\cdot)$ are continuous bounded functions for a.e.\ $(t,\omega)$. 
If a $\P(\widetilde\Omega;\PP)$-valued sequence $(\widetilde\PP^n)_{n\in\N}$ converges to some $\widetilde\PP$ in the stable topology, and if
\begin{equation}\label{p_unif}
\sup_{n\in\N}\,\E_{\widetilde\PP^n}\bigg[\int_0^T\|\bm b(t,\omega,\cdot)\|^p_{C(\R^d)}+\|\bm a(t,\omega,\cdot)\|^p_{C(\R^d)}+\|\bm\gamma(t,\omega,\cdot)\|^{2p}_{C(\R^d)}\d t\bigg]<\infty,
\end{equation}
then for any $\G_T$-measurable random variable $g:\,\Omega\to[-1,1]$, $0\le s<t\le T$, continuous $h^s:\,C([0,s];\R^d)^2\to[-1,1]$, and $\varphi\in C_c^\infty(\R^d)$, one has along a subsequence:\footnote{We abuse notation by writing $h^s(\bm w,\bm x)$ in place of $h^s(\bm w|_{[0,s]}, \bm x|_{[0,s]})$, for $(\bm w,\bm x) \in C([0,T];\R^d)^2$.}
\begin{equation*}
\begin{split}
&\lim_{n\to\infty}\E_{\widetilde\PP^n}\bigg[gh^s(\bm W,\bm X)\bigg(\!\varphi(\bm X_t)\!-\!\varphi(\bm X_s)
\!-\!\int_s^t \!\!\gen^{\bm b,\bm a}_{u,\omega}\varphi(\bm X_u)\d u\!-\!\int_s^t\!\nabla\varphi(\bm X_u)^\top\bm\gamma(u,\omega,\bm X_u)\d \bm B_u\!\bigg)\!\bigg] \\
&=\E_{\widetilde\PP}\bigg[gh^s(\bm W,\bm X)\bigg(\varphi(\bm X_t)\!-\!\varphi(\bm X_s)\!-\!\int_s^t \!\!\gen^{\bm b,\bm a}_{u,\omega}\varphi(\bm X_u)\d u\!-\!\int_s^t\!\nabla\varphi(\bm X_u)^\top\bm\gamma(u,\omega,\bm X_u)\d \bm B_u\!\bigg)\!\bigg].
\end{split}
\end{equation*}
\end{proposition}

\noindent\textbf{Proof.} Truncating $(\bm b,\bm a)$ to be uniformly bounded, exploiting the stable convergence, and removing the truncation with the help of \eqref{p_unif} we get:
\begin{equation*}
\begin{split}
&\,\lim_{n\to\infty}\,\E_{\widetilde\PP^n}\bigg[gh^s(\bm W,\bm X)\bigg(\varphi(\bm X_t)-\varphi(\bm X_s)-\int_s^t \gen^{\bm b,\bm a}_{u,\omega}\varphi(\bm X_u)\d u\bigg)\bigg] \\
&=\E_{\widetilde\PP}\bigg[gh^s(\bm W,\bm X)\bigg(\varphi(\bm X_t)-\varphi(\bm X_s)-\int_s^t \gen^{\bm b,\bm a}_{u,\omega}\varphi(\bm X_u)\d u\bigg)\bigg].
\end{split}
\end{equation*}

\smallskip

To justify the convergence of the It\^o integral part, we extend the probability space further to $\widehat\Omega := \widetilde\Omega \times C([0,T];\R)$, where the last factor accommodates the It\^o integral process
\begin{equation*}
M_t:=\int_0^t\nabla\varphi(\bm X_u)^\top\bm\gamma(u,\omega,\bm X_u)\d \bm B_u,\quad t\in[0,T]
\end{equation*} 
(which is a well-defined $\GG\vee\FF^{\bm X}$-martingale under any probability measure in $\P(\widetilde\Omega;\PP)$). Writing $\widehat{\PP}^n$ for the probability measure on
\[
\widetilde\Omega \times C([0,T];\R) = \Omega \times C([0,T];\R^d) \times C([0,T];\R^d) \times C([0,T];\R),
\]
we have that the marginal of $\widehat{\PP}^n$ on $\widetilde\Omega$ coincides with $\widetilde\PP^n$. Following the strategy of the proof of Proposition \ref{Prop_tightness} with $(\bm{X},M)$ in place of $\bm{X}$, we find a subsequential limit point $\widehat{\PP}$ with respect to the stable topology on $\P(\Omega \times C([0,T];\R^d) \stabtimes C([0,T];\R^d) \times C([0,T];\R))$.
For ease of exposition, we denote by $(\widehat{\PP}^n)_{n\in\N}$ a subsequence converging to $\widehat{\PP}$.

\medskip

Next, we claim that the process $M$ is actually a $\GG \vee\FF^{\bm{W}} \vee\FF^{\bm X}$-martingale under $\widehat{\PP}$. By \eqref{p_unif}, for each $t\in[0,T]$, the random variables $M_t$, $M_t^2$, and $M_t\bm B_t$ are uniformly integrable under $(\widehat{\PP}_n)_{n\in\N}$, since by the Burkholder-Davis-Gundy inequality,
\begin{equation*}
\begin{split}
\E_{\widehat{\PP}^n}\big[|M_t|^{2p}\big] &\leq C_{p,d,t}\,\E_{\widehat{\PP}^n}\bigg[
\int_0^t \big|\nabla\varphi(\bm X_u)^\top{\bm\gamma}(u,\omega,\bm X_u)\big|^{2p}\d u\bigg] \\ 
&\leq C_{p,d,t}\,\|\varphi\|^{2p}_{C_b^1(\R^d)}\,\E_{\widetilde\PP^n}\bigg[\int_0^t\|\bm\gamma(u,\omega,\cdot)\|^{2p}_{C(\R^d)}\d u\bigg],
\end{split}
\end{equation*}
for a constant $C_{p,d,t}$ depending on $p$, $d$ and $t$. The uniform integrability together with the stable convergence imply that, for $0 \le s < t \le T$ and any random variable $\psi$ which is the product of a bounded $(\Omega,\G_s)$-random variable and a bounded continuous function of $(\bm W_u,\bm X_u)_{u\in[0,s]}$, it holds: 
\begin{equation*}
\E_{\widehat{\PP}} \big[(M_t-M_s)\psi\big]=\lim_{n\to\infty}\,\E_{\widehat{\PP}^n}\big[(M_t-M_s)\psi\big]=0,
\end{equation*} 
so $M$ is a $\GG \vee\FF^{\bm{W}} \vee\FF^{\bm X}$-martingale under $\widehat{\PP}$. 

\medskip

We lastly argue that $M$ equals to the It\^o integral process $\int_0^t\nabla\varphi(\bm X_u)^\top\bm\gamma(u,\omega,\bm X_u)\d \bm B_u$.
To find the quadratic variation of $M$ under $\widehat{\PP}$ we compute, for $\psi$ as above,
\begin{equation*}
\begin{split}
\E_{\widehat{\PP}}\big[\big(M_t^2-M_s^2\big)\psi\big]=&\lim_{n\to\infty}\,\E_{\widehat{\PP}^n}\big[\big(M_t^2-M_s^2\big)\psi\big] \\
=&\lim_{n\to\infty}\,\E_{\widehat{\PP}^n}\bigg[\int_s^t\big|\nabla\varphi(\bm X_u)^\top\bm\gamma(u,\omega,\bm X_u)\big|^2\d u\,\psi\bigg] \\
=&\,\E_{\widehat{\PP}} \bigg[\int_s^t\big|\nabla\varphi(\bm X_u)^\top{\bm\gamma}(u,\omega,\bm X_u)\big|^2\d u\,\psi\bigg],
\end{split}
\end{equation*} 
where the last equation is due to stable convergence and \eqref{p_unif}. Similarly, we deduce that
\begin{equation*}
\begin{split}
\E_{\widehat{\PP}}\big[\big(M_t\bm B_t^\top-M_s\bm B_s^\top\big)\psi\big]
=&\lim_{n\to\infty}\,\E_{\widehat{\PP}^n}\big[\big(M_t\bm B_t^\top-M_s\bm B_s^\top\big)\psi\big] \\
=&\lim_{n\to\infty}\,\E_{\widehat{\PP}^n}\bigg[\int_s^t\nabla\varphi(\bm X_u)^\top{\bm\gamma}(u,\omega,\bm X_u)\d u\,\psi\bigg] \\
=&\,\E_{\widehat{\PP}}\bigg[\int_s^t\nabla\varphi(\bm X_u)^\top{\bm\gamma}(u,\omega,\bm X_u)\d u\,\psi\bigg].
\end{split}
\end{equation*} 
Hence, under $\widehat{\PP}$, we have:
\begin{equation*}
\frac{\d}{\d t}\langle M \rangle_t=\big|\nabla\varphi(\bm X_t)^\top{\bm\gamma}(t,\omega,\bm X_t)\big|^2,\quad \frac{\d}{\d t}\langle M,\bm B^\top\rangle_t=\nabla\varphi(\bm X_t)^\top{\bm\gamma}(t,\omega,\bm X_t).
\end{equation*}
At this point, we apply \cite[Chapter 3, Theorem 4.2]{karatzas-shreve} to find (possibly on a further extension of the underlying probability space) a standard $(d+1)$-dimensional Brownian motion $\widehat{\bm B}$ and an $\R^{(d+1)\times(d+1)}$-valued adapted process $\bm A$ such that
\begin{equation*}
\begin{split}
&M_t=\int_0^t\bm A_s^{d+1,\cdot}\d\widehat{\bm B}_s,\quad B^i_t=\int_0^t\bm A_s^{i,\cdot}\d\widehat{\bm B}_s,\;\;i=1,\ldots,d, \\
&\bm A_t \bm A_t^\top=\left(\begin{matrix}
\bm I_d & {\bm\gamma}(t,\omega,\bm X_t)^\top\nabla\varphi(\bm X_t) \\
\nabla\varphi(\bm X_t)^\top{\bm\gamma}(t,\omega,\bm X_t) & \big|\nabla\varphi(\bm X_t)^\top{\bm\gamma}(t,\omega,\bm X_t)\big|^2
\end{matrix}\right),
\end{split}
\end{equation*}
where $\bm A_s^{i,\cdot}$ denotes the $i$th row of $\bm A_s$.
Simple linear algebra then yields:
\begin{equation*}
\bm A_t^{d+1,\cdot}=\sum_{i=1}^d\nabla\varphi(\bm X_t)^\top\bm\gamma^{\cdot,i}(t,\omega,\bm X_t)\, \bm A_t^{i,\cdot},
\end{equation*}
with $\bm\gamma^{\cdot,i}$ denoting the $i$th column of $\bm\gamma$. This indeed reveals that $M$ is nothing but the It\^o integral process $\int_0^t\nabla\varphi(\bm X_u)^\top{\bm\gamma}(u,\omega,\bm X_u)\d \bm B_u$. By the uniform integrability of $M_t$ under $\{\widetilde\PP^n\}_{n\in\N}$ we have:
\begin{equation*}
\begin{split}
&\,\E_{\widetilde \PP}\bigg[gh^s(\bm W,\bm X)\int_s^t\nabla\varphi(\bm X_u)^\top{\bm\gamma}(u,\omega,\bm X_u)\d \bm B_u\bigg] \\
&\qquad =\E_{\widehat{\PP}}\big[gh^s(\bm W,\bm X)\big(M_t-M_s\big)\big]
=\lim_{n\to\infty}\E_{\widehat{\PP}^n}\big[gh^s(\bm W,\bm X)\big(M_t-M_s\big)\big] \\
&\qquad =\lim_{n\to\infty}\E_{{\widetilde \PP}^n}\bigg[gh^s(\bm W,\bm X)\int_s^t\nabla\varphi(\bm X_u)^\top{\bm\gamma}(u,\omega,\bm X_u)\d \bm B_u\bigg].
\end{split}
\end{equation*} 
The conclusion of the proposition readily follows. \qed


\section{Superposition from SPDE to SDE: general case} \label{sec_proof_SPDEtoSDE}

With preparations out of the way, in this section we establish the superposition principle in the sense of Theorem \ref{Thm_SPDEtoSDE_general} under gradually weaker additional conditions on the behavior of the coefficients in the $\bm x$-variable: smoothness in Proposition \ref{Thm_SPDEtoSDE_smooth}, boundedness in Proposition  \ref{Thm_SPDEtoSDE_bdd}, and local boundedness in Proposition \ref{Thm_SPDEtoSDE_locally bdd}. Similarly to \cite[proof of Theorem 2.6]{figalli2008existence} and \cite[proof of Theorem 2.5]{trevisan}, each weakening of the conditions is achieved by a suitable approximation-tightness-limit argument. 

\begin{proposition}[Theorem \ref{Thm_SPDEtoSDE_general}, smooth case]\label{Thm_SPDEtoSDE_smooth}
In the setting of Theorem \ref{Thm_SPDEtoSDE_general}, suppose that for a.e.~$(t,\omega)$ the functions $\bm b(t,\omega,\cdot)$, $\bm a(t,\omega,\cdot)$, and $\bm \gamma(t,\omega,\cdot)$ are bounded and twice continuously differentiable with bounded derivatives of the first and second order. If 
\begin{equation} \label{prop3.1_int}
\int_{0}^{T} \big\|\bm b(t,\omega,\cdot)\big\|^{p}_{C_b^1(\R^d)}+\big\|\bm a(t,\omega,\cdot)\big\|^{p}_{C_b^2(\R^d)}+\big\|\bm\gamma(t,\omega,\cdot)\big\|^{2p}_{C_b^2(\R^d)}\d t<\infty\;\;\text{a.s.},
\end{equation}
and $\mu_0$ admits a density $\rho\in L_2(\R^d)$ with $\E[\|\rho\|_2^2]<\infty$, then the conclusions of Theorem \ref{Thm_SPDEtoSDE_general} hold.
\end{proposition}

\noindent\textbf{Proof. Step 1 (Approximation).}
For $n \in \N$ and $t \in [0,T]$, define
\[
\lfloor t \rfloor_n = \max\{2^{-n}kT : k \in \{0,1,\ldots,2^n\}, \ 2^{-n}kT \le t\}.
\]
For a constant $D_{ p }<\infty$ as in Proposition  \ref{Thm_SPDEtoSDE_unismooth} and for $n \in \N$, let
\begin{align*}
\tau_n(\omega) := \inf & \Big\{t \in [0,T] :  \\
	&\int_{\lfloor t \rfloor_n}^{t} \big\|\bm b(s,\omega,\cdot)\big\|^{ p }_{C_b^1(\R^d)}\!+\!\big\|\bm a(s,\omega,\cdot)\big\|^{ p }_{C_b^2(\R^d)}\!+\!\big\|\bm \gamma(s,\omega,\cdot)\big\|^{2 p }_{C_b^2(\R^d)}\d s\!>\! D_{ p } \Big\}\!\wedge T
\end{align*}
and set $\mu^n_{\cdot}(\omega)=\mu_{\cdot\wedge\tau_n(\omega)}(\omega)$.
Then, $\tau_n$ is a $\GG$-stopping time for each $n$, and  $\lim_{n\to\infty} \tau_n=T$ a.s.~by \eqref{prop3.1_int}. Moreover, each process $\mu^n$ is $\mathbb G$-adapted and satisfies
\begin{equation*}
\d\langle \mu^n_t, \varphi\rangle = \big\langle \mu^n_t, \gen^{\bm b^n,\bm a^n}_{t,\omega}\varphi\big\rangle \d t + \big\langle \mu^n_t,(\nabla\varphi)^\top \bm \gamma^n(t,\omega,\cdot)\big\rangle \d\bm B_t, \quad t \in [0,T], \quad\varphi \in C_c^\infty(\R^d),
\end{equation*}
where $(\bm b^n,\bm a^n,\bm\gamma^n)(t,\omega,\bm x):=\mathbf{1}_{\{t\leq \tau_n(\omega)\}}\,(\bm b,\bm a, \bm \gamma)(t,\omega,\bm x)$. Thanks to Proposition \ref{Thm_SPDEtoSDE_unismooth} the superposition principle holds for every fixed $n$. Transferring to the canonical space $\widetilde\Omega$ introduced in \eqref{canonicalspace}, we find $\{\widetilde\PP^n:\,n\in\N\} \subset \P(\widetilde\Omega;\PP)$ (see Definition \ref{def:canonical}) such that, for each $n\in\N$, 
\begin{equation*}
\d\bm{X}_t = \bm{b}^n(t,\omega,\bm{X}_t)\d t + \bm{\sigma}^n(t,\omega,\bm{X}_t)\d\bm{W}_t + \bm{\gamma}^n(t,\omega,\bm{X}_t)\d\bm{B}_t \quad \widetilde\PP^n\text{-a.s.,}
\end{equation*}
with $\mu^n_t=\L_{\widetilde\PP^n}(\bm{X}_t\,|\,\G_T)=\L_{\widetilde\PP^n}(\bm{X}_t\,|\,\G_t)$ a.s.\ for all $t\in[0,T]$, and where $\bm \sigma^n$ is the symmetric square root of $\bm a^n-\bm \gamma^n({\bm\gamma^n})^\top$.

\medskip

\noindent\textbf{Step 2 (Tightness).} Next, we verify the tightness condition \eqref{tightness_cond} of Proposition \ref{Prop_tightness}:
\begin{equation*}
\begin{split}
\E_{\widetilde\PP^n}\bigg[\!\int_0^T\!\big|\bm b^n(t,\omega,\bm X_t)|^p\!+\!\big|\bm a^n(t,\omega,\bm X_t)\big|^p\d t\bigg]
&=\E\bigg[\!\int_0^T\Big\langle \mu_t^n,\big|\bm b^n(t,\omega,\cdot)\big|^p\!+\!\big|\bm a^n(t,\omega,\cdot)\big|^p\Big\rangle\d t\bigg] \\
&\leq\E\bigg[\int_0^T\Big\langle \mu_t,\big|\bm b(t,\omega,\cdot)\big|^p
+\big|\bm a(t,\omega,\cdot)\big|^p\Big\rangle\d t\bigg],
\end{split}
\end{equation*}
and the latter is finite by \eqref{integ}. In view of Proposition \ref{Prop_tightness}, we can find a limit point $\widetilde\PP \in \P(\widetilde\Omega;\PP)$ with respect to stable convergence. We relabel the subsequence so that $(\widetilde\PP^n)_{n\in\N}$ converges to $\widetilde\PP$.
Recall from Lemma \ref{lem:closed} that the set $\P(\widetilde\Omega;\PP)$ is closed in the stable topology, and in particular we have $\F^{\bm{X}}_t \indep \F^{\bm{W}}_T \vee \G_T\,|\,\F^{\bm{W}}_t \vee \G_t$ for each $t \in [0,T]$ under $\widetilde\PP$. By Lemma \ref{le:compatibility} (with $\HH=\FF^{\bm{X}}$), this implies in particular that $\F^{\bm{X}}_t \indep \G_T\,|\,\G_t$ for each $t \in [0,T]$, which we will use only at the very end of the proof.
Finally, since $\mu^n_t$ converges weakly a.s.\ to $\mu_t$ for each $t$, Lemma \ref{le:stable-condlaw} ensures that $\mu_t=\L_{\widetilde\PP}(\bm{X}_t\,|\,\G_T)=\L_{\widetilde\PP}(\bm{X}_t\,|\,\G_t)$ a.s.\ for each $t$.

\medskip

\noindent\textbf{Step 3 (Limit).} We proceed to derive an appropriate martingale problem for $\bm X$ under $\widetilde\PP$, in which we ``forget" the canonical Brownian motion $\bm{W}$ defined on $\widetilde\Omega$. To start with, we use It\^o's formula under $\widetilde\PP^n$ and obtain 
\begin{equation*}
\begin{split}
&\,\varphi(\bm X_t)-\varphi(\bm X_s)-\int_s^t \gen^{\bm b^n,\bm a^n}_{u,\omega}\varphi(\bm X_u)\d u
\\
&=\int_s^t\nabla\varphi(\bm X_u)^\top\bm\gamma^{n}(u,\omega,\bm X_u)\d \bm B_u+\int_s^t\nabla\varphi(\bm X_u)^\top\bm\sigma^{n}(u,\omega,\bm X_u)\d \bm W_u,\quad \varphi \in C_c^\infty(\R^d).
\end{split}
\end{equation*}
Now, fix a $\G_T$-measurable random variable $g:\,\Omega\to[-1,1]$, $0\le s<t\le T$, and a continuous function $h^s:\,C([0,s];\R^d)\to[-1,1]$. Since $\bm W$ is a Brownian motion independent of $\G_T$, we have (e.g., by Lemma \ref{le:fubini}):
\begin{equation*}
\E_{\widetilde\PP^n}\bigg[gh^s(\bm X)\bigg(\varphi(\bm X_t)-\varphi(\bm X_s)-\int_s^t \gen^{\bm b^n,\bm a^n}_{u,\omega}\varphi(\bm X_u)\d u-\int_s^t\nabla\varphi(\bm X_u)^\top\bm\gamma^{n}(u,\omega,\bm X_u)\d \bm B_u\bigg)\bigg]=0.
\end{equation*}

\smallskip

Further, for all $m\ge n \ge1$ we have $\tau_n\le\tau_m$, and thus
\begin{equation*}
\begin{split}
&\bigg|\E_{\widetilde\PP^m}\bigg[gh^s(\bm X)\int_s^t\nabla \varphi(\bm X_u)^\top \big(\bm b^m-\bm b^n\big)(u,\omega,\bm X_u)\d u\bigg]\bigg| \\
&\qquad \le \|\varphi\|_{C_b^1(\R^d)}\,\E\bigg[\int_s^t\big\langle \mu_u^m, |\bm b^m-\bm b^n|(u,\omega,\cdot)\big\rangle\d u\bigg] \\
&\qquad = \|\varphi\|_{C_b^1(\R^d)}\,\E\bigg[\int_s^t\big\langle \mu_u^m,  \mathbf{1}_{\{\tau_n< u\leq\tau_m\}}|\bm b(u,\omega,\cdot)|\big\rangle\d u\bigg] \\
&\qquad \leq \|\varphi\|_{C_b^1(\R^d)}\,\E\bigg[\int_s^t\big\langle \mu_u, \mathbf{1}_{\{\tau_n<u\}}|\bm b(u,\omega,\cdot)|\big\rangle\d u\bigg].
\end{split}
\end{equation*}
Similarly,
\begin{equation*}
\begin{split}
&\bigg|\E_{\widetilde\PP^m}\bigg[gh^s(\bm X)\int_s^t \nabla^2 \varphi(\bm X_u): (\bm a^m-\bm a^n)(u,\omega,\bm X_u)\d u\bigg]\bigg| \\
&\qquad \leq \|\varphi\|_{C_b^2(\R^d)}\,\E\bigg[\int_s^t\big\langle \mu_u, \mathbf{1}_{\{\tau_n<u\}}|\bm a(u,\omega,\cdot)|\big\rangle\d u\bigg].
\end{split}
\end{equation*}
For the stochastic integral part, we combine It\^o's isometry and Jensen's inequality to get
\begin{equation*}
\begin{split}
&\bigg|\E_{\widetilde\PP^m}\bigg[gh^s(\bm X)\int_s^t\nabla \varphi(\bm X_u)^\top (\bm \gamma^m-\bm \gamma^n)(u,\omega,\bm X_u)\d \bm B_u\bigg]\bigg| \\
&\qquad \leq  \|\varphi\|_{C_b^1(\R^d)}\,\E\bigg[\int_s^t\big\langle \mu_u^m, |\bm \gamma^m-\bm \gamma^n|(u,\omega,\cdot)^2\big\rangle\d u\bigg]^{1/2} \\
&\qquad \leq  \|\varphi\|_{C_b^1(\R^d)}\,\E\bigg[\int_s^t\big\langle \mu_u, \mathbf{1}_{\{\tau_n<u\}}|\bm \gamma(u,\omega,\cdot)|^2\big\rangle\d u\bigg]^{1/2}.
\end{split}
\end{equation*}
Putting everything together,  
\begin{equation*}
\begin{split}
& \bigg|\E_{\widetilde\PP^m}\bigg[gh^s(\bm X)\bigg(\!\varphi(\bm X_t)\!-\!\varphi(\bm X_s)\!-\!
\int_s^t \! \gen^{\bm b^n,\bm a^n}_{u,\omega}\varphi(\bm X_u)\d u
\!-\!\int_s^t \! \nabla\varphi(\bm X_u)^\top\bm\gamma^{n}(u,\omega,\bm X_u)\d \bm B_u\bigg)\bigg]\bigg| \\
&\qquad \leq \|\varphi\|_{C_b^2(\R^d)}\,\E\bigg[\int_s^t\big\langle \mu_u, \mathbf{1}_{\{\tau_n<u\}}\big(|\bm b(u,\omega,\cdot)|+|\bm a(u,\omega,\cdot)|\big)\big\rangle\d u\bigg] \\
&\qquad \quad\, + \|\varphi\|_{C_b^1(\R^d)}\,\E\bigg[\int_s^t\big\langle \mu_u, \mathbf{1}_{\{\tau_n<u\}}|\bm \gamma(u,\omega,\cdot)|^2\big\rangle\d u\bigg]^{1/2}.
\end{split} 
\end{equation*}

\smallskip

Take the $m\to\infty$ limit by means of Proposition \ref{prop_mtgconverge}, noting that \eqref{p_unif} holds for $(\bm b^n,\bm a^n,\bm \gamma^n)$, to get 
\begin{align}
&\,\bigg|\E_{\widetilde\PP}\bigg[gh^s(\bm X)\bigg(\!\varphi(\bm X_t)\!-\!\varphi(\bm X_s)\!-\!
\int_s^t \! \gen^{\bm b^n,\bm a^n}_{u,\omega}\varphi(\bm X_u)\d u
\!-\!\int_s^t \! \nabla\varphi(\bm X_u)^\top\bm\gamma^{n}(u,\omega,\bm X_u)\d \bm B_u\bigg)\bigg]\bigg| \nonumber \\
&\qquad \leq \|\varphi\|_{C_b^2(\R^d)}\,\E\bigg[\int_s^t\big\langle \mu_u, \mathbf{1}_{\{\tau_n<u\}}\big(|\bm b(u,\omega,\cdot)|+|\bm a(u,\omega,\cdot)|\big)\big\rangle\d u\bigg] \label{pf:bd1-1} \\
&\qquad \quad\, + \|\varphi\|_{C_b^1(\R^d)}\,\E\bigg[\int_s^t\big\langle \mu_u, \mathbf{1}_{\{\tau_n<u\}}|\bm \gamma(u,\omega,\cdot)|^2\big\rangle\d u\bigg]^{1/2}. \nonumber
\end{align}
As noted in Step 2, we have $\mu_t=\L_{\widetilde\PP}(\bm X_t\,|\,\G_t)$, and hence
\begin{equation*}
\begin{split}
&\bigg|\E_{\widetilde\PP}\bigg[gh^s(\bm X)\int_s^t\nabla \varphi(\bm X_u)^\top (\bm b^n-\bm b)(u,\omega,\bm X_u)+\frac{1}{2}\,\nabla^2 \varphi(\bm X_u): (\bm a^n-\bm a)(u,\omega,\bm X_u)\d u\bigg]\bigg| \\
&\qquad \leq \,\|\varphi\|_{C_b^2(\R^d)}\,\E\bigg[\int_s^t\big\langle \mu_u, \mathbf{1}_{\{\tau_n<u\}}\big(|\bm b(u,\omega,\cdot)|+|\bm a(u,\omega,\cdot)|\big)\big\rangle\d u\bigg], 
\end{split}
\end{equation*}
and, again by using It\^o isometry and Jensen's inequality,
\begin{equation*}
\begin{split}
&\,\bigg|\E_{\widetilde\PP}\bigg[gh^s(\bm X)\int_s^t\nabla \varphi(\bm X_u)^\top (\bm \gamma^n-\bm \gamma)(u,\omega,\bm X_u)\d \bm B_u\bigg]\bigg| \qquad\qquad\qquad\qquad\qquad\qquad\;\;\;\\
&\qquad \leq \|\varphi\|_{C_b^1(\R^d)}\,\E\bigg[\int_s^t\big\langle \mu_u, \mathbf{1}_{\{\tau_n<u\}}|\bm \gamma(u,\omega,\cdot)|^2\big\rangle\d u\bigg]^{1/2}.
\end{split}
\end{equation*} 
Consequently, returning to \eqref{pf:bd1-1}, we have:
\begin{equation*}
\begin{split}
&\,\bigg|\E_{\widetilde\PP}\bigg[gh^s(\bm X)\bigg(\varphi(\bm X_t)-\varphi(\bm X_s)-\int_s^t \gen^{\bm b,\bm a}_{u,\omega}\varphi(\bm X_u)\d u-\int_s^t\nabla\varphi(\bm X_u)^\top\bm\gamma(u,\omega,\bm X_u)\d \bm B_u\bigg)\bigg]\bigg| \\
&\qquad \leq 2\|\varphi\|_{C_b^2(\R^d)}\,\E\bigg[\int_s^t\big\langle \mu_u, \mathbf{1}_{\{\tau_n<u\}}\big(|\bm b(u,\omega,\cdot)|+|\bm a(u,\omega,\cdot)|\big)\big\rangle\d u\bigg] \\
&\qquad \quad+ 2  \|\varphi\|_{C_b^1(\R^d)}\,\E\bigg[\int_s^t\big\langle \mu_u, \mathbf{1}_{\{\tau_n<u\}}|\bm \gamma(u,\omega,\cdot)|^2\big\rangle\d u\bigg]^{1/2}.
\end{split} 
\end{equation*}
Finally, the a.s.~convergence $\tau_n \stackrel{n\to\infty}{\longrightarrow} T$ in conjunction with the dominated convergence theorem (recall \eqref{integ}) yield:
\begin{equation*}
\E_{\widetilde\PP}\bigg[gh^s(\bm X)\bigg(\varphi(\bm X_t)-\varphi(\bm X_s)-\int_s^t \gen^{\bm b,\bm a}_{u,\omega}\varphi(\bm X_u)\d u-\int_s^t\nabla\varphi(\bm X_u)^\top\bm\gamma(u,\omega,\bm X_u)\d \bm B_u\bigg)\bigg]=0.
\end{equation*}

\smallskip	

\noindent\textbf{Step 4 (Martingale problem to SDE).} We have thus shown that, for each $\varphi \in C^\infty_c(\R^d)$, the process
\begin{equation*}
\varphi(\bm X_t)-\varphi(\bm X_0)-\int_0^t \gen^{\bm b,\bm a}_{s,\omega}\varphi(\bm X_s)\d s-\int_0^t\nabla\varphi(\bm X_s)^\top\bm{\gamma}(s,\omega,\bm X_s)\d \bm B_s,\quad t\in[0,T]
\end{equation*} 
is a $(\G_T\vee\F^{\bm X}_t)_{t\in[0,T]}$-martingale under $\widetilde\PP$, and note that it is also adapted to the smaller filtration $\GG \vee \FF^{\bm{X}}=(\G_t\vee\F^{\bm{X}}_t)_{t \in [0,T]}$. Using linear and quadratic test functions as in \cite[Chapter 5, Proposition 4.6]{karatzas-shreve} we see further that
\begin{align*}
\bm{D}_t &:= \bm X_t-\bm X_0-\int_0^t\bm{b}(s,\omega,\bm{X}_s)\d s - \int_0^t\bm{\gamma}(s,\omega,\bm{X}_s)\d\bm{B}_s,\;t\in[0,T] \quad \text{and} \\
\bm{Q}_t &:= \bm{X}_t\bm{X}_t^\top - \int_0^t \big( \bm{X}_s\bm{b}(s,\omega,\bm{X}_s)^\top + \bm{b}(s,\omega,\bm{X}_s) \bm{X}_s^\top + \bm{a}(s,\omega,\bm{X}_s)\big)\d s \\
	&\quad\;\, - \int_0^t\bm{X}_s \big(\bm{\gamma}(s,\omega,\bm{X}_s) \d\bm{B}_s\big)^\top - \int_0^t\bm{\gamma}(s,\omega,\bm{X}_s) \d\bm{B}_s \, \bm{X}_s^\top,\;t\in[0,T]
\end{align*}
are both local martingales with respect to $(\G_T\vee\F^{\bm X}_t)_{t\in[0,T]}$ under $\widetilde\PP$.

\medskip

We claim next that the covariation $[\bm{D},\bm{B}]$ is identically zero (noting that quadratic variation does not depend on the underlying filtration). Indeed, $\bm{D}$ is a $(\G_T\vee\F^{\bm X}_t)_{t\in[0,T]}$-local martingale, so there is an increasing sequence of $(\G_T\vee\F^{\bm X}_t)_{t\in[0,T]}$-stopping times $\tau_n$ such that $\bm{D}_{t \wedge \tau_n}$ is a square-integrable martingale. Then, for $i,j=1,\ldots,d$ and $0\le s<t\le T$, the $\G_T$-measurability of $B^j_t$ implies
\begin{align*}
\E_{\widetilde\PP}[D^i_{t \wedge \tau_n} B^j_t \, |\, \G_s \vee \F^{\bm{X}}_s] = \E_{\widetilde\PP}\big[B^j_t \, \E_{\widetilde\PP}[D^i_{t \wedge \tau_n} \, |\, \G_T \vee \F^{\bm{X}}_s] \,|\, \G_s \vee \F^{\bm{X}}_s\big]
	&= D^i_{s \wedge \tau_n} \E_{\widetilde\PP}\big[B^j_t  \,|\, \G_s \vee \F^{\bm{X}}_s\big] \\
	&= D^i_{s \wedge \tau_n} B^j_s.
\end{align*}
We deduce that the covariation of $\bm{D}_{\cdot \wedge \tau_n}$ and $\bm{B}$ is identically zero for each $n$, and thus so is the covariation of $\bm{D}$ and $\bm{B}$.

\medskip

Using $[\bm{D},\bm{B}] \equiv 0$ and some simple stochastic calculus (in the smaller filtration $\GG \vee \FF^{\bm{X}}$), we deduce
\begin{align*}
\bm{Q}_t = \bm{Q}_0 + \int_0^t\bm{X}_s \d\bm{D}_s^\top + \int_0^t\d\bm{D}_s \,\bm{X}_s ^\top + [\bm{D}]_t - \int_0^t\bm{\sigma}(s,\omega,\bm{X}_s)^2\d s,\quad t\in[0,T]. 
\end{align*}
Since $\bm{D}$ and $\bm{Q}$ are continuous local martingales with respect to $\GG \vee \FF^{\bm{X}}$ under $\widetilde\PP$, we infer that the finite variation process
\[
[\bm{D}]_t - \int_0^t\bm{\sigma}(s,\omega,\bm{X}_s)^2\d s,\quad t\in[0,T]
\]
is also a continuous local martingale (in the same filtration), thus it is constant. Hence,
\[
[\bm{D}]_t = \int_0^t\bm{\sigma}(s,\omega,\bm{X}_s)^2\d s, \quad t\in [0,T].
\]
By \cite[Chapter 3, Theorem 4.2]{karatzas-shreve}, the $(\G_T \vee \F^{\bm{X}}_t)_{t \in [0,T]}$-local martingale $\bm{D}$ can be represented as
\begin{align}
\bm{D}_t = \int_0^t \bm{\sigma}(s,\omega,\bm{X}_s) \d\widehat{\bm W}_s, \quad t \in [0,T], \label{pf:Drep1}
\end{align}
where $\widehat{\bm W}$ is a $d$-dimensional Brownian motion. The construction of $\widehat{\bm W}$ in general requires enlarging the probability space $(\widetilde\Omega,\GG \vee \FF^{\bm{X}},\widetilde\PP)$, but, following the proof of \cite[Chapter 3, Theorem 4.2]{karatzas-shreve}, this additional randomness can be achieved merely by appending an independent Brownian motion $\widehat{\bm W}^0$; the desired Brownian motion $\widehat{\bm W}$ is then constructed to be adapted to the filtration $(\G_T \vee \F^{\bm{X}}_t \vee \F^{\widehat{\bm W}^0}_t)_{t \in [0,T]}$.
From the definition of $\bm{D}$ and its representation in \eqref{pf:Drep1} we deduce that $\bm{X}$ solves the desired SDE
\begin{equation*}
\bm X_t=\bm X_0+\int_0^t\bm b(s,\omega,\bm X_s)\d s+\int_0^t\bm\gamma(s,\omega,\bm X_s)\d\bm B_s+\int_0^t\bm\sigma(s,\omega,\bm X_s)\d \widehat{\bm W}_s.
\end{equation*}
%
%
This completes the proof that the conclusions (1) and (4) of Theorem \ref{Thm_SPDEtoSDE_general} are valid. Moreover, conclusion (2) has been shown in Step 2 above. To deduce the compatibility condition (3) of Theorem \ref{Thm_SPDEtoSDE_general}, it suffices to check condition (4) of Lemma \ref{le:compatibility} (with $\HH=\FF^{\bm X}$ and $\widehat{\bm W}$ in place of $\bm{W}$). We know that $\widehat{\bm W}$ is a Brownian motion with respect to the filtration $(\G_T \vee \F^{\bm{X}}_t \vee \F^{\widehat{\bm W}^0}_t)_{t \in [0,T]}$. Thus $\widehat{\bm W}$ is a Brownian motion in the smaller filtration $(\G_T \vee \F^{\bm{X}}_t \vee \F^{\widehat{\bm W}}_t)_{t \in [0,T]}$, with respect to which it is adapted.  This is condition (4a) of Lemma \ref{le:compatibility}.
Finally, recalling that $\widehat{\bm W}^0 \indep \G_T \vee \F^{\bm X}_T$ and $\F^{\bm X}_t \indep \G_T \,|\, \G_t$ for all $t \in [0,T]$, we easily deduce that $\F^{\widehat{\bm W}}_t \vee \F^{\bm X}_t \indep \G_T \,|\, \G_t$ for all $t \in [0,T]$, which is exactly condition (4b) of Lemma \ref{le:compatibility}. \qed


%
%


\begin{proposition}[Theorem \ref{Thm_SPDEtoSDE_general}, bounded case] \label{Thm_SPDEtoSDE_bdd}
In the setting of Theorem \ref{Thm_SPDEtoSDE_general}, if
\begin{equation}\label{prop3.2int}
\int_{0}^{T} \sup_{\bm x\in\R^d}|\bm b(t,\omega,\bm x)|^{p}+\sup_{\bm x\in\R^d}|\bm a(t,\omega,\bm x)|^{p}\,\d t<\infty\;\;\text{a.s.},
\end{equation}
then the conclusions of Theorem \ref{Thm_SPDEtoSDE_general} hold.
\end{proposition}

\noindent\textbf{Proof.} Our proof strategy is similar to the one in \cite[proof of Theorem 2.6]{figalli2008existence} and \cite[Section A.4, case of bounded coefficients]{trevisan}. We construct from $\mu$ a sequence of probability measure-valued processes $(\mu^n)_{n\in\N}$ associated to SPDEs with smooth and bounded coefficients and again perform a tightness-limit argument.

\medskip

\noindent\textbf{Step 1 (Approximation).} Let $\varrho$ be a probability density function belonging to $C^\infty_b(\R^d)$, symmetric about $0$, and such that $|D^k\varrho(\bm x)|\le M_{k}\varrho(\bm x)$ for some constant $M_k<\infty$ for each $k\in\N$. (For example, take $\varrho$ as a multiple of $e^{-\sqrt{1+|\bm x|^2}}$.) For each $n\in\N$, we set $\varrho_n(\bm x)\!=\!n^d\varrho(n\bm x)$ and define the probability measure-valued process $\mu_\cdot^n\!=\!\mu_\cdot \ast \varrho_n$, as well as 
\begin{equation}\label{def_coefficients_rhon}
\begin{split}
&  b^{\varrho_n}_i(t,\omega,\cdot)=\frac{\big(b_i(t,\omega,\cdot)\mu_t\big) *{\varrho_n}}{\mu_t * {\varrho_n}}=\frac{\int_{\R^d}\varrho_n(\cdot-\bm y)\,b_i(t,\omega,\bm y)\,\mu_t(\mathrm{d}\bm y)}{\int_{\R^d}\varrho_n(\cdot-\bm y)\,\mu_t(\mathrm{d}\bm y)}, \\
& a^{\varrho_n}_{ij}(t,\omega,\cdot)=\frac{\big(a_{ij}(t,\omega,\cdot)\mu_t\big) *{\varrho_n}}{\mu_t * {\varrho_n}}=\frac{\int_{\R^d}\varrho_n(\cdot-\bm y)\,a_{ij}(t,\omega,\bm y)\,\mu_t(\mathrm{d}\bm y)}{\int_{\R^d}\varrho_n(\cdot-\bm y)\,\mu_t(\mathrm{d}\bm y)}, \\
&\gamma^{\varrho_n}_{ij}(t,\omega,\cdot)=\frac{\big(\gamma_{ij}(t,\omega,\cdot)\mu_t\big) *{\varrho_n}}{\mu_t * {\varrho_n}}=\frac{\int_{\R^d}\varrho_n(\cdot-\bm y)\,\gamma_{ij}(t,\omega,\bm y)\,\mu_t(\mathrm{d}\bm y)}{\int_{\R^d}\varrho_n\,(\cdot-\bm y)\,\mu_t(\mathrm{d}\bm y)},
\end{split}
\end{equation}
for $i,j=1,\ldots,d$. Then, noting that $|\bm{\gamma}| \le \sqrt{d}\sqrt{|\bm{a}|}$, we have:
\begin{equation*}
\begin{split}
& \|\bm b^{\varrho_n}(t,\omega,\cdot)\|_{C_b^1(\R^d)}\leq 2nM_1\sup_{\bm x\in\R^d}|\bm b(t,\omega,\bm x)|, \\
& \|\bm a^{\varrho_n}(t,\omega,\cdot)\|_{C_b^2(\R^d)}\leq 2n^2(M_2 + 2M_1^2)\sup_{\bm x\in\R^d}|\bm a(t,\omega,\bm x)|, \\
& \|\bm \gamma^{\varrho_n}(t,\omega,\cdot)\|_{C_b^2(\R^d)} \leq 2n^2(M_2 + 2M_1^2)\sup_{\bm x\in\R^d}\!\sqrt{|\bm a(t,\omega,\bm x)|}.
\end{split}
\end{equation*}
It now follows from \eqref{prop3.2int} that, for each $n \in \N$,
\begin{equation*}
\int_0^T \|\bm b^{\varrho_n}(t,\omega,\cdot)\|^{ p }_{C_b^1(\R^d)}
+\|\bm a^{\varrho_n}(t,\omega,\cdot)\|^{ p }_{C_b^2(\R^d)}
+\|\bm \gamma^{\varrho_n}(t,\omega,\cdot)\|^{2 p }_{C_b^2(\R^d)}\d t<\infty\;\;\text{a.s.}
\end{equation*}
Moreover, by H\"older's inequality,
\begin{equation}\label{Hold1}
\begin{aligned}
		\int_{\R^d}\big|b^{\varrho_n}_i(t,\omega,\bm x)\big|^p\,\mu_t^n(\bm x)\d\bm x=&\int_{\R^d}\frac{\big|\int_{\R^d} \varrho_n(\bm x-\bm y)\,b_i(t,\omega,\bm y)\,\mu_t(\mathrm{d}\bm y)\big|^p}{\big|\int_{\R^d}\varrho_n(\bm x-\bm y)\,\mu_t(\mathrm{d}\bm y)\big|^{p-1}}\d\bm x \\
		\leq&\int_{\R^d}{\int_{\R^d} \varrho_n(\bm x-\bm y)\,\big|b_i(t,\omega,\bm y)\big|^p\,\mu_t(\mathrm{d}\bm y)}\d\bm x \\
		=&\int_{\R^d}\big|b_i(t,\omega,\bm y)\big|^p\,\mu_t(\mathrm{d}\bm y),
	\end{aligned}
\end{equation}
\begin{equation}\label{Hold2}
\int_{\R^d}\big|a^{\varrho_n}_{ij}(t,\omega,\bm x)\big|^p\,\mu_t^n(\bm x)\d\bm x
\le \int_{\R^d}\big|a_{ij}(t,\omega,\bm y)\big|^p\,\mu_t(\mathrm{d}\bm y), \qquad\qquad\qquad\;\;\,
\end{equation}
so that \eqref{integ} holds for $\bm b^{\varrho_n}$, $\bm a^{\varrho_n}$ and $\mu^n$. 

\medskip

In addition, a repeated application of Fubini's theorem gives
	\begin{equation*}
		\d\langle \mu^n_t, \varphi\rangle = \langle \mu^n_t, \gen^{\bm b^{\varrho_n},\bm a^{\varrho_n}}_{t,\omega}\varphi\rangle \d t + \big\langle \mu^n_t,(\nabla\varphi)^\top \bm \gamma^{\varrho_n}(t,\omega,\cdot)\big\rangle \d\bm B_t, \quad t \in [0,T], \quad \varphi \in C_c^\infty(\R^d),
	\end{equation*} 
with $\E[\|\mu_0^n\|^2_2]\le\sup_{x\in\R^d} \varrho_n(x)$ thanks to Jensen's inequality. Further, abbreviating $\bm a^{\varrho_n}(t,\omega,\cdot)$, $\bm \gamma^{\varrho_n}(t,\omega,\cdot)$, $\bm a(t,\omega,\cdot)$, $\bm \gamma(t,\omega,\cdot)$, and $\bm \sigma(t,\omega,\cdot)$ by $\bm a_t^{\varrho_n}$, $\bm \gamma_t^{\varrho_n}$, $\bm a_t$, $\bm \gamma_t$, and $\bm\sigma_t$, respectively, we find for any $\bm c =(c_1,\ldots,c_d)^\top\in\R^d$:
\begin{equation*}
\begin{aligned}
\bm c^\top \big(\bm a_t^{\varrho_n}-\bm\gamma^{\varrho_n}_t({\bm\gamma^{\varrho_n}_t})^\top\big)\bm c&=\frac{\big((\bm c^\top \bm a_t\bm c\,\mu_t)*\varrho_n\big)(\mu_t * {\varrho_n})-\big|(\bm c^\top\bm\gamma_t\mu_t)*\varrho_n\big|^2}{(\mu_t * {\varrho_n})^2} \\
&=\frac{\big((|\bm c^\top\bm\gamma_t|^2\mu_t+|\bm c^\top\bm \sigma_t|^2\mu_t)*\varrho_n\big)(\mu_t * {\varrho_n})-\big|(\bm c^\top\bm\gamma_t\mu_t)*\varrho_n\big|^2}{(\mu_t * {\varrho_n})^2} \\
&\geq\frac{\big(|\bm c^\top\bm\gamma_t|^2\mu_t*\varrho_n\big)(\mu_t * {\varrho_n})-\big|(\bm c^\top\bm\gamma_t\mu_t)*\varrho_n\big|^2}{(\mu_t * {\varrho_n})^2}\geq 0,
\end{aligned}
\end{equation*}
where the final estimate is due to the Cauchy-Schwarz inequality. Thus, $\bm a^{\varrho_n}-\bm\gamma^{\varrho_n}({\bm\gamma^{\varrho_n}})^\top$ is positive semidefinite and symmetric, and we can let $\bm \sigma^{\varrho_n}$ be its symmetric square root. All in all, the smoothed coefficients adhere to the conditions of Proposition \ref{Thm_SPDEtoSDE_smooth}, and thus the superposition principle applies for each $n$.
Transferring to the canonical space $\widetilde\Omega$ introduced in  \eqref{canonicalspace}, we find $\widetilde\PP^n\in\P(\widetilde\Omega;\PP)$ for $n\in\N$ (see Definition \ref{def:canonical}) such that, for each $n\in\N$, 
\begin{equation*}
\d\bm{X}_t = \bm{b}^{\varrho_n}(t,\omega,\bm{X}_t)\d t + \bm{\sigma}^{\varrho_n}(t,\omega,\bm{X}_t)\d\bm{W}_t + \bm{\gamma}^{\varrho_n}(t,\omega,\bm{X}_t)\d\bm{B}_t\quad \widetilde\PP^n\text{-a.s.,}
\end{equation*}
with $\mu^n_t=\L_{\widetilde\PP^n}(\bm{X}_t\,|\,\G_T)=\L_{\widetilde\PP^n}(\bm{X}_t\,|\,\G_t)$  for each $t\in[0,T]$, and with $\F^{\bm{X}}_t \indep \F^{\bm{W}}_T \vee \G_T \,|\,\F^{\bm{W}}_t \vee \G_t$ under $\widetilde\PP^n$ for each $t \in [0,T]$.
	
\medskip

\noindent\textbf{Step 2 (Tightness).} We aim to apply Proposition \ref{Prop_tightness} to $(\widetilde\PP^n)_{n\in\N}$. Since $\widetilde\PP^n \circ \bm{X}_0^{-1}$ converges weakly to $\E[\mu_0(\cdot)]$, we only need the existence of a uniform in $n$ bound on   
\begin{equation*}
		\E_{\widetilde\PP^n}\bigg[\int_0^T\big|\bm b^{\varrho_n}(t,\omega,\bm X_t)\big|^p+\big|\bm a^{\varrho_n}(t,\omega,\bm X_t)\big|^p\d t\bigg].
\end{equation*}
However, since $\mu^n_t=\L_{\widetilde\PP^n}(\bm{X}_t\,|\,\G_T)$ a.s., this follows directly from \eqref{Hold1}, \eqref{Hold2} and \eqref{integ}. Applying Proposition \ref{Prop_tightness} we obtain a limit point $\widetilde\PP \in \P(\widetilde\Omega;\PP)$ with respect to stable convergence. Once again, accepting an abuse of notation we relabel the subsequence so that $(\widetilde\PP^n)_{n\in\N}$  converges to $\widetilde\PP$.
Recall from Lemma \ref{lem:closed} that $\P(\widetilde\Omega;\PP)$ is closed in the stable topology, and in particular $\F^{\bm{X}}_t \indep \F^{\bm{W}}_T \vee \G_T\,|\,\F^{\bm{W}}_t \vee \G_t$ for each $t \in [0,T]$ under $\widetilde\PP$. By Lemma \ref{le:compatibility} (with $\HH=\FF^{\bm{X}}$), this implies that $\F^{\bm{X}}_t \indep \G_T\,|\,\G_t$ for each $t \in [0,T]$.
Lemma \ref{le:stable-condlaw} ensures that also $\mu_t=\L_{\widetilde\PP}(\bm{X}_t\,|\,\G_T)=\L_{\widetilde\PP}(\bm{X}_t\,|\,\G_t)$ a.s.~for each $t$.

\medskip

\noindent\textbf{Step 3 (Limit).} As in Step 3 in the proof of Proposition \ref{Thm_SPDEtoSDE_smooth}, our goal is now to derive an appropriate martingale problem for $\bm X$ under $\widetilde\PP$. With $\varphi$, $g$, $s$, $t$, and $h^s$ as therein,
\begin{equation}\label{equa1}
\E_{\widetilde\PP^n}\bigg[g h^s(\bm X)\bigg(\!\varphi(\bm X_t)-\varphi(\bm X_s)-\int_s^t\!\!\gen^{\bm b^{\varrho_n},\bm a^{\varrho_n}}_{u,\omega}\varphi(\bm X_u)\d u-\int_s^t\!\!\nabla\varphi(\bm X_u)^\top\bm\gamma^{\varrho_n}(u,\omega,\bm X_u)\d \bm B_u\!\bigg)\bigg]\!\!=\!0.
\end{equation}
Below we define a triple $(\widetilde{\bm{b}},\widetilde{\bm{a}},\widetilde{\bm{\gamma}}):\,[0,T]\times\Omega\times\R^d\to\R^d\times\R^{d\times d}\times\R^{d\times d}$ measurable with respect to the product of the $\mathbb G$-progressive $\sigma$-algebra on $[0,T]\times \Omega$ and the Borel $\sigma$-algebra on $\R^d$ such that all $(\widetilde{\bm{b}},\widetilde{\bm{a}},\widetilde{\bm{\gamma}})(t,\omega,\cdot)$ are continuous and compactly supported in $\R^d$, with their supremum norms being bounded uniformly by a constant $M<\infty$. For these, set
\begin{equation}\label{def_coefficients_rhontilde}
\begin{split}
& \widetilde{b}^{\varrho_n}_i(t,\omega,\cdot)=\frac{\big(\widetilde{b}_i(t,\omega,\cdot)\mu_t\big) *{\varrho_n}}{\mu_t * {\varrho_n}}=\frac{\int_{\R^d}\varrho_n(\cdot-\bm y)\,\widetilde{b}_i(t,\omega,\bm y)\,\mu_t(\mathrm{d}\bm y)}{\int_{\R^d}\varrho_n(\cdot-\bm y)\,\mu_t(\mathrm{d}\bm y)}, \\
& \widetilde{a}^{\varrho_n}_{ij}(t,\omega,\cdot)=\frac{\big(\widetilde{a}_{ij}(t,\omega,\cdot)\mu_t\big) *{\varrho_n}}{\mu_t * {\varrho_n}}=\frac{\int_{\R^d}\varrho_n(\cdot-\bm y)\,\widetilde{a}_{ij}(t,\omega,\bm y)\,\mu_t(\mathrm{d}\bm y)}{\int_{\R^d}\varrho_n(\cdot-\bm y)\,\mu_t(\mathrm{d}\bm y)}, \\
&	\widetilde{\gamma}^{\varrho_n}_{ij}(t,\omega,\cdot)=\frac{\big(\widetilde{\gamma}_{ij}(t,\omega,\cdot)\mu_t\big) *{\varrho_n}}{\mu_t * {\varrho_n}}=\frac{\int_{\R^d}\varrho_n(\cdot-\bm y)\,\widetilde{\gamma}_{ij}(t,\omega,\bm y)\,\mu_t(\mathrm{d}\bm y)}{\int_{\R^d}\varrho_n(\cdot-\bm y)\,\mu_t(\mathrm{d}\bm y)}.
\end{split}
\end{equation}

\smallskip

Next, we use the tower rule, a stochastic Fubini theorem (cf.~Lemma \ref{le:fubini}), Jensen's inequality, and It\^o isometry  to estimate an expectation similar to the one in \eqref{equa1} but with the new coefficients, for $0\le s<t\le T$: 
\begin{equation*}
\begin{split}
&\,\bigg|\E_{\widetilde\PP^n}\bigg[gh^s(\bm X)\bigg(\!\varphi(\bm X_t)\!-\!\varphi(\bm X_s)\!-\!\int_s^t\! \gen^{\widetilde{\bm b}^{\varrho_n},\widetilde{\bm a}^{\varrho_n}}_{u,\omega}\varphi(\bm X_u)\d u\!-\!\int_s^t\! \nabla\varphi(\bm X_u)^\top\widetilde{\bm\gamma}^{\varrho_n}(u,\omega,\bm X_u)\d \bm B_u\!\bigg)\bigg]\bigg| \\
&\qquad \leq\bigg|\E_{\widetilde\PP^n}\bigg[\int_s^t\E_{\widetilde\PP^n}\Big[gh^s(\bm X)\nabla\varphi(\bm X_u)^\top\big(\bm\gamma^{\varrho_n}-\widetilde{\bm\gamma}^{\varrho_n}\big)(u,\omega,\bm X_u)\Big|\G_T\Big] \d\bm B_u\bigg]\bigg| \\
&\qquad \quad\,+\E_{\widetilde\PP^n}\bigg[\int_s^t\big|\big(\gen^{\bm b^{\varrho_n},\bm a^{\varrho_n}}_{u,\omega}-\gen^{\widetilde{\bm b}^{\varrho_n},\widetilde{\bm a}^{\varrho_n}}_{u,\omega}\big)\varphi(\bm X_u)\big|\d u\bigg] \\
&\qquad  \leq \E_{\widetilde\PP^n}\bigg[\bigg(\int_s^t \Big|\E_{\PP^n}\Big[gh^s(\bm X)\nabla\varphi(\bm X_u)^\top\big(\bm\gamma^{\varrho_n}-\widetilde{\bm\gamma}^{\varrho_n}\big)(u,\omega,\bm X_u)\Big|\G_T\Big]\Big|^2 \d u\bigg)^{\!1/2}\bigg] \\
&\qquad \quad\,+\E_{\widetilde\PP^n}\bigg[\int_s^t\big|\big(\gen^{\bm b^{\varrho_n},\bm a^{\varrho_n}}_{u,\omega}-\gen^{\widetilde{\bm b}^{\varrho_n},\widetilde{\bm a}^{\varrho_n}}_{u,\omega}\big)\varphi(\bm X_u)\big|\d u\bigg] \\
&\qquad  \leq \E\bigg[\int_s^t\big\langle \mu_u^n,\big|(\nabla\varphi)^\top\big(\bm\gamma^{\varrho_n}-\widetilde{\bm\gamma}^{\varrho_n}\big)(u,\omega,\cdot)\big|\big\rangle^2 \d u\bigg]^{1/2} \\
&\qquad \quad\,+\E\bigg[\int_s^t\big\langle\mu_u^n,\big|\big(\gen^{\bm b^{\varrho_n},\bm a^{\varrho_n}}_{u,\omega}-\gen^{\widetilde{\bm b}^{\varrho_n},\widetilde{\bm a}^{\varrho_n}}_{u,\omega}\big)\varphi\big|\big\rangle \d u\bigg],
\end{split}
\end{equation*}
where we have exploited $|g|\le 1$ and $|h^s| \le 1$. Moreover, we control the latter term via the pointwise bound
\begin{align*}
\big\langle\mu_u^n,\big|\big(\bm b^{\varrho_n}\!-\!\widetilde{\bm b}^{\varrho_n} \big)(u,\omega,\cdot)\,\cdot \nabla\varphi\big|\big\rangle &= \int_{\R^d}\bigg|\int_{\R^d}\varrho_n(\bm x\!-\!\bm y)\,(\bm b-\widetilde{\bm b})(u,\omega,\bm y)\,\cdot \nabla\varphi(\bm x)\,\mu_u(\mathrm{d}\bm y)\bigg|\d\bm x  \\
&\leq \|\varphi\|_{C_b^1(\R^d)}\, \big\langle\mu_u,|\bm b - \widetilde{\bm b}|(u,\omega,\cdot)\big\rangle,
\end{align*}
where we have relied on the definitions of $ \bm b^{\varrho_n}$ and $\widetilde{\bm b}^{\varrho_n}$. Putting this together with an analogous bound for the $\bm a$-terms we arrive at 
\begin{equation*}
\E\bigg[\int_s^t\big\langle\mu_u^n,\big|\big(\gen^{\bm b^{\varrho_n},\bm a^{\varrho_n}}_{u,\omega}-\gen^{\widetilde{\bm b}^{\varrho_n},\widetilde{\bm a}^{\varrho_n}}_{u,\omega}\big)\varphi\big|\big\rangle \d u\bigg]
\leq \|\varphi\|_{C_b^2(\R^d)}\,\E\bigg[\int_s^t\big\langle\mu_u, \big|\gen^{\bm b-\widetilde{\bm b},\bm a-\widetilde{\bm a}}_{u,\omega}\big|\big\rangle\d u\bigg],
\end{equation*}
with the shorthand notation
\begin{align}
\big|\gen^{\bm b-\widetilde{\bm b},\bm a-\widetilde{\bm a}}_{u,\omega}\big|(\bm x)=|\bm b(u,\omega,\bm x)-\widetilde{\bm b}(u,\omega,\bm x)|+\big|\bm a(u,\omega,\bm x)-\widetilde{\bm a}(u,\omega,\bm x)\big|. \label{Lshorthand}
\end{align}
Similarly, for the $\bm \gamma$-terms we have the pointwise bound
\begin{align*}
\big\langle \mu_u^n,\big|(\nabla\varphi)^\top\big(\bm\gamma^{\varrho_n}-\widetilde{\bm\gamma}^{\varrho_n}\big)(u,\omega,\cdot)\big|\big\rangle \leq \|\varphi\|_{C_b^1(\R^d)} \big\langle\mu_u,|\bm\gamma-\widetilde{\bm\gamma}|(u,\omega,\cdot)\big\rangle .
\end{align*}
All in all, we conclude that
\begin{equation}\label{equa2}
\begin{split}
&\,\bigg|\E_{\widetilde\PP^n}\bigg[gh^s(\bm X)\bigg(\!\varphi(\bm X_t)\!-\!\varphi(\bm X_s)\!-\!\int_s^t\! \gen^{\widetilde{\bm b}^{\varrho_n},\widetilde{\bm a}^{\varrho_n}}_{u,\omega}\varphi(\bm X_u)\d u\!-\!\int_s^t\! \nabla\varphi(\bm X_u)^\top\widetilde{\bm\gamma}^{\varrho_n}(u,\omega,\bm X_u)\d \bm B_u\!\bigg)\bigg]\bigg| \\
&\quad \leq \|\varphi\|_{C_b^2(\R^d)}\,\bigg(\E\bigg[\int_s^t\big\langle\mu_u,\big|\bm\gamma-\widetilde{\bm\gamma}\big|(u,\omega,\cdot)\big\rangle^2\d u\bigg]^{1/2}+\E\bigg[\int_s^t\big\langle\mu_u, \big|\gen^{\bm b-\widetilde{\bm b},\bm a-\widetilde{\bm a}}_{u,\omega}\big|\big\rangle\d u\bigg]\bigg).
\end{split}
\end{equation}

We aim to take the $n\to\infty$ limit on the left-hand side of \eqref{equa2}. To this end, note that
\begin{align}
\bigg|\E_{\widetilde\PP^n}\bigg[gh^s(\bm X)\int_s^t &\big(\gen^{\widetilde{\bm b}^{\varrho_n},\widetilde{\bm a}^{\varrho_n}}_{u,\omega}-\gen^{\widetilde{\bm b},\widetilde{\bm a}}_{u,\omega}\big)\varphi(\bm X_u)\d u\bigg]\bigg| 
\leq \,\int_s^t\E_{\widetilde\PP^n}\Big[\big|\big(\gen^{\widetilde{\bm b}^{\varrho_n},\widetilde{\bm a}^{\varrho_n}}_{u,\omega}-\gen^{\widetilde{\bm b},\widetilde{\bm a}}_{u,\omega}\big)\varphi(\bm X_u)\big|\Big]\d u \nonumber \\
&\qquad\qquad  \leq \|\varphi\|_{C_b^2(\R^d)}\,\int_s^t\E\Big[\big\langle\mu_u^n, \big|\gen^{\widetilde{\bm b}^{\varrho_n}-\widetilde{\bm b},\widetilde{\bm a}^{\varrho_n}-\widetilde{\bm a}}_{u,\omega}\big|\big\rangle\Big]\d u, \label{equa3}
\end{align}
and also
\begin{align}
&\,\bigg|\E_{\widetilde\PP^n}\bigg[gh^s(\bm X)\int_s^t\nabla\varphi(\bm X_u)^\top\big(\widetilde{\bm\gamma}^{\varrho_n}(u,\omega,\bm X_u)-\widetilde{\bm\gamma}(u,\omega,\bm X_u)\big)\d \bm B_u\bigg]\bigg| \nonumber \\
&\qquad \leq \E_{\widetilde\PP^n}\bigg[\int_s^t\E_{\widetilde\PP^n}\Big[gh^s(\bm X)\nabla\varphi(\bm X_u)^\top\big(\widetilde{\bm\gamma}^{\varrho_n}-\widetilde{\bm\gamma}\big)(u,\omega,\bm X_u)\big|\G_T\Big]^2\d u\bigg]^{1/2} \qquad\qquad\qquad\qquad\; \nonumber \\
&\qquad  \leq \|\varphi\|_{C_b^2(\R^d)}\,\E\bigg[\int_s^t\big\langle\mu_u^n,\big|\widetilde{\bm\gamma}^{\varrho_n}-\widetilde{\bm\gamma}\big|(u,\omega,\cdot)\big\rangle^2\d u\bigg]^{1/2}.  \label{equa4}
\end{align}
In addition, the definition of $\widetilde{a}^{\varrho_n}_{ij}$ yields
\begin{equation*}
\begin{split}
&\,\big\langle\mu_u^n,\big|\widetilde{a}^{\varrho_n}_{ij}- \widetilde{a}_{ij}\big|(u,\omega,\cdot)\big\rangle \\
&\quad =\int_{\R^d}\int_{\R^d}\bigg|\frac{\int_{\R^d} \varrho^n(\bm x-\bm y)\,\widetilde{a}_{ij}(u,\omega,\bm y)\,\mu_t(\mathrm{d}\bm y)}{\int_{\R^d}\varrho^n(\bm x-\bm y)\,\mu_t(\mathrm{d}\bm y)} - \widetilde{a}_{ij}(u,\omega,\bm x)\bigg|\,\varrho^n(\bm x-\bm z)\,\mu_t(\mathrm{d}\bm z)\d\bm x \\
&\quad \leq\int_{\R^d}{\int_{\R^d}\big|\widetilde{a}_{ij}(u,\omega,\bm y)-\widetilde{a}_{ij}(u,\omega,\bm x)\big|\,\varrho^n(\bm x-\bm y)\,\mu_t(\mathrm{d}\bm y)}\d\bm x.
\end{split}
\end{equation*}
Fubini's theorem, the continuity of $\widetilde{a}_{ij}$, and the dominated convergence theorem (recall $|\widetilde{a}_{ij}|\le M$) imply that the latter upper bound tends to $0$ as $n\to\infty$, and so do the corresponding $\widetilde{\bm b}$- and $\widetilde{\bm\gamma}$-terms. Another application of the dominated convergence theorem gives the convergence to $0$ of the upper bounds in \eqref{equa3}, \eqref{equa4} as $n\to\infty$. Thus, the $n\to\infty$ limit of the left-hand side in \eqref{equa2} is the same as that of 
\begin{equation*}
\bigg|\E_{\widetilde\PP^n}\bigg[gh^s(\bm X)\bigg(\!\varphi(\bm X_t)\!-\!\varphi(\bm X_s)\!-\!\int_s^t\! \gen^{\widetilde{\bm b},\widetilde{\bm a}}_{u,\omega}\varphi(\bm X_u)\d u\!-\!\int_s^t\! \nabla\varphi(\bm X_u)^\top\widetilde{\bm\gamma}(u,\omega,\bm X_u)\d \bm B_u\!\bigg)\bigg]\bigg|.
\end{equation*}

\smallskip

At this point, we can pass to the $n\to\infty$ limit by means of Proposition \ref{prop_mtgconverge}, which in view of \eqref{equa2} results in 
\begin{equation}\label{equa5}
\begin{split}
&\, \bigg|\E_{\widetilde\PP}\bigg[gh^s(\bm X)\bigg(\!\varphi(\bm X_t)\!-\!\varphi(\bm X_s)\!-\!\int_s^t\! \gen^{\widetilde{\bm b},\widetilde{\bm a}}_{u,\omega}\varphi(\bm X_u)\d u\!-\!\int_s^t\! \nabla\varphi(\bm X_u)^\top\widetilde{\bm\gamma}(u,\omega,\bm X_u)\d \bm B_u\!\bigg)\bigg]\bigg| \\
& \leq \|\varphi\|_{C_b^2(\R^d)}\,\bigg(\E\bigg[\int_s^t\big\langle\mu_u,\big|\bm\gamma-\widetilde{\bm\gamma}\big|(u,\omega,\cdot)\big\rangle^2\d u\bigg]^{1/2}+\E\bigg[\int_s^t\big\langle\mu_u, \big|\gen^{\bm b-\widetilde{\bm b},\bm a-\widetilde{\bm a}}_{u,\omega}\big|\big\rangle\d u\bigg]\bigg).
\end{split}
\end{equation}
Now, as noted in Step 2, we have $\mu_u=\L_{\widetilde\PP}(\bm X_u|\G_T)=\L_{\widetilde\PP}(\bm X_u|\G_u)$ a.s. Using this and the same method as in \eqref{equa3} and \eqref{equa4}, we may change $(\widetilde{\bm b},\widetilde{\bm a},\widetilde{\bm \gamma})$ to $(\bm b,\bm a,\bm\gamma)$ on the left-hand side of \eqref{equa5}:
\begin{equation}\label{equa_consist}
\begin{split}
&\, \bigg|\E_{\PP}\bigg[gh^s(\bm X)\bigg(\!\varphi(\bm X_t)\!-\!\varphi(\bm X_s)\!-\!\int_s^t\! \gen^{\bm b,\bm a}_{u,\omega}\varphi(\bm X_u)\d u\!-\!\int_s^t\! \nabla\varphi(\bm X_u)^\top\bm\gamma(u,\omega,\bm X_u)\d \bm B_u\!\bigg)\bigg]\bigg| \\
& \leq 2\,\|\varphi\|_{C_b^2(\R^d)}\,\bigg(\E\bigg[\int_s^t\big\langle\mu_u,\big|\bm\gamma-\widetilde{\bm\gamma}\big|(u,\omega,\cdot)\big\rangle^2\d u\bigg]^{1/2}+\E\bigg[\int_s^t\big\langle\mu_u, \big|\gen^{\bm b-\widetilde{\bm b},\bm a-\widetilde{\bm a}}_{u,\omega}\big|\big\rangle\d u\bigg]\bigg).
\end{split}	
\end{equation}
	
\smallskip	
	
We now specify the way we choose the triple $(\widetilde{\bm{b}},\widetilde{\bm{a}},\widetilde{\bm{\gamma}})$. Take $\widetilde{a}_{ij}(t,\omega,\cdot)$ as an example. For any desired bound $M<\infty$ on $|\widetilde{a}_{ij}|$ and any fixed $\varepsilon>0$, we define for each $K\in\N$: 
\begin{equation*}
\begin{split}
& C^K_c(\R^d)=\big\{\eta\in C_c(\R^d):\;\|\eta\|_{C(\R^d)}\leq M,\;\;\|\eta\|_{\mathrm{Lip}(\R^d)}\leq K,\;\;\mathrm{supp}(\eta)\subset[-K,K]^d\big\}, \\
& f_K:\;[0,T]\times\Omega\times C^K_c(\R^d)\to\R, \\
&\qquad\,\,(t,\omega,\eta)\mapsto \big(\big\|\eta-\mathbf{1}_{\{|a_{ij}(t,\omega,\cdot)|\leq M\}}\, a_{ij}(t,\omega,\cdot)\big\|_{L_1(\mu_t(\omega))}-\varepsilon\big)_+, \\
& 	m_K(t,\omega)=\inf_{\eta\in C^K_{c}(\R^d)} f_K(t,\omega,\eta).
\end{split}
\end{equation*}	
Then $C^K_{c}(\R^d)$ is compact in the Polish space $C([-K,K]^d)$, and $f(t,\omega,\eta)$ is continuous in $\eta$. Therefore, by the measurable maximum theorem \cite[Theorem 18.19]{guide2006infinite}, $m_K$ is measurable and there exists a measurable selection $\eta^*_K$ such that
\begin{equation*}
f_K\big(t,\omega,\eta^*_K(t,\omega)\big)=m_K(t,\omega),\quad(t,\omega)\in[0,T]\times\Omega. 
\end{equation*}
Further, we let $S_0:=\varnothing$ and
\begin{equation*}
S_K:=\Big\{(t,\omega):\exists\,\eta\in C_c^{K}(\R^d)\;\;\text{s.t.}\;\;\big\|\eta-\mathbf{1}_{\{|a_{ij}(t,\omega,\cdot)|\leq M\}}\, a_{ij}(t,\omega,\cdot)\big\|_{L_1(\mu_t(\omega))}\leq\varepsilon\Big\},\quad K\in\N 
\end{equation*}
and notice that every $S_K$ is measurable as the preimage of $\{0\}$ under the measurable $m_K$. Moreover, $(S_K)_{K\in\N}$ is a non-decreasing set sequence with $\bigcup_{K\in \N} S_K=[0,T]\times\Omega$ because every $L_1(\mu_t(\omega))$-function can be approximated in the $L_1(\mu_t(\omega))$-norm by compactly supported Lipschitz functions. Finally, we define $\widetilde{a}_{ij}$ according to
\begin{equation}\label{equa_tilde}
\widetilde{a}_{ij}(t,\omega,\bm x)=\eta^*_K(t,\omega)(\bm x),\quad (t,\omega)\in S_K\backslash S_{K-1},\quad K\in\N. 
\end{equation}

\smallskip

The measurability of $\widetilde{a}_{ij}$ follows from the measurability of the restrictions of $\widetilde{a}_{ij}$ to $S_K\backslash S_{K-1}$, which in turn are measurable as compositions of the measurable mappings $(\eta^*_K,\mathrm{Id}):\,S_K\backslash S_{K-1}\times\R^d\to C([-K,K]^d)\times\R^d$ and $C([-K,K]^d)\times\R^d\to\R^d:(\eta,\bm x)\mapsto\eta(\bm x)$. In addition, we have
\begin{equation*}
\big\langle\mu_t(\omega),\big|\widetilde{a}_{ij}-a_{ij}\,\mathbf{1}_{\{|a_{ij}|\leq M\}}\big|(t,\omega,\cdot)\big\rangle\leq\varepsilon,\quad(t,\omega)\in[0,T]\times\Omega,
\end{equation*}
and thus
\begin{equation*}
\begin{split}
&\;\E\Big[\big\langle\mu_t,\big|a_{ij}-\widetilde{a}_{ij}\big|(t,\omega,\cdot)\big\rangle\Big] 
\leq\varepsilon + \E\Big[\big\langle\mu_t,|a_{ij}|\,\mathbf{1}_{\{|a_{ij}|>M\}}(t,\omega,\cdot)\big\rangle\Big].
\end{split}
\end{equation*}
Constructing $\widetilde{b}_i$ and $\widetilde{\gamma}_{ij}$ in a similar fashion we see that the right-hand side of \eqref{equa_consist} can be made arbitrarily small by choosing $\varepsilon>0$ small enough and $M<\infty$ large enough.

\medskip

In conclusion,
\begin{equation*}
\E_{\PP}\bigg[gh^s(\bm X)\bigg(\!\varphi(\bm X_t)\!-\!\varphi(\bm X_s)\!-\!\int_s^t\! \gen^{\bm b,\bm a}_{u,\omega}\varphi(\bm X_u)\d u\!-\!\int_s^t\! \nabla\varphi(\bm X_u)^\top\bm\gamma(u,\omega,\bm X_u)\d \bm B_u\!\bigg)\bigg]=0,
\end{equation*}
which means that
\begin{equation*}
\varphi(\bm X_t)-\varphi(\bm X_0)-\int_0^t \gen^{\bm b,\bm a}_{u,\omega}\varphi(\bm X_u)\d u-\int_0^t \nabla\varphi(\bm X_u)^\top\bm\gamma(u,\omega,\bm X_u)\d \bm B_u,\quad t\in[0,T]
\end{equation*} 
is a $(\G_T\vee\F^{\bm X}_t)_{t\in[0,T]}$-martingale. We complete the proof by following Step 4 in the proof of Proposition \ref{Thm_SPDEtoSDE_smooth}. \qed

\begin{proposition}[Theorem \ref{Thm_SPDEtoSDE_general}, locally bounded case] \label{Thm_SPDEtoSDE_locally bdd}
In the setting of Theorem \ref{Thm_SPDEtoSDE_general}, if 
\begin{equation} \label{loc_bnd}
\int_{0}^{T} \sup_{\bm x\in A} |\bm b(t,\omega,\bm x)|^{p}
+\sup_{\bm x\in A} |\bm a(t,\omega,\bm x)|^{p}\,\d t<\infty\;\;\text{a.s.},
\end{equation}
for each bounded Borel set $A\subset\R^d$, then the conclusions of Theorem  \ref{Thm_SPDEtoSDE_general} hold.
\end{proposition}

\noindent\textbf{Proof.} Our proof strategy is similar to that of \cite[Section A.4, case of locally bounded coefficients]{trevisan}. This time, we approximate the coefficients with the help of $C^\infty(\R^d)$ cutoff functions $\chi_R:\,\R^d\to[0,1]$ which satisfy $\chi_R(\bm x)=1$ for $|\bm x|<R$, $\chi_R(\bm x)=0$ for $|\bm x|>2R$, and $|\nabla \chi_R|\leq 2R^{-1}$, $|\nabla^2 \chi_R|\leq 4R^{-2}$ everywhere. 
Define $\bm\pi^R=(\pi^R_1,\ldots,\pi^R_d):\R^d\to\R^d$ by $\bm\pi^R(\bm x):=\chi_R(\bm x)\,\bm x$, so that $\bm\pi^R(\bm x)\to\bm x$ as $R\to\infty$. We note for later use that
\begin{align}
\sup_{R > 0}\|D\bm\pi^R\|_{C^1_b(\R^d)} < \infty. \label{pf:piRbound1}
\end{align}
Here we write $D\bm\pi^R$ for the Jacobian matrix of $\bm\pi^R$, so that its $(i,j)$ entry is $\partial_j \pi^R_i$.

\medskip

\noindent\textbf{Step 1 (Approximation).} We let $\bm\pi^R(\bm x):=\chi_R(\bm x)\,\bm x$ and introduce the pushforward measures $\mu^R_t$ defined by 
\begin{equation*}
\langle\mu^R_t,\varphi\rangle=\langle\mu_t,\varphi\circ\bm\pi^R\rangle,\quad\varphi\in C_c^\infty(\R^d).
\end{equation*}
A straightforward computation yields, for  $\varphi \in C_c^\infty(\R^d)$, 
\begin{align}
& \d\langle \mu^R_t, \varphi\rangle = \sum_{i=1}^d  \Big\langle\mu_t,(\partial_i\varphi)(\bm\pi^R(\cdot)) \,\Big((\nabla \pi^R_i)^\top\bm{b}(t,\omega,\cdot) + \frac12 \bm{a}(t,\omega,\cdot) : \nabla^2 \pi^R_i \Big) \Big\rangle\d t \nonumber \\
&\qquad\qquad\;\;\; + \frac12 \sum_{i,j=1}^d \Big\langle \mu_t, (\nabla\pi^R_i)^\top\bm a(t,\omega,\cdot)(\nabla\pi^R_j)\,(\partial_{ij}\varphi)\big(\bm \pi^R(\cdot)\big)\Big\rangle \d t \nonumber \\
&\qquad\qquad\;\;\; +\big\langle\mu_t,(\nabla\varphi)\big(\bm \pi^R(\cdot)\big)^\top ( D\bm\pi^R)\,\bm\gamma(t,\omega,\cdot)\big\rangle\d\bm B_t.  \label{def_eq_R0}
\end{align}
Define new coefficients $(\bm{b}^R,\bm{a}^R,\bm{\gamma}^R)$ by
\begin{align}
& b^R_i(t,\omega,\cdot):=\E_{\mu_t(\omega)}\big[ \nabla \pi^R_i(\bm X)^\top\bm{b}(t,\omega,\bm X) + \frac12 \bm{a}(t,\omega,\bm X) : \nabla^2 \pi^R_i(\bm X) \,\big| \, \bm\pi^R(\bm X)=\cdot\big], \nonumber \\
& a^R_{ij}(t,\omega,\cdot):=\E_{\mu_t(\omega)}\big[\nabla\pi^R_i(\bm X)^\top\bm a(t,\omega,\bm X)\nabla\pi^R_j(\bm X) \, \big| \, \bm\pi^R(\bm X)=\cdot\big], \label{def_coefficient_R} \\
& \bm \gamma^R(t,\omega,\cdot):=\E_{\mu_t(\omega)}\big[(D\bm\pi^R)(\bm X)\,\bm\gamma(t,\omega,\bm X) \,\big|\, \bm\pi^R(\bm X)=\cdot\big], \nonumber
\end{align}
where the notation is understood to mean, for each fixed $(t,\omega)$, that we are taking expectations of functions of a random variable $\bm X \sim \mu_t(\omega)$, conditionally on the value of $\bm\pi^R(\bm X)$.
Then \eqref{def_eq_R0} becomes
\begin{equation*}
\d\langle \mu^R_t, \varphi\rangle = \big\langle \mu^R_t, \gen^{\bm b^R,\bm a^R}_{t,\omega}\varphi\big\rangle \d t + \big\langle \mu^R_t,(\nabla\varphi)^\top \bm \gamma^R(t,\omega,\cdot)\big\rangle \d\bm B_t, \quad t \in [0,T],\quad \varphi \in C_c^\infty(\R^d).
\end{equation*}

\smallskip

Notice that the new coefficients $\bm b^R$ and $\bm a^R$ obey 
\begin{equation*}
\begin{split}
&\sup_{\bm x\in\R^d}|\bm b^R(t,\omega,\bm x)|\leq \|D\bm\pi^R\|_{C_b^1(\R^d)}\,\sup_{|\bm x|\leq 2R}\big(|\bm b(t,\omega,\bm x)|+|\bm a(t,\omega,\bm x)|\big), \\
& \sup_{\bm x\in\R^d}|\bm a^R(t,\omega,\bm x)|\leq \|D\bm\pi^R\|^2_{C_b(\R^d)}\sup_{|\bm x|\leq 2R}|\bm a(t,\omega,\bm x)|.
\end{split}
\end{equation*}
Together with \eqref{loc_bnd}, these show that \eqref{prop3.2int} holds for $\bm b^R$, $\bm a^R$. Moreover, Jensen's inequality and \eqref{pf:piRbound1} yield \eqref{integ} for $\mu^R$, $\bm b^R$, $\bm a^R$. Lastly, we compute
\begin{equation*}
\begin{split}
&\,\big(\bm a^{R}-\bm\gamma^{R}({\bm\gamma^{R}})^\top\big)(t,\omega,\cdot) \\
&=\E_{\mu_t(\omega)}\big[(D\bm\pi^R)(\bm X)(\bm a-\bm \gamma\bm\gamma^\top)(t,\omega,\bm X)(D\bm\pi^R)^\top(\bm X)\big|\bm\pi^R(\bm X)=\cdot\big] \\
&\;\;\;\,+\E_{\mu_t(\omega)}\big[\big((D\bm\pi^R)(\bm X)\,\bm \gamma(t,\omega,\bm X)\big)\big((D\bm\pi^R)(\bm X)\,\bm\gamma(t,\omega,\bm X)\big)^\top\big|\bm\pi^R(\bm X)=\cdot\big] \\
&\;\;\;\,-\E_{\mu_t(\omega)}\big[(D\bm\pi^R)(\bm X)\,\bm \gamma(t,\omega,\bm X)\big|\bm\pi^R(\bm X)=\cdot\big]\,\E_{\mu_t(\omega)}\big[(D\bm\pi^R)(\bm X)\,\bm \gamma(t,\omega,\bm X)|\bm\pi^R(\bm X)=\cdot\big]^\top,
\end{split}
\end{equation*}
where the latter two terms give the conditional covariance matrix of $(D\bm\pi^R)(\bm X)\,\bm \gamma(t,\omega,\bm X)$. Hence, $\bm a^{R}-\bm\gamma^{R}({\bm\gamma^{R}})^\top$ is positive definite and symmetric, and we can let $\bm\sigma^R$ be its symmetric square root. All in all, $\mu^R$ and the coefficients $\bm b^R$, $\bm a^R$, $\bm\gamma^R$ comply with the conditions of Proposition \ref{Thm_SPDEtoSDE_bdd}, and the superposition principle applies for each $R$.
Transferring to the canonical space $\widetilde\Omega$ introduced in  \eqref{canonicalspace}, we find $\widetilde\PP^R\in\P(\widetilde\Omega;\PP)$ for $R>0$ (see Definition \ref{def:canonical}) such that, for each $R$, 
\begin{equation*}
\d\bm{X}_t = \bm{b}^{\varrho_n}(t,\omega,\bm{X}_t)\d t + \bm{\sigma}^{\varrho_n}(t,\omega,\bm{X}_t)\d\bm{W}_t + \bm{\gamma}^{\varrho_n}(t,\omega,\bm{X}_t)\d\bm{B}_t\quad \widetilde\PP^R\text{-a.s.,}
\end{equation*}
with $\mu^R_t=\L_{\widetilde\PP^R}(\bm{X}_t\,|\,\G_T)=\L_{\widetilde\PP^R}(\bm{X}_t\,|\,\G_t)$ a.s.\ for each $t\in[0,T]$, and with $\F^{\bm{X}}_t \indep \F^{\bm{W}}_T \vee \G_T \,|\,\F^{\bm{W}}_t \vee \G_t$ under $\widetilde\PP^R$ for each $t \in [0,T]$.

\medskip

\noindent\textbf{Step 2 (Tightness).} To apply Proposition \ref{Prop_tightness} it suffices to provide a uniform bound on
\begin{equation}\label{tight_bnd}
\E_{\widetilde\PP^R}\bigg[\int_0^T\big|\bm b^{R}(t,\omega,\bm X_t)\big|^p+\big|\bm a^{R}(t,\omega,\bm X_t)\big|^p\d t\bigg].
\end{equation}
By Fubini's theorem, the tower rule, the definition of $\mu^R$, and Jensen's inequality,
\begin{equation*}
\begin{split}
&\;\E_{\widetilde\PP^R}\bigg[\int_0^T\big|\bm b^{R}(t,\omega,\bm X_t)\big|^p
+\big|\bm a^{R}(t,\omega,\bm X_t)\big|^p\d t\bigg] \\
&\quad=\E\bigg[\int_0^T\big\langle\mu_t^R,\big|\bm b^{R}(t,\omega,\cdot)\big|^p
+\big|\bm a^{R}(t,\omega,\cdot)\big|^p\big\rangle\d t\bigg] \\
&\quad=\E\bigg[\int_0^T\big\langle\mu_t,\big|\bm b^{R}\big(t,\omega,\bm\pi^R(\cdot)\big)\big|^p
+\big|\bm a^{R}\big(t,\omega,\bm\pi^R(\cdot)\big)\big|^p\big\rangle\d t\bigg] \\
&\quad\leq C_p\max\Big(\|D\bm \pi^R\|_{C_b^1(\R^d)}^p,\|D\bm \pi^R\|_{C_b(\R^d)}^{2p}\Big)\,\E\bigg[\int_0^T\big\langle\mu_t,\big|\bm b(t,\omega,\cdot)\big|^p+\big|\bm a(t,\omega,\cdot)\big|^p\big\rangle\d t\bigg].
\end{split}
\end{equation*}
Thanks to \eqref{integ} and since $\|D\bm \pi^R\|_{C_b^1(\R^d)}$ remains uniformly bounded as $R\in\N$ varies, the expectation in \eqref{tight_bnd} is uniformly bounded.
Applying Proposition \ref{Prop_tightness} we obtain a limit point $\widetilde\PP \in \P(\widetilde\Omega;\PP)$ with respect to stable convergence. Once again we relabel the subsequence so that $(\widetilde\PP^R)_{R > 0}$  converges to $\widetilde\PP$.
Recall from Lemma \ref{lem:closed} that $\P(\widetilde\Omega;\PP)$ is closed in the stable topology, and in particular $\F^{\bm{X}}_t \indep \F^{\bm{W}}_T \vee \G_T\,|\,\F^{\bm{W}}_t \vee \G_t$ for each $t \in [0,T]$.
By Lemma \ref{le:compatibility} (with $\HH=\FF^{\bm{X}}$), this implies that $\F^{\bm{X}}_t \indep \G_T\,|\,\G_t$ for each $t \in [0,T]$. Lemma \ref{le:stable-condlaw} ensures that also $\mu_t=\L_{\widetilde\PP}(\bm{X}_t\,|\,\G_T)=\L_{\widetilde\PP}(\bm{X}_t\,|\,\G_t)$ a.s.\ for each $t$.

\medskip

\noindent\textbf{Step 3 (Limit).} With $\varphi$, $g$, $s$, $t$, $h^s$ as in Step 3 of the proof of Proposition \ref{Thm_SPDEtoSDE_smooth}, we have: 
\begin{equation*}
\E_{\widetilde\PP^R}\bigg[g h^s(\bm X)\bigg(\!\varphi(\bm X_t)-\varphi(\bm X_s)-\int_s^t\!\!\gen^{\bm b^R,\bm a^R}_{u,\omega}\varphi(\bm X_u)\d u-\int_s^t\!\!\nabla\varphi(\bm X_u)^\top\bm\gamma^R(u,\omega,\bm X_u)\d \bm B_u\!\bigg)\bigg]\!\!=\!0.
\end{equation*}
As in Step 3 of the proof of Proposition \ref{Thm_SPDEtoSDE_bdd}, we define below a triple $(\widetilde{\bm{b}},\widetilde{\bm{a}},\widetilde{\bm{\gamma}}):\,[0,T]\times\Omega\times\R^d\to\R^d\times\R^{d\times d}\times\R^{d\times d}$ measurable with respect to the product of the $\mathbb G$-progressive $\sigma$-algebra on $[0,T]\times \Omega$ and the Borel $\sigma$-algebra on $\R^d$ such that all $(\widetilde{\bm{b}},\widetilde{\bm{a}},\widetilde{\bm{\gamma}})(t,\omega,\cdot)$ are continuous and compactly supported in $\R^d$, with their supremum norms being bounded uniformly by a constant $M<\infty$. For these, it holds: 
\begin{equation}\label{equa6}
\begin{split}
&\,\bigg|\E_{\widetilde\PP^R}\bigg[gh^s(\bm X)\bigg(\varphi(\bm X_t)-\varphi(\bm X_s)-\int_s^t \gen^{\widetilde{\bm b},\widetilde{\bm a}}_{u,\omega}\varphi(\bm X_u)\d u-\int_s^t \nabla\varphi(\bm X_u)^\top\widetilde{\bm\gamma}(u,\omega,\bm X_u)\d \bm B_u\bigg)\bigg]\bigg| \\
&=\!\bigg|\E_{\widetilde\PP^R}\bigg[gh^s(\bm X)\bigg(\!\int_s^t\!\! \big(\gen^{\bm b^R,\bm a^R}_{u,\omega}-\gen^{\widetilde{\bm b},\widetilde{\bm a}}_{u,\omega}\big)\varphi(\bm X_u)\d u
+\!\int_s^t\!\!\nabla\varphi(\bm X_u)^\top\big({\bm\gamma}^R-\widetilde{\bm\gamma}\big)(u,\omega,\bm X_u)\d \bm B_u\!\bigg)\bigg]\bigg|.
\end{split}
\end{equation}
In addition, we note:
\begin{equation*}
\begin{split}
\bigg|\E_{\widetilde\PP^R}\bigg[gh^s(\bm X)\int_s^t\big(\gen^{\bm b^R,\bm a^R}_{u,\omega}-\gen^{\widetilde{\bm b},\widetilde{\bm a}}_{u,\omega}\big)\varphi(\bm X_u)\d u\bigg]\bigg|
\leq\E_{\widetilde\PP^R}\bigg[\int_s^t \big|\big(\gen^{\bm b^R,\bm a^R}_{u,\omega}-\gen^{\widetilde{\bm b},\widetilde{\bm a}}_{u,\omega}\big)\varphi(\bm X_u)\big|\d u\bigg]\; \\
=\E\bigg[\int_s^t\big\langle\mu_u^R,\big| (\bm b^R-\widetilde{\bm b})(u,\omega,\cdot)\cdot\nabla\varphi+\frac{1}{2}\,(\bm a^R-\widetilde{\bm a})(u,\omega,\cdot): \nabla^2\varphi\big|\big\rangle\d u\bigg].
\end{split}
\end{equation*}
By plugging in the definitions of $\mu^R$, $\bm b^R$, $\bm a^R$ and using Jensen's inequality we bound this further:
\begin{equation*}
\begin{split}
&\big\langle\mu_u^R,\big| ( \bm b^R-\widetilde{\bm b})(u,\omega,\cdot)\cdot \nabla\varphi\big|\big\rangle  \\
& \ \ = \big\langle\mu_u,\big| ( \bm b^R-\widetilde{\bm b})(u,\omega,\bm\pi^R(\cdot)) \cdot \nabla\varphi (\bm\pi^R(\cdot))\big|\big\rangle \\
& \ \ \leq \|\varphi\|_{C_b^1(\R^d)}\, \Big\langle\mu_u,\Big(\sum_{i=1}^d\big|(\nabla \pi^R_i)^\top\bm{b}(u,\omega,\cdot) + \frac12 \bm{a}(u,\omega,\cdot) : \nabla^2 \pi^R_i  - \widetilde{b}_i(u,\omega,\bm\pi^R(\cdot))\big|^2 \Big)^{1/2}\Big\rangle
\end{split}
\end{equation*}
for $(u,\omega) \in [0,T] \times \Omega$, and, similarly,
\begin{equation*}
\begin{split}
& \big\langle\mu_u^R,\big| (\bm a^R-\widetilde{\bm a})(u,\omega,\cdot): \nabla^2\varphi\big|\big\rangle \qquad\qquad\quad\;\;\, \\
& \ \ \leq \|\varphi\|_{C_b^2(\R^d)}\, \big\langle\mu_u,\big|(D\bm\pi^R)\,\bm a(u,\omega,\cdot)\,(D\bm\pi^R)^\top-\widetilde{\bm a}(u,\omega,\bm\pi^R(\cdot))\big|\big\rangle .
\qquad 
\end{split}
\end{equation*}
From the definition of $\bm\pi^R$ it is clear that $D\bm\pi^R$ converges pointwise to the identity matrix, and all second derivatives of $\bm\pi^R$ converge pointwise to zero.
Recalling \eqref{pf:piRbound1}, \eqref{integ}, and that $(\widetilde{\bm b},\widetilde{\bm a})$ are bounded and continuous, we may pass to the $R\to\infty$ limit by means of the dominated convergence theorem and get
\begin{equation*}
\begin{split}
&\limsup_{R\to\infty}\bigg|\E_{\widetilde\PP^R}\bigg[gh^s(\bm X)\int_s^t\big(\gen^{\bm b^R,\bm a^R}_{u,\omega}-\gen^{\widetilde{\bm b},\widetilde{\bm a}}_{u,\omega}\big)\varphi(\bm X_u)\d u\bigg]\bigg| \\
&\qquad \leq \|\varphi\|_{C_b^2(\R^d)}\,\E\bigg[\int_s^t\big\langle\mu_u,\big|\gen^{\bm b-\widetilde{\bm b},\bm a-\widetilde{\bm a}}_{u,\omega}\big|\big\rangle\d u\bigg],
\end{split}
\end{equation*}
with the shorthand notation introduced in \eqref{Lshorthand}.

\smallskip

For the stochastic integral term, a similar argument upon an application of the tower rule, Jensen's inequality, and It\^o isometry yields
\begin{equation*}
\begin{split}
&\,\bigg|\E_{\widetilde\PP^R}\bigg[gh^s(\bm X)\int_s^t \nabla\varphi(\bm X_u)^\top({\bm\gamma}^R-\widetilde{\bm\gamma})(u,\omega,\bm X_u)\d \bm B_u\bigg]\bigg| \\
&\qquad \leq \E\bigg[\int_s^t\big\langle\mu_u^R,\big|(\nabla\varphi)^\top(\bm\gamma^R-\widetilde{\bm \gamma})(u,\omega,\cdot)\big|\big\rangle^2\d u\bigg]^{1/2} \\
&\qquad \leq \|\varphi\|_{C_b^1(\R^d)}\,\E\bigg[\int_s^t \big\langle\mu_u,\big|(D\bm\pi^R)\,\bm\gamma(u,\omega,\cdot)-\widetilde{\bm\gamma}(u,\omega,\bm\pi^R(\cdot))\big|\big\rangle^2\d u\bigg]^{1/2},
\end{split}
\end{equation*}
so that by the dominated convergence theorem,
\begin{equation*}
\begin{split}
&\limsup_{R\to\infty}\bigg|\E_{\widetilde\PP^R}\bigg[gh^s(\bm X)\int_s^t \nabla\varphi(\bm X_u)^\top({\bm\gamma}^R-\widetilde{\bm\gamma})(u,\omega,\bm X_u)\d \bm B_u\bigg]\bigg| \\
&\qquad \le \|\varphi\|_{C_b^1(\R^d)}\,\E\bigg[\int_s^t\big\langle\mu_u,\big|(\bm\gamma-\widetilde{\bm\gamma})(u,\omega,\cdot)\big|\big\rangle^2\d u\bigg]^{1/2}.
\end{split}
\end{equation*}
We now take $R\to \infty$ in \eqref{equa6} and use Proposition \ref{prop_mtgconverge} to deduce that
\begin{equation*}
\begin{split}
&\,\bigg|\E_{\widetilde\PP}\bigg[gh^s(\bm X)\bigg(\varphi(\bm X_t)-\varphi(\bm X_s)-\int_s^t \gen^{\widetilde{\bm b},\widetilde{\bm a}}_{u,\omega}\varphi(\bm X_u)\d u-\int_s^t \nabla\varphi(\bm X_u)^\top\widetilde{\bm\gamma}(u,\omega,\bm X_u)\d \bm B_u\bigg)\bigg]\bigg| \\
&\leq \|\varphi\|_{C_b^2(\R^d)} \bigg(\E\bigg[\int_s^t\big\langle\mu_u, \big|\gen^{\bm b-\widetilde{\bm b},\bm a-\widetilde{\bm a}}_{u,\omega}\big|\big\rangle\d u\bigg]+\E\bigg[\int_s^t\big\langle\mu_u,\big|(\bm\gamma-\widetilde{\bm\gamma})(u,\omega,\cdot)\big|\big\rangle^2\d u\bigg]^{1/2}\bigg).
\end{split}
\end{equation*}
Arguing as in the paragraph between \eqref{equa5} and \eqref{equa_consist} we change $(\widetilde{\bm b},\widetilde{\bm a},\widetilde{\bm\gamma})$ to $(\bm b,\bm a,\bm\gamma)$: 
\begin{equation}\label{equa_consistag'}
\begin{split}
&\,\bigg|\E_{\widetilde\PP}\bigg[gh^s(\bm X)\bigg(\varphi(\bm X_t)-\varphi(\bm X_s)-\int_s^t \gen^{\bm b,\bm a}_{u,\omega}\varphi(\bm X_u)\d u-\int_s^t \nabla\varphi(\bm X_u)^\top\bm\gamma(u,\omega,\bm X_u)\d \bm B_u\bigg)\bigg]\bigg| \\
&\leq 2\|\varphi\|_{C_b^2(\R^d)} \bigg(\E\bigg[\int_s^t\big\langle\mu_u, \big|\gen^{\bm b-\widetilde{\bm b},\bm a-\widetilde{\bm a}}_{u,\omega}\big|\big\rangle\d u\bigg]+\E\bigg[\int_s^t\big\langle\mu_u,\big|(\bm\gamma-\widetilde{\bm\gamma})(u,\omega,\cdot)\big|\big\rangle^2\d u\bigg]^{1/2}\bigg).
\end{split}
\end{equation}

\smallskip

It remains to construct $(\widetilde{\bm b},\,\widetilde{\bm a},\,\widetilde{\bm \gamma})$. We pick an approximation level $\varepsilon>0$ and claim that we can find some $M<\infty$ such that $b_{i}\,\mathbf{1}_{\{|b_{i}|\le M\}}$ $a_{ij}\,\mathbf{1}_{\{|a_{ij}|\le M\}}$, and $\gamma_{ij}\,\mathbf{1}_{\{|\gamma_{ij}|\le M\}}$ are $\varepsilon$-good approximations of $b_i,\,a_{ij},\, \gamma_{ij}$ in the distance on the right-hand side of \eqref{equa_consistag'}. Indeed,
\begin{equation*}
\begin{split}
&\;\E\bigg[\int_s^t\big\langle\mu_u,\big|b_i(u,\omega,\cdot)-b_i(u,\omega,\cdot)\,\mathbf{1}_{\{|b_i(u,\omega,\cdot)|\le M\}}\big|\big\rangle\d u\bigg] \\
&\qquad =\E\bigg[\int_s^t\big\langle\mu_u,\big|b_i(u,\omega,\cdot)\,\mathbf{1}_{\{|b_i(u,\omega,\cdot)|> M\}}\big|\big\rangle\d u\bigg]
\leq M^{1-p}\,\E\bigg[\int_s^t\big\langle\mu_u,\big|b_i(u,\omega,\cdot)\big|^p\big\rangle\d u\bigg],
\end{split}
\end{equation*}
the corresponding expression for $a_{ij}$, and 
\begin{equation*}
\begin{split}
&\;\E\bigg[\int_s^t\big\langle\mu_u,\big|\gamma_{ij}(u,\omega,\cdot)-\gamma_{ij}(u,\omega,\cdot)\,\mathbf{1}_{\{|\gamma_{ij}(u,\omega,\cdot)|\le M\}}\big|\big\rangle^2\d u\bigg] \\
&\leq\E\bigg[\int_s^t\big\langle\mu_u,\big|\gamma_{ij}(u,\omega,\cdot)\,\mathbf{1}_{\{|\gamma_{ij}(u,\omega,\cdot)|>M\}}\big|^2\big\rangle\d u\bigg]
\leq M^{2-2p}\,\E\bigg[\int_s^t\big\langle\mu_u,\big|\gamma_{ij}(u,\omega,\cdot)\big|^{2p}\big\rangle\d u\bigg]
\end{split}
\end{equation*}
can all be made arbitrarily small by choosing a large enough $M<\infty$. We now define $(\widetilde{\bm b},\,\widetilde{\bm a},\,\widetilde{\bm \gamma})$ by applying the procedure in Step 3 of the proof of Proposition \ref{Thm_SPDEtoSDE_bdd} to $b_i\,\mathbf{1}_{\{|b_i|\le M\}}$, $a_{ij}\,\mathbf{1}_{\{|a_{ij}|\le M\}}$, $\gamma_{ij}\,\mathbf{1}_{\{|\gamma_{ij}| \le M\}}$, thus making the right-hand side of \eqref{equa_consistag'} arbitrarily small. Consequently,
\begin{equation*}
\varphi(\bm X_t)-\varphi(\bm X_0)-\int_0^t \gen^{\bm b,\bm a}_{u,\omega}\varphi(\bm X_u)\d u-\int_0^t \nabla\varphi(\bm X_u)^\top\bm\gamma(u,\omega,\bm X_u)\d \bm B_u,\quad t\in[0,T]
\end{equation*}
is a $(\G_T\vee\F^{\bm X}_t)_{t\in[0,T]}$-martingale. We complete the proof exactly as in Step 4 of the proof of Proposition \ref{Thm_SPDEtoSDE_smooth}. \qed

%

\bigskip

\noindent\textbf{Proof of Theorem \ref{Thm_SPDEtoSDE_general}.} We are finally ready to prove the general case, Theorem \ref{Thm_SPDEtoSDE_general}. Our strategy is similar that of \cite[Section A.4, general case]{trevisan}. We use the same convolution construction as in the proof of Proposition \ref{Thm_SPDEtoSDE_bdd} to reduce the general case to the locally bounded case covered by Proposition \ref{Thm_SPDEtoSDE_locally bdd}.

\medskip

\noindent\textbf{Step 1 (Approximation).} We define $\mu^n$, $\bm b^{\varrho_n}$, $\bm a^{\varrho_n}$, and $\bm \gamma^{\varrho_n}$ according to the formulas in \eqref{def_coefficients_rhon} and the line preceding it, with a strictly positive $\varrho$ as there. Next, we check that $\bm b^{\varrho_n}$, $\bm a^{\varrho_n}$, $\bm \gamma^{\varrho_n}$ satisfy the conditions in Proposition \ref{Thm_SPDEtoSDE_locally bdd}. By Jensen's inequality, we have: 
\begin{equation*}
\begin{split}
\big|\bm b^{\varrho_n}(t,\omega,\bm x)\big|
= \bigg|\frac{\int_{\R^d}\varrho_n(\bm x-\bm y) \, \bm b(t,\omega,\bm y) \,\mu_t(\mathrm{d}\bm y)}{\int_{\R^d}\varrho_n(\bm x-\bm y)\,\mu_t(\mathrm{d}\bm y)}\bigg| 
\leq\|\varrho_n\|_{C_b(\R^d)}\frac{\int_{\R^d} |\bm b(t,\omega,\bm y)|\,\mu_t(\mathrm{d}\bm y)}{\int_{\R^d}\varrho_n(\bm x-\bm y)\,\mu_t(\mathrm{d}\bm y)}.
\end{split}
\end{equation*}
Since $\int_{\R^d} \varrho_n(\bm x-\bm y)\,\mu_t(\mathrm{d}\bm y)$ is a strictly positive continuous function of $(t,\bm x)$ a.s., it admits a lower bound $c_A(\omega)>0$ on $[0,T]\times A$ for any fixed bounded Borel $A\subset\R^d$. Hence, by Jensen's inequality,
\begin{equation*}
\int_0^T \sup_{\bm x\in A}\big|\bm b^{\varrho_n}(t,\omega,\bm x)\big|^p\d t
\leq  c_A(\omega)^{-p}\,\|\varrho_n\|_{C_b(\R^d)}^p
\int_0^T\int_{\R^d} |\bm b(t,\omega,\bm y)|^p\,\mu_t(\mathrm{d}\bm y)\d t.
\end{equation*}
Carrying out the same estimates for $\bm a^{\varrho_n}$ and recalling \eqref{integ} we deduce that  
\begin{equation*}
	\int_{0}^{T}\sup_{\bm x\in A}\big|\bm b^{\varrho_n}(t,\omega,\bm x)\big|^{p}
	+\sup_{\bm x\in A}\big|\bm a^{\varrho_n}(t,\omega,\bm x)\big|^{p}\d t<\infty\;\;\text{a.s.},
\end{equation*}
for each $n$.
The other conditions in Proposition \ref{Thm_SPDEtoSDE_locally bdd} can be verified for $\mu^n$, $\bm b^{\varrho_n}$, $\bm a^{\varrho_n}$, and $\bm \gamma^{\varrho_n}$ exactly as in Step 1 of the proof of Proposition \ref{Thm_SPDEtoSDE_bdd}. Therefore, the superposition principle of Proposition \ref{Thm_SPDEtoSDE_bdd} allows us to find probability measures $\widetilde\PP^n \in \P(\widetilde\Omega;\PP)$ for $n\in\N$ such that, under $\widetilde\PP^n$, 
\begin{equation*}
\d\bm{X}_t = \bm{b}^{\varrho_n}(t,\omega,\bm{X}_t)\d t + \bm{\sigma}^{\varrho_n}(t,\omega,\bm{X}_t)\d\bm{W}_t + \bm{\gamma}^{\varrho_n}(t,\omega,\bm{X}_t)\d\bm{B}_t\quad \widetilde\PP^n\text{-a.s.,}
\end{equation*}
with $\mu^n_t=\L_{\widetilde\PP^n}(\bm{X}_t\,|\,\G_T)=\L_{\widetilde\PP^n}(\bm{X}_t\,|\,\G_t)$ a.s.\ for each $t\in[0,T]$, and with $\F^{\bm{X}}_t \indep \F^{\bm{W}}_T \vee \G_T \,|\,\F^{\bm{W}}_t \vee \G_t$ under $\widetilde\PP^n$ for each $t \in [0,T]$.

\medskip

\noindent\textbf{Step 2 (Tightness).} Step 2 in the proof of Proposition  \ref{Thm_SPDEtoSDE_bdd} can be repeated literally. 

\medskip

\noindent\textbf{Step 3 (Limit).} We derive the inequality \eqref{equa_consist} as in the proof of Proposition \ref{Thm_SPDEtoSDE_bdd} and pick the coefficients $\widetilde{\bm b},\widetilde{\bm a}, \widetilde{\bm\gamma}$ therein as in Step 3 of the proof of Proposition  \ref{Thm_SPDEtoSDE_locally bdd}.~This renders
\begin{equation*}
\varphi(\bm X_t)-\varphi(\bm X_0)-\int_0^t \gen^{\bm b,\bm a}_{u,\omega}\varphi(\bm X_u)\d u-\int_0^t \nabla\varphi(\bm X_u)^\top\bm\gamma(u,\omega,\bm X_u)\d \bm B_u,\quad t\in[0,T]
\end{equation*}
a $(\G_T\vee\F^{\bm X}_t)_{t\in[0,T]}$-martingale, for $\varphi \in C^\infty_c(\R^d)$. We then complete the proof exactly as in Step 4 of the proof of Proposition \ref{Thm_SPDEtoSDE_smooth}. \qed

\section{Superposition from Fokker-Planck equation on $\P(\R^d)$ to SPDE} \label{sec_proof_PDEtoSPDE}

\noindent\textbf{Proof of Theorem \ref{Thm_PDEtoSPDE}.} The main idea of the proof is to recast $(P_t)_{t\in[0,T]}$ as probability measures on $\R^\infty$ and to use the superposition principle of \cite[Theorem 7.1]{trevisan}. 

\medskip

\noindent\textbf{Step 1 (Reduction from $\P(\R^d)$ to $\R^\infty$).} We start by establishing the following lemma which yields a countable collection of suitable test functions. 
\begin{lemma}\label{Lem_sep}
There exists a set $\{\varphi_n:\,n\in\N\}\subset C_c^\infty(\R^d)$ such that for all $\varphi\in C_c^\infty(\R^d)$ one can find a sequence $(n_k)_{k\in\N}$ such that $(\varphi_{n_k},\nabla\varphi_{n_k},\nabla^2 \varphi_{n_k})_{k\in\N}$ converges uniformly to $(\varphi,\nabla\varphi,\nabla^2 \varphi)$.
\end{lemma}

\noindent\textbf{Proof.} It suffices to show that the space $C_c^\infty(\R^d)$ is separable under the norm
\begin{equation*}
\|\varphi\|_*:=\max_{|\alpha|\leq 2}\,\|D^\alpha \varphi\|_{C_b(\R^d)},
\end{equation*}
where we have adopted multi-index notation.
We first prove that the space $C_c^\infty(\R^d)$ is separable under the uniform norm. By the Stone-Weierstrass theorem, $C([-K,K]^d)$ is separable under the uniform norm for every $K\in\N$, e.g., the polynomials with rational coefficients form a countable dense set. Therefore, as a subspace of a separable space, $C_c^\infty((-K,K)^d)$ is also separable under the uniform norm for every $K\in\N$. Suppose $S_K$ is a countable dense subset of $C_c^\infty((-K,K)^d)$. Then $\bigcup_{K\in\N} S_K$ forms a countable dense subset of $C_c^\infty(\R^d)$ under the uniform norm.
	
\medskip	
	
Further, for each $m\in\N$, the space $C_c^\infty(\R^d)^{m}$ is separable under the norm
	\begin{equation*}
		\|\bm\varphi\|_{\infty}:=\max_{1\le i\le m}\, \|\varphi_i\|_{C_b(\R^d)},\quad \bm\varphi=(\varphi_1,\ldots,\varphi_m)\in C_c^\infty(\R^d)^{m}.
	\end{equation*}
	Consider the mapping $\M:C_c^\infty(\R^d)\to C_c^\infty(\R^d) ^{d^2+d+1},\, \varphi\to(\varphi,\nabla\varphi,\nabla^2\varphi)$. Then, $\|\varphi\|_*=\|\M(\varphi)\|_{\infty}$, so $\M$ is an isometry from $(C_c^\infty(\R^d),\|\cdot\|_*)$ to $(\M(C_c^\infty(\R^d)),\|\cdot\|_\infty)$. Since $(\M(C_c^\infty(\R^d)),\|\cdot\|_\infty)$ is a subspace of a separable space, it is itself separable, which proves that $(C_c^\infty(\R^d),\|\cdot\|_*)$ is separable. \qed

\medskip

Henceforth, we fix $\{\varphi_n : n \in \N\} \subset C^\infty_c(\R^d)$ as in Lemma \ref{Lem_sep} and define the mapping $\SS:\,\P(\R^d)\to\R^\infty,\, m\mapsto(\langle  m, \varphi_n\rangle)_{n\in\N}$. It is straightforward to see that $\mathcal{S}$ is an injection. The next step is to lift the PDE \eqref{FPE} to $\P(\R^\infty)$.

\medskip

\noindent\textbf{Step 2 (Identifying the PDE in $\P(\mathbb{R}^\infty)$).} Let $Q_t=P_t \circ \SS^{-1}$, i.e., $Q_t$ is the law of $\SS(\mu_t)$ when $\L(\mu_t)=P_t$.
Let us write $C^\infty_c(\R^\infty)$ for the set of smooth cylindrical functions on $\R^\infty$, that is, functions of the form $f(\bm{z})=\widetilde{f}(z_1,\ldots,z_k)$ for some $k \in \N$ and $\widetilde{f} \in C^\infty_c(\R^k)$.
Applying \eqref{FPE} with $\bm\varphi=(\varphi_n)_{n\in \N}$ we find for every $f\in C^\infty_c(\R^\infty)$ and $t \in [0,T]$:
\begin{align}
&\int_{\R^\infty} f(\bm z)\, (Q_t-Q_0)(\! \d\bm z) \nonumber \\
&\qquad = \int_0^t \int_{\R^\infty} \Bigg[ \sum_{i=1}^\infty \partial_if(\bm z)\,\beta_i(s,\bm z) + \frac12\sum_{i,j=1}^\infty \partial_{ij} f(\bm z)\,\alpha_{ij}(s,\bm z) \Bigg]\,Q_s(\! \d\bm z) \d s, \label{pf:Ptequa1}
\end{align}
where $\bm \beta : [0,T] \times \R^\infty \to \R^\infty$ and $\bm \alpha : [0,T] \times \R^\infty \to \R^\infty \times \R^\infty$ are given by
\begin{align*}
\beta_i(t, \bm z) &= \big\langle \SS^{-1}(\bm z), \gen^{\bm b, \bm a}_{t,\SS^{-1}(\bm z)}\varphi_i\big\rangle, \\
\alpha_{ij}(t, \bm z) &= \big\langle \SS^{-1}(\bm z), (\nabla \varphi_i)^\top \bm{\gamma}(t,\SS^{-1}(\bm z),\cdot)\big\rangle \cdot \big\langle \SS^{-1}(\bm z), (\nabla \varphi_j)^\top \bm{\gamma}(t,\SS^{-1}(\bm z),\cdot)\big\rangle.
\end{align*}
Note that the summations in \eqref{pf:Ptequa1} are in fact finite sums, since $f$ is cylindrical, and that we only need to define $\bm \beta(t,\cdot)$ and $\bm \alpha(t,\cdot)$ on $\SS(\P(\R^d))$, since $Q_s$ is supported on $\SS(\P(\R^d))$.

\medskip

\noindent\textbf{Step 3 (Using the superposition principle in $\R^\infty$).} We use \cite[Theorem 7.1]{trevisan} upon checking that, for $i,j \in \N$,
\begin{align*}
& \int_0^T \int_{\R^\infty} |\beta_i(t,\bm z)|^p\,Q_t(\!\d \bm z)\d t
=\int_0^T\int_{\P(\R^d)}\big|\big\langle m,\gen^{\bm b,\bm a}_{t, m}\varphi_i\big\rangle\big|^p \,  P_t(\!\d m)\d t \\
&\qquad \leq 2^{p-1}\|\varphi_i\|^p_{C_b^2(\R^d)}\int_{0}^{T}\int_{\P(\R^d)} \big\|\bm b(t, m,\cdot)\big\|^p_{L_p( m)}+\big\|\bm a(t, m,\cdot)\big\|^p_{L_p( m)}\,P_t(\!\d m)\d t<\infty,
\end{align*}
and, similarly,
\begin{align*}
&\,\int_0^T\int_{\R^\infty} |\alpha_{ij}(t,\bm z)|^p \, Q_t(\!\d \bm z)\d t  \\
&\qquad =\int_0^T\int_{\P(\R^d)}\Big|\big\langle  m,(\nabla \varphi_i)^\top \bm\gamma(t, m,\cdot)\big\rangle\cdot\big\langle  m, (\nabla \varphi_j)^\top \bm\gamma(t, m,\cdot)\big\rangle\Big|^p\, P_t(\!\d m)\d t \\
&\qquad \leq \|\varphi_i\|^{p}_{C_b^1(\R^d)}\|\varphi_j\|^{p}_{C_b^1(\R^d)}\int_{0}^{T}\int_{\P(\R^d)}\big\|\bm \gamma(t, m,\cdot)\big\|^{2p}_{L_p( m)}P_t(\!\d  m)\d t<\infty.
\end{align*}
By \cite[Theorem 7.1]{trevisan}, there exists a solution $Q \in \P(C([0,T];\R^\infty))$ of the martingale problem associated to the system of SDEs
\begin{align}\label{MV:reSPDE-CN}
\d Z^i_t&= \big\langle \mathcal{S}^{-1}(\bm Z_t),\gen^{\bm b,\bm a}_{t,\mathcal{S}^{-1}(\bm Z_t)}\varphi_i\big\rangle \d t+\big\langle \mathcal{S}^{-1}(\bm Z_t),(\nabla\varphi_i)^\top \bm\gamma(t,\mathcal{S}^{-1}(\bm Z_t),\cdot)\big\rangle \d\bm B_t,\quad i\in\N 
\end{align}
for which the corresponding marginal flow $(\L(\bm Z_t))_{t\in[0,T]}$ equals $(Q_t)_{t \in [0,T]}$. 

\medskip

\noindent\textbf{Step 4 (Mapping back to $C([0,T];\P(\R^d))$).} We claim that $Q$ gives rise to a weak solution of \eqref{MV:reSPDE-CN}. Indeed, by following \cite[Chapter 5, proof of Proposition 4.6]{karatzas-shreve} the claim can be reduced to a martingale representation theorem in the form of \cite[Chapter 3, Theorem 4.2]{karatzas-shreve}, but in our case with countably infinitely many local martingales. To deal with this discrepancy we choose, for each $t\in[0,T]$, the smallest $n_t\ge d$ such that the top $n_t\times d$ submatrix of the matrix $(\langle \mathcal{S}^{-1}(\bm Z_t),(\nabla\varphi_i)^\top \bm\gamma(t,\mathcal{S}^{-1}(\bm Z_t),\cdot)\rangle)_{i\in\N}$ has the same (finite) rank as the full matrix. We can now carry out the constructions in \cite[Chapter 3, proof of Theorem 4.2]{karatzas-shreve} on the sets $\{t\in[0,T]:\,n_t=n\}$ for $n\ge d$ separately (picking the square root of the diffusion matrix as the one in \eqref{MV:reSPDE-CN} padded by $n-d$ zero columns), and combine them in the natural way to define $\bm B$. Thanks to \eqref{MV:reSPDE-CN}, we have for $(\mu_t:=\SS^{-1}(\bm Z_t))_{t\in[0,T]}$:
\begin{equation*}
\d\langle \mu_t, \varphi_i\rangle = \langle \mu_t, \gen^{\bm b,\bm a}_{t,\mu_t}\varphi_i\rangle \d t + \big\langle \mu_t,(\nabla\varphi_i)^\top \bm \gamma(t,\mu_t,\cdot)\big\rangle \d\bm B_t,\quad i\in\N.
\end{equation*}

\smallskip

For an arbitrary $\varphi\in C_c^\infty(\R^d)$, we use Lemma \ref{Lem_sep} to find a sequence $(\varphi_{n_k})_{k\in\N}$ such that $(\varphi_{n_k},\nabla\varphi_{n_k},\nabla^2\varphi_{n_k})$ converges uniformly to $(\varphi,\nabla\varphi,\nabla^2 \varphi)$. Then, $\langle \mu_t, \varphi_{n_k}\rangle\stackrel{k\to\infty}{\longrightarrow}\langle \mu_t, \varphi\rangle$ a.s.~for every $t\in[0,T]$. In addition,
\begin{equation*}
\begin{split}
&\,\bigg|\int_0^t \langle \mu_s, \gen^{\bm b,\bm a}_{s,\mu_s}\varphi_{n_k}\rangle \d s
-\int_0^t \langle \mu_s, \gen^{\bm b,\bm a}_{s,\mu_s}\varphi\rangle \d s\bigg| \\
&\qquad \le  \|\varphi_{n_k}-\varphi\|_{C^2_b(\R^d)}\,\int_0^t \big\langle \mu_s,\big|\bm b(s,\mu_s,\cdot)\big|+\big|\bm a(s,\mu_s,\cdot)\big|\big\rangle\,\mathrm{d}s
\stackrel{k\to\infty}{\longrightarrow}0\;\;\mathrm{a.s.}
\end{split}
\end{equation*}
for every $t\in[0,T]$, since \eqref{integ2} and Jensen's inequality imply that the latter integral has finite expectation and is therefore finite a.s. For the stochastic integral term, we apply the It\^o isometry:
\begin{align*} 
&\E\Big[\Big(\int_0^t\big\langle \mu_s,(\nabla\varphi^{n_k})^\top \bm \gamma(s,\mu_s,\cdot)\big\rangle \d\bm B_s-\int_0^t\big\langle \mu_s,(\nabla\varphi)^\top \bm \gamma(s,\mu_s,\cdot)\big\rangle \d\bm B_s \Big)^2\Big] \\
&\qquad = \E\bigg[\int_0^t\big|\big\langle \mu_s,(\nabla\varphi^{n_k}-\nabla\varphi)^\top
\bm \gamma(s,\mu_s,\cdot)\big\rangle\big|^2 \d s\bigg] \\
&\qquad \le  \|\varphi_{n_k}-\varphi\|^2_{C_b^1(\R^d)}\,\E\bigg[\int_0^t \big\langle \mu_s,\big|
\bm \gamma(s,\mu_s,\cdot)\big|\big\rangle^2 \d s\bigg]\stackrel{k\to\infty}{\longrightarrow}0,
\end{align*}
where the latter expectation is finite by \eqref{integ2} and Jensen's inequality. All in all,
\begin{equation*}
\langle\mu_t,\varphi\rangle=\langle\mu_0,\varphi\rangle+\int_0^t \langle \mu_s, \gen^{\bm b,\bm a}_{s,\mu_s}\varphi\rangle \d s+\int_0^t\big\langle \mu_s,(\nabla\varphi)^\top \bm \gamma(s,\mu_s,\cdot)\big\rangle \d\bm B_s\quad\mathrm{a.s.,}
\end{equation*}
for every $t\in[0,T]$. In view of the a.s.~continuity of both sides in $t$ (for the left-hand side, we exploit that uniform limits of continuous functions are continuous), we conclude that 
\begin{equation*}
\d\langle \mu_t, \varphi\rangle = \big\langle \mu_t, \gen_{t,\mu_t}^{\bm b,\bm a}\varphi\big\rangle \d t + \big\langle \mu_t,(\nabla\varphi)^\top \bm \gamma(t,\mu_t,\cdot)\big\rangle \d\bm B_t,\quad t \in [0,T], \quad\varphi \in C_c^\infty(\R^d),
\end{equation*}
as desired. We note lastly that we may take the probability space here to be the canonical one, $C([0,T];\P(\R^d)) \times C([0,T];\R^d)$, housing the process $\mu$ and the Brownian motion $\bm B$, equipped with its canonical filtration $\GG=(\G_t)_{t \in [0,T]}$. In particular, $\G_T$ is then countably generated, as the Borel $\sigma$-algebra on a Polish space. \qed


\section{Proof of the mimicking result, Corollary \ref{thm_mim}}\label{sec_application}

This section is devoted to the proof of Corollary \ref{thm_mim}. 

\medskip

\noindent\textbf{Proof of Corollary \ref{thm_mim}.}
The assumption $\F^{\bm X}_t \indep \F^{\bm W}_T \vee \G_T \,|\, \F^{\bm W}_t \vee \G_t$ for all $t$ implies $\F^{\bm X}_t \indep \G_T \,|\, \G_t$ for all $t$, by Lemma \ref{le:compatibility}. This in turn implies that $\mu_t=\L(\bm X_t \,|\,\G_T)=\L(\bm X_t \,|\,\G_t)$ a.s.
The existence of a jointly measurable version of $(t,m,\bm x) \mapsto \E[\bm b_t \,|\,\mu_t=m,\bm X_t=\bm x]$ follows by general arguments, and similarly for $\bm\sigma$: Let $\tau$ be an independent uniform random variable in $[0,T]$, and take a measurable version of $(t,m, \bm x) \mapsto \E[\bm b_\tau \,|\, \tau=t,\mu_\tau=m,\bm X_\tau= \bm x]$. See \cite[Section 5]{BS} for full details.

\medskip

We turn to the main line of proof. By It\^o's formula, for all $\varphi\!\in\! C_c^2(\R^d)$ and $t \in [0,T]$,
\begin{equation*}
\begin{aligned}
\varphi(\bm X_t)-\varphi(\bm X_0)=&\,\int_0^t \nabla \varphi(\bm X_s)\cdot\bm b_s + \frac{1}{2}\,\nabla^2 \varphi(\bm X_s):\big(\bm \sigma_s \bm \sigma_s^\top+\widehat{\bm \gamma}\widehat{\bm \gamma}^\top(s,\mu_s,\bm X_s)\big)\d s \\
&+\int_0^t\nabla \varphi(\bm X_s)^\top\bm \sigma_s\d\bm W_s+{\int_0^t\nabla \varphi(\bm X_s)^\top \widehat{\bm \gamma}(s,\mu_s,\bm X_s) \d\bm B_s}.
\end{aligned}
\end{equation*}
Next, we take the conditional expectation with respect to $\G_T$, using $\bm W\indep\G_T$, the assumption that $\F^{\bm{X}}_t \indep \F^{\bm{W}}_T \vee \G_T \,|\, \F^{\bm{W}}_t\vee\G_t$ for all $t \in [0,T]$, and the stochastic Fubini theorem (see Lemma \ref{le:fubini}) to get
\begin{align*}
\langle \mu_t,\varphi \rangle-\langle \mu_0,\varphi \rangle =& \,\int_0^t \E\Big[\nabla \varphi(\bm X_s)\cdot\bm b_s + \frac{1}{2}\,\nabla^2 \varphi(\bm X_s):\big(\bm \sigma_s \bm \sigma_s^\top+\widehat{\bm \gamma}\widehat{\bm \gamma}^\top(s,\mu_s,\bm X_s)\big)\Big|\G_s\Big]\d s
			\\
			&+\int_0^t\E\big[\nabla \varphi(\bm X_s)^\top\widehat{\bm \gamma}(s,\mu_s,\bm X_s)\big|\G_s\big]\d\bm B_s
			\\
			=&\,\int_0^t \E\Big[\nabla \varphi(\bm X_s)\cdot\bm b_s + \frac{1}{2}\,\nabla^2 \varphi(\bm X_s):\big(\bm \sigma_s \bm \sigma_s^\top+\widehat{\bm \gamma}\widehat{\bm \gamma}^\top(s,\mu_s,\bm X_s)\big)\Big|\G_s\Big]\d s
			\\
			&+\int_0^t\big\langle \mu_s,(\nabla\varphi)^\top\widehat{\bm \gamma}(s,\mu_s,\cdot)\big\rangle\d\bm B_s.
\end{align*}
Thus, for any $f\in C_c^2(\R^k)$ and $\bm\varphi=(\varphi_1,\ldots,\varphi_k)\in C_c^2(\R^d)^k$, It\^o's formula yields
\begin{align}
&\;f\big(\langle \mu_t,\bm\varphi \rangle\big)-f\big(\langle \mu_0,\bm\varphi \rangle\big) \nonumber  \\
&=\frac{1}{2}\int_0^t\sum_{i,j=1}^k\partial_{ij} f\big(\langle \mu_s,\bm \varphi \rangle\big)\,\big\langle \mu_s,(\nabla\varphi_i)^\top\widehat{\bm \gamma}(s,\mu_s,\cdot)\big\rangle\cdot\big\langle \mu_s,(\nabla\varphi_j)^\top\widehat{\bm \gamma}(s,\mu_s,\cdot)\big\rangle\d s \nonumber  \\
&\quad\,+\int_0^t\sum_{i=1}^k\partial_if\big(\langle \mu_s,\bm \varphi \rangle\big)\,\E\Big[\nabla \varphi_i(\bm X_s)\cdot\bm b_s + \frac{1}{2}\,\nabla^2 \varphi_i(\bm X_s):\big(\bm \sigma_s \bm \sigma_s^\top+\widehat{\bm \gamma}\widehat{\bm \gamma}^\top(s,\mu_s,\bm X_s)\big)\Big|\G_s\Big]\d s \nonumber \\
&\quad\,+\int_0^t\sum_{i=1}^k\partial_if\big(\langle \mu_s,\bm \varphi \rangle\big)\,\big\langle \mu_s,(\nabla\varphi_i)^\top\widehat{\bm \gamma}(s,\mu_s,\cdot)\big\rangle\d\bm B_s. \label{equ1}
\end{align}
We now take the expectation on both sides of \eqref{equ1} and recall that $\mu_s$ is $\G_s$-measurable, so that the expectation of the second term on the right-hand side of \eqref{equ1} evaluates to
\begin{equation*}
\begin{split}
&\,\int_0^t\sum_{i=1}^k\,\E\Big[\partial_if\big(\langle \mu_s,\bm \varphi \rangle\big)\,\E\Big[\nabla \varphi_i(\bm X_s)\cdot\bm b_s + \frac{1}{2}\,\nabla^2 \varphi_i(\bm X_s):\big(\bm \sigma_s \bm \sigma_s^\top+\widehat{\bm \gamma}\widehat{\bm \gamma}^\top(s,\mu_s,\bm X_s)\big)\Big|\G_s\Big]\Big]\d s \\
&=\int_0^t\sum_{i=1}^k\,\E\Big[\partial_if\big(\langle \mu_s,\bm \varphi \rangle\big)\Big(\nabla \varphi_i(\bm X_s)\cdot\bm b_s + \frac{1}{2}\,\nabla^2 \varphi_i(\bm X_s):\big(\bm \sigma_s \bm \sigma_s^\top+\widehat{\bm \gamma}\widehat{\bm \gamma}^\top(s,\mu_s,\bm X_s)\big)\!\Big)\Big]\d s.
\end{split}
\end{equation*}
In view of the tower rule with respect to the $\sigma$-algebra generated by $\mu_s$ and $\bm X_s$, 
\begin{equation*}
\begin{split}
&\;\E\Big[\partial_if\big(\langle \mu_s,\bm \varphi \rangle\big)\Big(\nabla \varphi_i(\bm X_s)\cdot\bm b_s + \frac{1}{2}\,\nabla^2 \varphi_i(\bm X_s):\big(\bm \sigma_s \bm \sigma_s^\top+\widehat{\bm \gamma}\widehat{\bm \gamma}^\top(s,\mu_s,\bm X_s)\big)\!\Big)\Big] \\
&=\E\Big[\partial_if\big(\langle \mu_s,\bm \varphi \rangle\big)\Big(\nabla \varphi_i(\bm X_s)\cdot\widehat{\bm b}(s,\mu_s,\bm X_s) + \frac{1}{2}\,\nabla^2 \varphi_i(\bm X_s):{\big(\widehat{\bm \sigma}\widehat{\bm\sigma}^\top+\widehat{\bm \gamma}\widehat{\bm\gamma}^\top\big)}(s,\mu_s,\bm X_s)\Big)\Big] \\
&=\E\Big[\partial_if\big(\langle \mu_s,\bm\varphi \rangle\big)\,\Big\langle \mu_s,\nabla \varphi_i\cdot\widehat{\bm b}(s,\mu_s,\cdot) + \frac{1}{2}\,\nabla^2 \varphi_i:{\big(\widehat{\bm \sigma}\widehat{\bm\sigma}^\top+\widehat{\bm \gamma}\widehat{\bm\gamma}^\top\big)}(s,\mu_s,\cdot)\Big\rangle\Big].
\end{split}
\end{equation*}
Now let $P_t := \L(\mu_t)$ for $t \in [0,T]$, set $\widehat{\bm a} = \widehat{\bm \sigma}\widehat{\bm\sigma}^\top+\widehat{\bm \gamma}\widehat{\bm\gamma}^\top$, and define $\gen^{\widehat{\bm b},\widehat{\bm a}}_{t,m}$ from $(\widehat{\bm b}, \widehat{\bm a})$ as in \eqref{def:gen-MV}.
Take expectations in \eqref{equ1} to get
\begin{equation*}
\begin{split}
&\,\int_{\P(\R^d)} f\big(\langle m,\bm\varphi\rangle \big) \,(P_t-P_0)(\!\d m)
= \int_0^t \int_{\P(\R^d)} \Bigg[\sum_{i=1}^k \partial_i f\big(\langle m,\bm\varphi\rangle \big) \langle m, \gen^{\widehat{\bm b},\widehat{\bm a}}_{s,m}\varphi_i \rangle \\
&\qquad\quad + \frac12\sum_{i,j=1}^k \partial_{ij} f\big(\langle m,\bm\varphi\rangle \big) \langle m, (\nabla\varphi_i)^\top \widehat{\bm \gamma}(s,m,\cdot) \rangle \cdot \langle m, (\nabla\varphi_j)^\top \widehat{\bm \gamma}(s,m,\cdot) \rangle \Bigg]  P_s(\!\d m)\d s.
\end{split}
\end{equation*}

We aim to apply Theorem \ref{Thm_PDEtoSPDE}, and for this purpose we claim that
\begin{equation*}
\int_0^T\int_{\P(\R^d)}\big\langle m,\big|\widehat{\bm b}(t,m,\cdot)\big|^p
		+\big|(\widehat{\bm \sigma}\widehat{\bm \sigma}^\top+\widehat{\bm \gamma}\widehat{\bm \gamma}^\top)(t,m,\cdot)\big|^p\big\rangle\, P_t(\!\d m)\d t<\infty.
\end{equation*}
Indeed, using $\mu_t=\L(\bm X_t\,|\,\G_T)$ and the definition of $\widehat{\bm b}$ along with Jensen's inequality, 
\begin{equation*}
\begin{split}
&\,\int_0^T\int_{\P(\R^d)}\big\langle m,\big|\widehat{\bm b}(t, m,\cdot)\big|^p\big\rangle\, P_t(\!\d m)\d t=\E\bigg[\int_0^T\big\langle\mu_t,\big|\widehat{\bm b}(t,\mu_t,\cdot)\big|^p\big\rangle\d t\bigg] \\
&=\E\bigg[\int_0^T\E\big[\big|\widehat{\bm b}(t,\mu_t,\bm X_t)\big|^p \,\big|\, \G_t \big] \d t\bigg] 
= \E\bigg[\int_0^T \big|\E[\bm b_t \,|\, \mu_t,\bm X_t ]\big|^p \d t\bigg] 
\le \E\bigg[\int_0^T|\bm b_t|^p\d t\bigg]<\infty,
\end{split}
\end{equation*}
and, similarly, 
\begin{align*}
&\int_0^T \int_{\P(\R^d)} \big\langle m,\big|(\widehat{\bm \sigma}\widehat{\bm \sigma}^\top+\widehat{\bm \gamma}\widehat{\bm \gamma}^\top)(t, m,\cdot)\big|^p\big\rangle\, P_t(\!\d m)\d t \\
&\qquad\leq\E\bigg[\int_0^T |\bm \sigma_t\bm \sigma_t^\top + \widehat{\bm \gamma}\widehat{\bm \gamma}^\top(t,\mu_t,\bm X_t)|^p \d t\bigg]  <\infty.
\end{align*}
We may now apply Theorem \ref{Thm_PDEtoSPDE} to construct a filtered probability space $(\widetilde{\Omega},\widehat{\GG},\widetilde{\PP})$, with $\widehat\G_T$ countably generated, supporting a $d$-dimensional $\widehat{\GG}$-Brownian motion $\widehat{\bm B}$ and a $\widehat{\GG}$-adapted $\P(\R^d)$-valued process $(\widehat{\mu}_t)_{t\in[0,T]}$ solving
\begin{equation}
\begin{split}
&\d\langle\widehat{\mu}_t,\varphi\rangle=\Big\langle\widehat{\mu}_t, \nabla\varphi\cdot\widehat{\bm b}(t,\widehat{\mu}_t,\cdot)+\frac{1}{2}\,\nabla^2\varphi:\big(\widehat{\bm\sigma}\widehat{\bm\sigma}^\top+\widehat{\bm\gamma}\widehat{\bm\gamma}^\top\big)(t,\widehat{\mu}_t,\cdot)\Big\rangle\d t \\
&\qquad\qquad\;\;+\big\langle \widehat{\mu}_t,(\nabla\varphi)^\top\widehat{\bm \gamma}(t,\widehat{\mu}_t,\cdot)\big\rangle\d\widehat{\bm B}_t
\end{split}
\end{equation}
with $\L(\widehat{\mu}_t)=P_t=\L(\mu_t)$, $t\in[0,T]$. 

\medskip

Next, we apply Theorem \ref{Thm_SPDEtoSDE_general} with the choice of coefficients as in Remark \ref{re:MVcase}:
\[
(t,\omega,x) \mapsto \big(\widehat{\bm b}(t,\widehat\mu_t(\omega),x),\widehat{\bm a}(t,\widehat\mu_t(\omega),x),\widehat{\bm \gamma}(t,\widehat\mu_t(\omega),x)\big).
\]
This yields an extension $(\widehat{\Omega}=\widetilde{\Omega}\times\Omega',\widehat{\mathbb F},\widehat{\PP})$ of the probability space $(\widetilde{\Omega},\widehat{\mathbb G},\widetilde{\PP})$ supporting $\widehat\FF$-Brownian motions $\widehat{\bm W}$ and $\widehat{\bm B}$, with $\widehat{\bm W}$ independent of $\widehat\G_T$, and a continuous $\widehat\FF$-adapted $d$-dimensional process $\widehat{\bm{X}}$ which satisfies
\begin{equation*}
\d\widehat{\bm{X}}_t = \widehat{\bm{b}}(t,\widehat{\mu}_t,\widehat{\bm{X}}_t)\d t + \widehat{\bm{\sigma}}(t,\widehat{\mu}_t,\widehat{\bm{X}}_t)\d\widehat{\bm{W}}_t + \widehat{\bm{\gamma}}(t,\widehat{\mu}_t,\widehat{\bm{X}}_t)\d\widehat{\bm{B}}_t,
\end{equation*} 
with $\widehat{\mu}_t = \L(\widehat{\bm{X}}_t\,|\,\widehat\G_T) = \L(\widehat{\bm{X}}_t\,|\,\widehat\G_t)$ a.s.~and such that $\F^{\widehat{\bm X}}_t \indep \F^{\widehat{\bm W}}_T \vee \widehat\G_T \,|\,\F^{\widehat{\bm W}}_t \vee \widehat\G_t$ for each $t \in [0,T]$. Finally, for all $g\in C_b(\R^d)$, $h\in C_b(\P(\R^d))$, and $t \in [0,T]$ we have:
\begin{equation}
\begin{split}
\E\big[g(\bm X_t)\,h(\mu_t)\big] =\E\big[\E\big[g(\bm X_t)\big|\G_t]\,h(\mu_t)\big] =\E\big[\langle \mu_t,g\rangle\,h(\mu_t)\big]
	=\E\big[\langle \widehat{\mu}_t,g\rangle\, h(\widehat{\mu}_t)\big]\\
	=\E\big[\E\big[g(\widehat{\bm X}_t)\big|\widehat\G_t]\,h(\widehat\mu_t)\big] = \E\big[g(\widehat{\bm X}_t)\,h(\widehat{\mu}_t)\big].
\end{split}
\end{equation}
Thus $(\widehat{\bm X}_t,\widehat \mu_t)\overset{d}{=}(\bm X_t,\mu_t)$ for all $t\in[0,T]$, which completes the proof. \qed


\section{Application to controlled McKean-Vlasov dynamics} \label{MFG_appl} 

The original mimicking theorem of \cite{Gyo,BS} can be used to prove both (1) the existence of optimal Markovian controls in classical stochastic optimal control problems and (2) the equivalence of open-loop and Markovian formulations of said control problems. Early references on this topic \cite{ElDuJe,HaLe} use an alternative but related approach often called \emph{Krylov's Markov selection} rather than the mimicking theorem (though see also \cite{KuSt}). 
To carry out this approach one typically first establishes the existence of an optimal control in a \emph{weak} or \emph{relaxed} sense, in which controls are allowed to include additional randomization; this relaxation facilitates compactness arguments. With a \emph{weak} optimal control in hand, the second step is to apply the mimicking theorem to project away the additional randomness and obtain a Markovian control which, under suitable convexity assumptions, achieves a lower cost. 

\medskip

Let us sketch a simple illustration of this second step. Suppose we are given a filtered probability space $(\Omega,\FF,\PP)$ supporting a controlled SDE of the form
\begin{align*}
\d \bm{X}_t = \bm{\alpha}_t \d t + \d {\bm{W}}_t,
\end{align*}
where $\bm{W}$ is an $\FF$-Brownian motion and $\bm{\alpha}$ an $\FF$-progressively measurable and square-integrable process. 
Then, defining the Markovian control $\widehat{\bm  \alpha}(t,x) = \E[\bm{\alpha}_t \,|\, {\bm{X}}_t=x]$, 
the mimicking theorem of \cite[Corollary 3.7]{BS} ensures that there exists a weak solution of
\[
\d\widehat{\bm X}_t = \widehat{\bm \alpha}(t,\widehat{\bm X}_t)\d t + \d \widehat{\bm W}_t
\]
such that $\widehat{\bm X}_t \stackrel{d}{=} {\bm{X}}_t$ for all $t$. If $f=f(x,a)$ and $g=g(x)$ are  suitably integrable cost functions, and $a \mapsto f(x,a)$ is convex for each $x$, then applying the identity $\widehat{\bm X}_t \stackrel{d}{=} {\bm{X}}_t$ followed by Jensen's inequality yields 
\begin{align*}
\E\left[ \int_0^Tf(\widehat{\bm X}_t,\widehat{\bm \alpha}(t,\widehat{\bm X}_t))\d t + g(\widehat{\bm X}_T) \right] \le \E\left[ \int_0^Tf({\bm{X}}_t,\bm{\alpha}_t)\d t + g({\bm{X}}_T) \right]. 
\end{align*}
In other words, starting from any \emph{open-loop} control $\bm{\alpha}$ as above (i.e., progressively measurable with respect to some filtration with respect to which $\bm W$ is a Brownian motion), we can construct a Markovian control achieving a lower cost. In particular, the optimal value over open-loop controls equals that over Markovian controls, and if the open-loop problem admits an optimizer then so does the Markovian problem.

\medskip

This procedure has been applied in the setting of mean field control in \cite{lacker2017limit}. The cost functions therein depend nonlinearly on the law of the state process, i.e., $f({\bm{X}}_t,\bm{\alpha}_t)$ and $g({\bm{X}}_T)$ are replaced with $f(\L({\bm{X}}_t),\bm{\alpha}_t)$ and $g(\L({\bm{X}}_t))$. The argument given above applies essentially without change, because the construction of a Markovian control does not alter the time-$t$ marginal laws. Similarly, this method works well for proving the existence of Markovian equilibria for mean field games \cite{lacker2015mean,BeCaDi}.
It should be stressed that the results cited in this paragraph are limited to settings without common noise.

\medskip

For mean field control problems with common noise, the situation is substantially more complex because the measure flow involved is stochastic. The following explains how our mimicking result, Corollary \ref{thm_mim}, can be used in this setting, described next. We are given a Polish space $A$ (the control space), an initial distribution $\lambda_0 \in \P(\R^d)$, and measurable functions
\begin{align*}
(\bm{b},\bm{\sigma},f) &: [0,T] \times \R^d \times \P(\R^d) \times A \to \R^d \times \R^{d \times d} \times \R, \\
\bm\gamma &: [0,T] \times \R^d \times \P(\R^d) \to \R^{d \times d}, \\
g &: \R^d \times \P(\R^d) \to \R.
\end{align*}
The common noise coefficient $\bm\gamma$ is uncontrolled, as is the case almost universally in the literature for technical reasons.
We assume a form of \emph{Roxin's condition}, namely that the subset of $\R^d \times \R^{d \times d} \times \R$ given by
\begin{align}
\begin{split}
&\left\{\big(\bm{b}(t,\bm x,m,a),\bm{\sigma}\bm{\sigma}^\top(t,\bm x,m,a),z\big) : z \in \R, \, a \in A, \, z \ge f(t,\bm x,m,a)\right\} \\
&\quad \ \text{is closed and convex, for each } (t,\bm x,m) \in [0,T] \times \R^d \times \P(\R^d).
\end{split} \label{roxin}
\end{align}
For example, this holds if $\bm{\sigma}$ is uncontrolled, $\bm{b}$ is linear in $a$, and $f$ is convex in $a$ (with $A$ a convex subset of a vector space). Alternatively, this includes \emph{relaxed control} setups in which $A = \P(\widetilde A)$ for some other Polish space $\widetilde A$ and $(\bm{b},\bm{\sigma}\bm{\sigma}^\top,f)$ is linear in the sense that $(\bm{b},\bm{\sigma}\bm{\sigma}^\top,f)(t,\bm x,m,a) = \int_{\widetilde A}(\widetilde{\bm b},\widetilde{\bm \sigma}\widetilde{\bm \sigma}^\top,\widetilde f)(t,\bm x,m,\widetilde a)\,a(d\widetilde a)$.

\medskip

The mean field control problem with common noise is, roughly speaking, to choose a control $\bm{\alpha}$ to minimize the objective
\begin{align*}
\E\left[ \int_0^T f(t,{\bm{X}}_t,\mu_t,\bm{\alpha}_t)\d t + g({\bm{X}}_T,\mu_T)\right],
\end{align*}
where the state process $\bm X$ is given by
\begin{align*}
\d {\bm{X}}_t &= b(t,{\bm{X}}_t,\mu_t,\bm{\alpha}_t)\d t + \bm{\sigma}(t,{\bm{X}}_t,\mu_t,\bm{\alpha}_t)\d {\bm{W}}_t + \gamma(t,{\bm{X}}_t,\mu_t)\d \bm{B}_t, \\
\mu_t &= \L({\bm{X}}_t\,|\,\G_t).
\end{align*}
Here $\bm{W}$ is independent of the filtration $\GG$, and $\bm{B}$ is a $\GG$-Brownian motion.

\medskip

The few recent papers on mean field control with common noise, such as \cite{PhWe,DjPoTa} have proposed various notions of admissible controls, usually with the goal of deriving a dynamic programming principle and an associated Hamilton-Jacobi-Bellman equation. See also \cite[Chapter I.6]{CD2} for the case without common noise. 
The following definition of a weak control is essentially \cite[Definition 2.1]{DjPoTa}, which is closely related to the weak solution concept introduced for mean field games in the earlier papers \cite{CDL,La}.

\begin{definition} \label{def:weak control}
A \emph{weak control} is a tuple $\mathcal{R}=(\Omega,\FF,\GG,\PP,\bm B,\bm W,\mu,\bm X,\bm{\alpha})$ such that:
\begin{enumerate}[(1)]
\item $(\Omega,\FF,\PP)$ is a filtered probability space and $\GG$ a subfiltration of $\FF$.
\item $\bm W$ and $\bm B$ are independent $\FF$-Brownian motions of dimension $d$.
\item $\bm B$ is $\GG$-adapted.
\item $\mu$ is a continuous $\GG$-adapted $\P(\R^d)$-valued process satisfying $\mu_t=\L({\bm{X}}_t\,|\,\G_T)=\L({\bm{X}}_t\,|\,\G_t)$ a.s., for each $t \in [0,T]$.
\item $\bm{\alpha}$ is an $\FF$-progressively measurable $A$-valued process.
\item $\bm X$ is a continuous $\FF$-adapted $\R^d$-valued process satisfying 
\begin{align*}
\d \bm{X}_t = \bm b(t,{\bm{X}}_t,\mu_t,\bm{\alpha}_t)\d t + \bm{\sigma}(t,{\bm{X}}_t,\mu_t,\bm{\alpha}_t)\d{\bm{W}}_t + \bm\gamma(t,{\bm{X}}_t,\mu_t)\d \bm{B}_t,\quad {\bm{X}}_0 \sim \lambda_0.
\end{align*}
\item For some $p > 1$, we have:
\begin{align*}
& \E\left[\int_0^T |\bm b(t,{\bm{X}}_t,\mu_t,\bm{\alpha}_t)|^p + |\bm{\sigma}\bm{\sigma}^\top(t,{\bm{X}}_t,\mu_t,\bm{\alpha}_t)|^p + |\bm\gamma(t,{\bm{X}}_t,\mu_t)|^{2p} \d t\right] < \infty, \\
& \E\left[\int_0^T |f(t,{\bm{X}}_t,\mu_t,\bm{\alpha}_t)|+ |g({\bm{X}}_T,\mu_T)|\d t\right] < \infty.
\end{align*}
\item ${\bm{X}}_0$, $\bm W$, and $\G_T$ are independent.
\item \label{MFcontrol-compatibility} $({\bm{X}}_s,\bm{\alpha}_s)_{s\in[0,t]} \indep \F^{\bm{W}}_T \vee \G_T\,|\,\F^{\bm{W}}_t \vee \G_t$, for each $t \in [0,T]$.
\end{enumerate}
The \emph{cost} of the weak control $\mathcal{R}$ is the quantity
\begin{align*}
J(\mathcal{R}) := \E\left[ \int_0^T f(t,{\bm{X}}_t,\mu_t,\bm{\alpha}_t)\d t + g({\bm{X}}_T,\mu_T)\right],
\end{align*}
which is well-defined thanks to property (7). A \emph{Markovian control} is a weak control for which there exists a measurable function $\widehat{\bm \alpha} : [0,T] \times \R^d \times \P(\R^d) \to A$ such that $\bm{\alpha}_t=\widehat{\bm \alpha}(t,{\bm{X}}_t,\mu_t)$ a.s.\ for a.e.\ $t \in [0,T]$.  
\end{definition}

\begin{remark} \label{re:MFCcompatibility}
In \cite[Definition 2.3]{DjPoTa2}, a weak control is defined similarly. The only significant difference is that instead of (9) they require that $\F_t \vee \F^{\bm W}_T \indep \G_T\,|\,\G_t$ for each $t \in [0,T]$. 
But, by Lemma \ref{le:compatibility}, this is equivalent to $\F_t \indep \F^{\bm W}_T \vee \G_T \,|\, \F^{\bm W}_t \vee \G_t$ for all $t \in [0,T]$, which implies (9). If $\FF=\FF^{\bm X} \vee \FF^{\bm \alpha} \vee \FF^{\bm W} \vee \GG$ is the minimal filtration, these conditions are equivalent.
\end{remark}

Our main result in this section shows that a weak control can be turned into a superior Markovian control without changing the marginal flow $(\L(\bm{X}_t,\mu_t))_{t \in [0,T]}$. 


\begin{theorem} \label{thm:MFcontrol}
Assume the convexity condition \eqref{roxin}. Then, for every weak control $\mathcal{R}=(\Omega,\FF,\GG,\PP,\bm B,\bm W,\mu,\bm X,\bm{\alpha})$ there is a Markovian control $\widehat{\mathcal{R}}=(\widehat \Omega,\widehat \FF,\widehat \GG,\widehat \PP,\widehat{\bm B},\widehat{\bm W},\widehat \mu,\widehat{\bm X},\widehat{\bm \alpha})$ satisfying $\widehat \PP \circ (\widehat{\bm X}_t,\widehat \mu_t)^{-1} = \PP \circ ({\bm{X}}_t,\mu_t)^{-1}$ for all $t \in [0,T]$, and also $J(\widehat{\mathcal{R}}) \le J(\mathcal{R})$.
In particular, the optimal cost over Markovian controls equals the optimal cost over weak controls, and the existence of an optimal weak control implies the existence of an optimal Markovian control.
\end{theorem}

\noindent\textbf{Proof.} Let $\mathcal{R}=(\Omega,\FF,\GG,\PP,\bm B,\bm W,\mu,\bm X,\bm{\alpha})$ be a weak control. Using the convexity assumption \eqref{roxin}, we may find a measurable\footnote{This requires a measurable selection argument, cf.\ \cite[Theorem A.9]{HaLe} or \cite[Lemma 3.1]{DuSt}.} function $\widehat{\bm \alpha} : [0,T] \times \R^d \times \P(\R^d) \to A$ such that
\begin{align}
\bm{b}(t,{\bm{X}}_t,\mu_t,\widehat{\bm \alpha}(t,{\bm{X}}_t,\mu_t)) &= \E[\bm{b}(t,{\bm{X}}_t,\mu_t,\bm{\alpha}_t) \, | \, {\bm{X}}_t,\mu_t], \nonumber \\
\bm{\sigma}\bm{\sigma}^\top(t,{\bm{X}}_t,\mu_t,\widehat{\bm \alpha}(t,{\bm{X}}_t,\mu_t)) &= \E[\bm{\sigma}\bm{\sigma}^\top(t,{\bm{X}}_t,\mu_t,\bm{\alpha}_t) \, | \, {\bm{X}}_t,\mu_t], \nonumber \\
f(t,{\bm{X}}_t,\mu_t,\widehat{\bm \alpha}(t,{\bm{X}}_t,\mu_t)) &\le \E[f(t,{\bm{X}}_t,\mu_t,\bm{\alpha}_t) \,|\, {\bm{X}}_t,\mu_t] \label{pf:mfcontrol-1}
\end{align}
a.s., for a.e.\ $t \in [0,T]$.  By Corollary \ref{thm_mim}, there exists a filtered probability space $(\widehat\Omega,\widehat{\mathbb F},\widehat{\PP})$ with a subfiltration $\widehat{\mathbb G}\subset\widehat{\mathbb F}$, two $\widehat\FF$-Brownian motions $\widehat{\bm W}$ and $\widehat{\bm B}$ with $\widehat{\bm W}$ independent of $\widehat{\G}_T$ and $\widehat{\bm B}$ adapted to $\widehat\GG$, and a weak solution $\widehat{\bm X}$ of
\begin{align*}
\d\widehat{\bm X}_t &= \bm b\big(t,\widehat{\bm X}_t,\widehat{\mu}_t,\widehat{\bm \alpha}(t,\widehat{\bm X}_t,\widehat\mu_t)\big)\d t+\bm{\sigma}\big(t,\widehat{\bm X}_t,\widehat{\mu}_t,\widehat{\bm \alpha}(t,\widehat{\bm X}_t,\widehat\mu_t)\big)\d\widehat{\bm W}_t+ \bm\gamma(t,\widehat{\bm X}_t,\widehat{\mu}_t)\d\widehat{\bm B}_t, \\
\widehat{\mu}_t &= \L\big(\widehat{\bm X}_t \,|\, \widehat{\G}_T\big) = \L\big(\widehat{\bm X}_t \,|\, \widehat{\G}_t\big), \quad t \in [0,T]
\end{align*}
such that $(\widehat{\bm X}_t,\widehat \mu_t)\overset{d}{=}(\bm {\bm{X}}_t,\mu_t)$ for all $t\in[0,T]$. Moreover, $\F^{\widehat{\bm X}}_t \indep \F^{\widehat{\bm W}}_T \vee \widehat\G_T\,|\,\F^{\widehat{\bm W}}_t \vee \widehat\G_t$ for each $t \in [0,T]$. Let $\widehat{\bm \alpha}_t := \widehat{\bm \alpha}(t,\widehat{\bm X}_t,\widehat\mu_t)$ and $\widehat{\mathcal{R}} :=(\widehat \Omega,\widehat \FF,\widehat \GG,\widehat \PP,\widehat{\bm B},\widehat{\bm W},\widehat \mu,\widehat{\bm X},\widehat{\bm \alpha})$.

\medskip

It is immediate that $\widehat{\mathcal{R}}$ satisfies conditions (1--6) of Definition \ref{def:weak control}. Condition (7) follows from Jensen's inequality and \eqref{pf:mfcontrol-1}. To check (8), note first that since $\bm{X}_0 \indep \G_T$ we have $\mu_0=\L(\bm{X}_0\,|\,\G_T)=\L(\bm{X}_0)$. Since $\widehat\mu_0 \stackrel{d}{=}\mu_0$, this implies $\L(\widehat{\bm X}_0\,|\, \widehat\G_T)=\widehat\mu_0=\L(\bm{X}_0)=\L(\widehat{\bm X}_0)$, which shows that $\widehat{\bm X}_0 \indep \widehat\G_T$. 
Next, use Lemma \ref{le:compatibility} (the implication (1) $\Rightarrow$ (2) with $\HH=\FF^{\widehat{\bm X}}$, $\FF^{\widehat{\bm W}}$ and $\widehat\GG$) to deduce that $\widehat{\bm{W}}$ is independent of $\sigma(\widehat{\bm X}_0) \vee \widehat\G_T$. These two facts combined yield (8), the mutual independence of $\widehat{\bm X}_0$, $\widehat{\bm W}$, and $\widehat\G_T$.
To check condition (9), simply use the conditional independence of $\F^{\widehat{\bm X}}_t \indep \F^{\widehat{\bm W}}_T \vee \widehat\G_T\,|\,\F^{\widehat{\bm W}}_t \vee \widehat\G_t$ for each $t \in [0,T]$ along with the expression $\widehat{\bm \alpha}_t = \widehat{\bm \alpha}(t,\widehat{\bm X}_t,\widehat\mu_t)$.

\medskip

We have thus shown that $\widehat{\mathcal{R}}$ is a Markovian control.
Finally, we use \eqref{pf:mfcontrol-1} along with Fubini's theorem and the fact that $(\widehat{\bm X}_t,\widehat \mu_t)\overset{d}{=}(\bm {\bm{X}}_t,\mu_t)$ for all $t\in[0,T]$ to conclude that
\begin{align*}
\qquad\qquad\qquad J(\widehat{\mathcal{R}}) &= \E\left[\int_0^T f(t,\widehat{\bm X}_t,\widehat\mu_t,\widehat{\bm \alpha}(t,\widehat{\bm X}_t,\widehat\mu_t))\d t + g(\widehat{\bm X}_T,\widehat\mu_T)\right] \\
	&= \E\left[\int_0^T f(t,{\bm X}_t,\mu_t,\widehat{\bm \alpha}(t,{\bm X}_t,\mu_t))\d t + g({\bm X}_T,\mu_T)\right] \\
	&\le \E\left[\int_0^T f(t,{\bm X}_t,\mu_t,\bm{\alpha}_t)\d t + g({\bm X}_T,\mu_T)\right] = J(\mathcal{R}). \qquad\qquad\qquad\qed
\end{align*}

\smallskip

Theorem \ref{thm:MFcontrol} is useful in clarifying the relationships between various natural formulations of the mean field control problem, a delicate matter in the case of common noise \cite{DjPoTa2}. In addition, the very recent results of \cite{DjPoTa2} characterize limits of $n$-player optimal control problems in terms of \emph{weak controls}; Theorem \ref{thm:MFcontrol} then implies that the limits can be characterized instead in terms of \emph{Markovian controls}, at least if one is only interested in the convergence of the time-$t$ marginals of the measure-valued process.

\medskip

Similar arguments may be possible for mean field games with common noise, see \cite{CD2} for a recent comprehensive overview. General results on existence, uniqueness, and limit theory require working with a notion of \emph{weak equilibrium} which is well-suited for weak convergence arguments \cite{CDL,La,CD2}. A reduction from such a weak equilibrium to some form of \emph{Markovian equilibrium} is a much more delicate task than for the control problems described in this section, and we leave it for future work.


\appendix

\section{Proof of compatibility lemma} \label{se:misc}

Here we give the proof of Lemma \ref{le:compatibility}. We first provide a few alternative characterizations of the compatibility condition, also known as \emph{immersion}, or the \emph{H-hypothesis}, in the literature. Given two filtrations $\FF^1=(\F^1_t)_{t \in [0,T]}$ and $\FF^2=(\F^2_t)_{t \in [0,T]}$ defined on the same space, with $\F^1_t \subset \F^2_t$ for all $t$, the following are well known to be equivalent \cite[Theorem 3]{BrYo}:
\begin{itemize}
\item $\F^2_t \indep \F^1_T \,|\, \F^1_t$, for each $t \in [0,T]$.
\item $\FF^1$ is immersed in $\FF^2$, in the sense that every bounded $\FF^1$-martingale is also an $\FF^2$-martingale.
\item $\E[Z \, | \, \F^2_t] = \E[Z \, | \, \F^1_t]$ a.s., for each $t \in [0,T]$ and each bounded $\F^1_T$-measurable random variable $Z$.
\end{itemize}
Note that if $\FF^2=\FF^1 \vee \FF^3$ for some auxiliary filtration $\FF^3=(\F^3_t)_{t \in [0,T]}$, 
then the first bullet point can be rewritten as $\F^3_t \indep \F^1_T\,|\,\F^1_t$ for all $t \in [0,T]$.
We will use these characterizations in the proof of Lemma \ref{le:compatibility}.

\medskip 

\noindent\textbf{Proof of Lemma \ref{le:compatibility}.}
We first show that (1) $\Longrightarrow$ (2a): We use that (1) and (2a) are respectively equivalent to
\begin{enumerate}[(1')]
\item[(1$'$)] $\E[Z\,|\,\H_t \vee \F^{\bm{W}}_t \vee \G_t]=\E[Z\,|\,\F^{\bm{W}}_t \vee \G_t]$ a.s., for each $t \in [0,T]$ and each bounded $\F^{\bm{W}}_T \vee \G_T$-measurable random variable $Z$.
\item[(2a$'$)] $\E[Z\,|\,\H_t \vee \F^{\bm{W}}_t \vee \G_t]=\E[Z\,|\,\H_0 \vee \F^{\bm{W}}_t \vee \G_t]$ a.s., for each $t \in [0,T]$ and each bounded $\H_0 \vee \F^{\bm{W}}_T \vee \G_T$-measurable random variable $Z$.
\end{enumerate}
Fix $t \in [0,T]$. Consider $Z$ of the form $Z=Z^H_0Z_T^WZ_T^G$, where $Z^H_0$, $Z_T^W$, and $Z_T^G$ are bounded random variables, measurable with respect to $\H_0$, $\F^{\bm{W}}_T$, and $\G_T$, respectively. Then factoring out $Z^H_0$ and applying (1$'$) yields
\begin{align*}
\E[Z\,|\,\H_t \vee \F^{\bm{W}}_t \vee \G_t] = Z^H_0\E[Z_T^WZ_T^G\,|\,\H_t \vee \F^{\bm{W}}_t \vee \G_t] 
	= Z^H_0\E[Z_T^WZ_T^G\,|\, \F^{\bm{W}}_t \vee \G_t].
\end{align*}
This shows that $\E[Z\,|\,\H_t \vee \F^{\bm{W}}_t \vee \G_t]$ is $\H_0 \vee \F^{\bm{W}}_t \vee \G_t$-measurable, which in turn implies $\E[Z\,|\,\H_t \vee \F^{\bm{W}}_t \vee \G_t]=\E[Z\,|\,\H_0 \vee \F^{\bm{W}}_t \vee \G_t]$ a.s. As this holds for all $Z$ of the form $Z=Z^H_0Z_T^WZ_T^G$ described above, by a monotone class argument the same holds for any bounded $\H_0 \vee \F^{\bm{W}}_T \vee \G_T$-measurable random variable $Z$. This proves (2a$'$), which in turn implies (2a).

\medskip

In showing both that (1) $\Longrightarrow$ (2b) and that (2) $\Longrightarrow$ (1) we rely on the fact (mentioned above) that (1), (2a), and (2b) are respectively equivalent to
\begin{enumerate}[(1'')]
\item[(1$''$)] Every bounded $\FF^{\bm{W}} \vee \GG$-martingale is also a $\HH \vee \FF^{\bm{W}} \vee \GG$-martingale.
\item[(2a$''$)] Every bounded $(\H_0 \vee \F^{\bm{W}}_t \vee \G_t)_{t \in [0,T]}$-martingale is also a $\HH \vee \FF^{\bm{W}} \vee \GG$-martingale.
\item[(2b$''$)] Every bounded $\GG$-martingale is also a $\HH \vee \GG$-martingale.
\end{enumerate}
Now, to prove that (1) $\Longrightarrow$ (2b), we instead show that (1$''$) $\Longrightarrow$ (2b$''$). Since $\F^{\bm{W}}_T \indep \G_T$, it is straightforward to check that every bounded $\GG$-martingale is also an $\FF^{\bm{W}} \vee \GG$-martingale and thus, by (1$''$), also an $\HH \vee \FF^{\bm{W}} \vee \GG$-martingale. But any $\HH \vee \FF^{\bm{W}} \vee \GG$-martingale which is adapted to the smaller filtration $\HH \vee \GG$ must be also an $\HH \vee \GG$-martingale. This shows that (1$''$) $\Longrightarrow$ (2b$''$).

\medskip

To prove that (1) $\Longrightarrow$ (2c), note first that $\F^{\bm{W}}_0$ is the trivial $\sigma$-algebra. Hence, (1) implies that $\H_0 \indep \F^{\bm{W}}_T \vee \G_T \,|\, \G_0$. For random variables $Z^H_0$, $Z^W_T$, and $Z^G_T$, measurable with respect to $\H_0$, $\F^{\bm{W}}_T$, and $\G_T$, respectively, we thus have
\begin{align*}
\E\big[Z^H_0Z^W_TZ^G_T\big] = \E\big[\E[Z^H_0\,|\,\G_0]Z^W_TZ^G_T\big] = \E[Z^W_T]\,\E\big[\E[Z^H_0\,|\,\G_0]Z^G_T\big] 
	= \E[Z^W_T]\,\E[ Z^H_0 Z^G_T].
\end{align*}
Indeed, the second step follows from the assumed independence of $\bm{W}$ and $\G_T$, and the other two steps from the conditional independence $\H_0 \indep \F^{\bm{W}}_T \vee \G_T \,|\, \G_0$. This proves (2c).

\medskip


Next, we show that (2) $\Longrightarrow$ (1). With (1$''$) and (2a$''$) in mind,
it suffices to check that every bounded $\FF^{\bm{W}} \vee \GG$-martingale is also an $(\H_0 \vee \F^{\bm{W}}_t \vee \G_t)_{t \in [0,T]}$-martingale, since then (2a$''$) will show that it is therefore also an $\HH \vee \FF^{\bm{W}} \vee \GG$-martingale. In light of the equivalences summarized before the proof, it suffices to instead show that $\H_0 \indep \F^{\bm{W}}_T \vee \G_T \,|\, \F^{\bm{W}}_t \vee \G_t$ for each $t \in [0,T]$. To prove this, note that for any bounded $\H_0$-measurable random variable $Z$ we may use the independence of $\bm{W}$ and $\H_0 \vee \G_T$ from (2c) to deduce
\begin{align*}
\E[Z \, | \, \F^{\bm{W}}_T \vee \G_T] &= \E[Z \, | \,  \G_T] = \E[Z \, | \, \G_0],
\end{align*}
where the last step follows from (2b). The same argument applies with $t$ in place of $T$:
\begin{align*}
\E[Z \, | \, \F^{\bm{W}}_t \vee \G_t] &= \E[Z \, | \,  \G_t] = \E[Z \, | \, \G_0].
\end{align*}
This shows that $\E[Z \, | \, \F^{\bm{W}}_T \vee \G_T] = \E[Z \, | \, \F^{\bm{W}}_t \vee \G_t]$. We deduce that $\H_0 \indep \F^{\bm{W}}_T \vee \G_T \,|\, \F^{\bm{W}}_t \vee \G_t$, which completes the proof of (1).

\medskip

It is clear that (3) $\Longrightarrow$ (1). To prove that (1) $\Longrightarrow$ (3), note that $\F^{\bm{W}}_T = \F^{\bm{W}}_t \vee \sigma(\bm{W}_s-\bm{W}_t : s \in [t,T])$. With this in mind, the independence of $\bm{W}$ and $\G_T$ then easily implies that $\F^{\bm{W}}_T \indep \G_T$ given $\F^{\bm{W}}_t \vee \G_t$, which is enough to deduce (3) from (1).

\medskip

We turn to (3) $\Longrightarrow$ (4a). In this and the remaining steps of the proof we write
\[
\F^{\bm W}_{> t} := \sigma(\bm{W}_s-\bm{W}_t : s \in [t,T]), \qquad t\in [0,T].
\]
It suffices to verify that $\G_T \vee \F^{\bm{W}}_t \vee \H_t \indep \F^{\bm W}_{> t}$ for each $t \in [0,T]$. Fix $t \in [0,T]$, and let $Z$ be a bounded $\F^{\bm W}_{> t}$-measurable random variable. Independence of $\bm{W}$ and $\G_T$ yields $\E[Z \,|\,\F^{\bm{W}}_t \vee \G_t] = \E[Z ]$. Using (3) we have $Z \indep \H_t \vee \G_T \,|\, \F^{\bm{W}}_t \vee \G_t$, and thus
\begin{align*}
\E[Z \, | \, \F^{\bm{W}}_t \vee \H_t \vee \G_T] = \E[ Z \, | \, \F^{\bm{W}}_t \vee \G_t] = \E[Z].
\end{align*}

\smallskip

Next, we prove (1) $\Longrightarrow$ (4b). As discussed above, (1) implies (1$''$), that every bounded $\FF^{\bm{W}} \vee \GG$-martingale is also an $\HH \vee \FF^{\bm{W}} \vee \GG$-martingale. Because $\bm{W} \indep \G_T$, it is straightforward to check that every bounded $\GG$-martingale is also an $\FF^{\bm{W}} \vee \GG$-martingale. Thus, every bounded $\GG$-martingale is an $\HH \vee \FF^{\bm{W}} \vee \GG$-martingale, which is equivalent to (4b).

\medskip

We next show that (4) $\Longrightarrow$ (1). Fix $t \in [0,T]$ and bounded random variables $Z_t^H$, $Z^W_t$, $Z_{t+}^W$, $Z^G_T$, $Z^G_t$, which are measurable with respect to $\H_t$, $\F^{\bm{W}}_t$, $\F^{\bm W}_{> t}$, $\G_T$, and $\G_t$, respectively. The independence of $\bm{W}$ and $\G_T$ easily implies
\begin{align*}
\E[Z^W_{t+}Z^G_T\,|\,\F^{\bm{W}}_t \vee \G_t] = \E[Z^W_{t+}]\,\E[Z^G_T\,|\,\G_t].
\end{align*}
From (4a) we know that $Z^W_{t+} \indep \G_T \vee \H_t \vee \F^{\bm{W}}_t$, and we compute
\begin{align*}
\E\big[Z^H_tZ^W_{t+}Z^G_TZ^W_tZ^G_t \big] &= \E[Z^W_{t+}]\,\E\big[Z^H_tZ^G_TZ^W_tZ^G_t \big] \\
	&= \E[Z^W_{t+}]\,\E\big[ Z^H_tZ^W_t \E[Z^G_T\,|\,\G_t]Z^G_t \big] \\
	&= \E\big[ Z^H_tZ^W_t \E[Z^W_{t+}Z^G_T\,|\,\F^{\bm{W}}_t \vee \G_t]Z^G_t \big] \\
	&= \E\big[ \E[Z^H_t \, | \, \F^{\bm{W}}_t \vee \G_t]\,\E[Z^W_{t+}Z^G_T\,|\,\F^{\bm{W}}_t \vee \G_t]Z^W_tZ^G_t \big].
\end{align*}
Indeed, the second step follows from (4b). This shows that $\H_t$ and $\G_T \vee \F^{\bm W}_{> t}$ are conditionally independent given $\F^{\bm{W}}_t \vee \G_t$. Since $\F^{\bm{W}}_T = \F^{\bm{W}}_t \vee \F^{\bm W}_{> t}$, we deduce (1).

\medskip

We proceed to (5)  $\Longrightarrow$ (3). Fix a $t \in [0,T]$. Since $\H_t \vee \F^{\bm W}_T \indep \G_T\,|\,\G_t$ by (5),
\begin{align*}
\E[Z\,|\,\H_t \vee \F^{\bm W}_T \vee \G_t] &= \E[Z\,|\,\G_t]
\end{align*}
for all bounded $\G_T$-measurable random variables $Z$. Using independence of $\bm{W}$ and $\G_T$, for such $Z$ we easily deduce $\E[Z\,|\,\G_t] = \E[Z\,|\,\F^{\bm W}_t \vee \G_t]$, and thus
\begin{align*}
\E[Z\,|\,\H_t \vee \F^{\bm W}_T \vee \G_t] &= \E[Z\,|\,\F^{\bm W}_t \vee \G_t].
\end{align*}
Hence, $\G_T \indep \H_t \vee \F^{\bm W}_T\,|\,\F^{\bm W}_t \vee \G_t$. Writing $\F^{\bm W}_T=\F^{\bm W}_{>t} \vee \F^{\bm W}_t$, we see that $\G_T \indep \H_t \vee \F^{\bm W}_{> t}\,|\,\F^{\bm W}_t \vee \G_t$. But $\F^{\bm W}_{> t} \indep \F_t$ and $\HH \vee \FF^{\bm W} \vee \GG \subset \FF$ imply that also $\F^{\bm W}_{> t} \indep \H_t\,|\,\F^{\bm W}_t \vee \G_t$. Consequently, the three $\sigma$-algebras $\F^{\bm W}_{> t}$, $\H_t$, and $\G_T$ are conditionally independent given $\F^{\bm W}_t \vee \G_t$. This proves (3).

\medskip

Finally, we show (4) $\Longrightarrow$ (5). Fix a $t \in [0,T]$.  
Consider bounded random variables $Z^{H,W}_t$, $Z^W_{t+}$, $Z^G_T$, and $Z^G_t$, measurable with respect to $\H_t \vee \F^{\bm W}_t$, $\F^{\bm W}_{>t}$, $\G_T$, and $\G_t$, respectively. 
Then
\begin{align*}
\E[Z^{H,W}_tZ^W_{t+}Z^G_TZ^G_t] &= \E[Z^W_{t+}]\,\E[Z^{H,W}_tZ^G_TZ^G_t] \\
	&= \E[Z^W_{t+}]\,\E[Z^{H,W}_t\E[Z^G_T\,|\,\G_t]Z^G_t] \\
	&= \E[Z^{H,W}_tZ^W_{t+}\E[Z^G_T\,|\,\G_t]Z^G_t].
\end{align*}
Indeed, the first and third equality use $\F^{\bm W}_{> t} \indep \G_T \vee \F^{\bm W}_t \vee \H_t$, which is (4a), and the second equality uses $\H_t \vee \F^{\bm W}_t \indep \G_T \,|\, \G_t$, which is (4b). We conclude that $\H_t \vee \F^{\bm W}_T \indep \G_T\,|\,\G_t$. \qed
%
%

\section{Derivation of SPDE from SDE} \label{sec_app_SDEtoSPDE}

In this section we give the proof of Proposition \ref{pr:SDEtoSPDE}. We first state a stochastic version of Fubini's theorem, which can be found in various references, such as \cite[Lemma 1.11]{HSS} or \cite[Lemma A.1]{LaSv}, in different degrees of generality. We include a short proof in our setting for the sake of completeness. 

\begin{lemma} \label{le:fubini}
Use the notation and assumptions of Lemma \ref{le:compatibility}, and assume one (and thus all) of the equivalent conditions (1)--(4) therein holds. Assume also that $(\Omega,\FF,\PP)$ supports an $\FF$-Brownian motion $\bm{B}$ of dimension $d$ which is adapted to $\GG$. Then, for any $\GG \vee \HH$-progressively measurable process $\bm{\eta}$ of dimension $d$ satisfying $\E\big[\int_0^T|\bm{\eta}_t|^2\d t\big] < \infty$, we have:
\begin{align*}
\E\left[\left. \int_0^t \bm{\eta}_s \cdot \d\bm{B}_s \,\right|\, \G_t\right] = \int_0^t \E[\bm{\eta}_s\,|\,\G_s] \cdot \d\bm{B}_s \quad\text{and}\quad \E\left[\left. \int_0^t \bm{\eta}_s \cdot \d\bm{W}_s \,\right|\, \G_t\right] = 0
\end{align*}
a.s., for each $t \in [0,T]$.
\end{lemma}

\noindent\textbf{Proof.} By linearity and a routine approximation in $L^2([0,T] \times \Omega, \d t \times \mathrm{d}\PP)$ by simple predictable processes, it suffices to prove the claims in the case $\bm{\eta}_u=Z\,\mathbf{1}_{[r,s]}(u)$, where $0\le r \le s \le t$ and $Z$ is $\G_r \vee \H_r$-measurable. For the first claim, use property (2b) from Lemma \ref{le:compatibility} to deduce that $Z \indep \G_T \,|\,\G_r$, which yields
\begin{align*}
\E\left[\left. \int_0^t \bm{\eta}_u \cdot \d\bm{B}_u \,\right|\, \G_t\right] &= \E[ Z(\bm{B}_s-\bm{B}_r) \,|\, \G_t] = \E[ Z \,|\, \G_t](\bm{B}_s-\bm{B}_r) = \E[ Z \,|\, \G_r](\bm{B}_s-\bm{B}_r).
\end{align*}
For the second claim, apply property (4a) from Lemma \ref{le:compatibility} to get
\begin{align*}
\qquad\;\;\E\left[\left. \int_0^t \bm{\eta}_u \cdot \d\bm{W}_u \,\right|\, \G_t\right] &= \E[ Z(\bm{W}_s-\bm{W}_r) \,|\, \G_t] = \E[\bm{W}_s-\bm{W}_r]\,\E[ Z \,|\, \G_t] = 0. \qquad\;\;\qed
\end{align*}

\smallskip

\noindent\textbf{Proof of Proposition \ref{pr:SDEtoSPDE}.}
First note that we may view $\L(\bm X\,|\,\G_T)$ as a (regular version of) the conditional law of the $C([0,T];\R^d)$-valued random variable $\bm X$ given $\G_T$, which explains why $(\mu_t:=\L(\bm{X}_t\,|\,\G_T))_{t \in [0,T]}$ admits a continuous version.
We use Lemma \ref{le:compatibility} to deduce that $\F^{\bm{X}}_t \indep \G_T \,|\, \G_t$, for each $t \in [0,T]$, which then implies claim (1). Indeed, this conditional independence yields $\mu_t=\L(\bm{X}_t\,|\,\G_t)$ a.s. More generally, we deduce that
\begin{align}
\mu_s = \L(\bm{X}_s\,|\,\G_t)\;\;\text{a.s.,} \quad 0 \le s \le t \le T. \label{pf:SDEtoSPDE1}
\end{align}
Now, starting from the SDE \eqref{SDE_gen}, apply It\^o's formula to $\varphi(\bm{X}_t)$, where $\varphi \in C^\infty_c(\R^d)$, to get
\begin{align}
\begin{split}
\varphi(\bm{X}_t) =\ & \varphi(\bm{X}_0) + \int_0^t\gen_{s,\omega}^{\bm{b},\bm{a}}\varphi(\bm{X}_s)\d s + \int_0^t(\nabla\varphi)^\top(\bm{X}_s)\bm{\gamma}(s,\omega,\bm{X}_s)\d \bm{B}_s \\
	&+  \int_0^t(\nabla\varphi)^\top(\bm{X}_s)\bm{\sigma}(s,\omega,\bm{X}_s)\d \bm{W}_s, \quad t \in [0,T].
\end{split} \label{pf:SDEtoSPDE2}
\end{align}
For fixed $t \in [0,T]$, we will take conditional expectations with respect to $\G_t$. Noting that $\mu_t=\L(\bm{X}_t\,|\,\G_t)$ a.s.~and $\mu_0=\L(\bm{X}_0\,|\,\G_t)$ a.s., the first two terms become
\begin{align*}
\E[\varphi(\bm{X}_t)\,|\,\G_t] = \langle \mu_t, \varphi\rangle, \quad \E[\varphi(\bm{X}_0)\,|\,\G_t] = \langle \mu_0, \varphi\rangle.
\end{align*}
In view of Fubini's theorem and \eqref{pf:SDEtoSPDE1}, the third term yields
\begin{align*}
\E\left[\left.\int_0^t\gen_{s,\omega}^{\bm{b},\bm{a}}\varphi(\bm{X}_s)\d s \,\right|\,\G_t\right] &= \int_0^t\E[\gen_{s,\omega}^{\bm{b},\bm{a}}\varphi(\bm{X}_s) \, | \, \G_t]\d s = \int_0^t\langle\mu_s, \gen_{s,\omega}^{\bm{b},\bm{a}}\varphi \rangle \d s.
\end{align*}
For the $\d \bm{W}_s$ and $\d \bm{B}_s$ integrals we use a stochastic Fubini theorem in the form of Lemma \ref{le:fubini}, with $\HH=\FF^{\bm{X}}$, to get
\begin{align*}
\E\left[\left. \int_0^t(\nabla\varphi)^\top(\bm{X}_s)\bm{\sigma}(s,\omega,\bm{X}_s)\d \bm{W}_s \,\right|\,\G_t\right] &= 0
\end{align*}
and
\begin{align*}
\E\left[\left. \int_0^t(\nabla\varphi)^\top(\bm{X}_s)\bm{\gamma}(s,\omega,\bm{X}_s)\d \bm{B}_s \,\right|\,\G_t\right] &= \int_0^t\E[(\nabla\varphi)^\top(\bm{X}_s)\bm{\gamma}(s,\omega,\bm{X}_s) \,|\, \G_s]\d \bm{B}_s \\
	&= \int_0^t\langle \mu_s, (\nabla\varphi)^\top \bm{\gamma}(s,\omega,\cdot)\rangle \d \bm{B}_s.
\end{align*}
Note that \eqref{integrability-SDEtoSPDE} ensures that all of the above integrals and expectations are well-defined.
Putting it together, taking conditional expectations in \eqref{pf:SDEtoSPDE2} leads to
\begin{align*}
\langle \mu_t, \varphi\rangle &= \langle \mu_0, \varphi\rangle + \int_0^t\langle\mu_s, L_{s,\omega}^{\bm{b},\bm{a}}\varphi \rangle \d s + \int_0^t\langle \mu_s, (\nabla\varphi)^\top \bm{\gamma}(s,\omega,\cdot)\rangle \d \bm{B}_s.
\end{align*}
In differential form, this is exactly the SPDE \eqref{SPDE}. \qed

\section{Strong existence for \eqref{SDE_gen} in the setting of Proposition \ref{Thm_SPDEtoSDE_unismooth}}\label{sec_appendix}

Let $\bm X^{(0)}\equiv\bm X_0$ on $[0,T/N]$ and define $\bm X^{(n)}$ on $[0,T/N]$ by
\begin{equation*}
\bm X^{(n)}_t=\bm X_0+\int_0^t \bm b\big(s,\omega,\bm X^{(n-1)}_s\big)\d s+\int_0^t 
\bm\sigma\big(s,\omega,\bm X^{(n-1)}_s\big)\d\bm W_s+\int_0^t \bm\gamma\big(s,\omega,\bm X^{(n-1)}_s\big)\d\bm B_s.
\end{equation*} 
Setting
\begin{equation*}
\Delta_n=\E\Big[\sup_{[0,T/N]}\big|\bm X^{(n+1)}-\bm X^{(n)}\big|^{ p }\Big],
\end{equation*}
we have:
\begin{equation*}
\begin{split}
			\Delta_n\leq &\;3^{ p -1}\,\E\bigg[\sup_{t\in[0,T/N]}\bigg|\int_0^{t} \bm b\big(s,\omega,\bm X^{(n)}_s\big)-\bm b\big(s,\omega,\bm X^{(n-1)}_s\big)\d s\bigg|^{ p }\,\bigg]
			\\
			&+3^{ p -1}\,\E\bigg[\sup_{t\in[0,T/N]}\bigg|\int_0^{t}\bm \sigma\big(s,\omega,\bm X^{(n)}_s\big)-\bm \sigma\big(s,\omega,\bm X^{(n-1)}_s\big)\d\bm W_s\bigg|^{ p }\,\bigg]
			\\
			&+3^{ p -1}\,\E\bigg[\sup_{t\in[0,T/N]}\bigg|\int_0^{t} \bm \gamma\big(s,\omega,\bm X^{(n)}_s\big)-\bm \gamma\big(s,\omega,\bm X^{(n-1)}_s\big)\d\bm B_s\bigg|^{ p }\,\bigg].
\end{split}
\end{equation*}
By Jensen's inequality, 
\begin{equation*}
\begin{split}
&\E\bigg[\sup_{t\in[0,T/N]}\bigg|\int_0^{t} \bm b\big(s,\omega,\bm X^{(n)}_s\big)
-\bm b\big(s,\omega,\bm X^{(n-1)}_s\big)\d s\bigg|^{ p }\,\bigg] \\
& \ \ \leq \bigg(\frac{T}{N}\bigg)^{ p -1}\,\E\bigg[\sup_{[0,T/N]} \big|\bm X^{(n)}-\bm X^{(n-1)}\big|^{ p }\int_0^{T/N}\big\|\bm b(s,\omega,\cdot)\big\|^{ p }_{C_b^1(\R^d)}\d s\bigg]
\leq D_{ p }\,T^{ p -1}\,\Delta_{n-1}.
\end{split}
\end{equation*}
Moreover, thanks to the Burkholder-Davis-Gundy and the Jensen inequalities, 
\begin{equation*}
\begin{split}
&\E\bigg[\sup_{t\in[0,T/N]} \bigg|\int_0^{t} \bm \sigma\big(s,\omega,\bm X^{(n)}_s\big)
-\bm \sigma\big(s,\omega,\bm X^{(n-1)}_s\big)\d\bm W_s\bigg|^{ p }\,\bigg] \\
& \ \ \leq C_{ p ,d}\,\E\bigg[\bigg(\int_0^{T/N}\Big|\bm \sigma\big(s,\omega,\bm X^{(n)}_s\big)
-\bm \sigma\big(s,\omega,\bm X^{(n-1)}_s\big)\!\Big|^2\d s\bigg)^{ p /2}\bigg] \\
& \ \ \leq C_{ p ,d}\,\bigg(\frac{T}{N}\bigg)^{\frac{ p -1}{2}}\,\E\bigg[\sup_{[0,T/N]} \big|\bm X^{(n)}-\bm X^{(n-1)}\big|^{ p }\,\bigg(\int_0^{T/N}\big\|\bm \sigma(t,\omega,\cdot)\big\|^{2 p }_{\mathrm {Lip}(\R^d)}\d t\bigg)^{1/2}\bigg] \\
& \ \ \leq \sqrt{(8d^2)^{ p }D_{ p }\,T^{ p -1}}\,C_{ p ,d}\,\Delta_{n-1},
\end{split}
\end{equation*}
where $C_{ p ,d} < \infty$ is a constant depending only on $ p $ and $d$.
Another application of the Burkholder-Davis-Gundy and the Jensen inequalities yields
\begin{equation*}
\E\bigg[\sup_{t\in[0,T/N]} \bigg|\int_0^{t} \bm \gamma\big(s,\omega,\bm X^{(n)}_s\big)
-\bm \gamma\big(s,\omega,\bm X^{(n-1)}_s\big)\d\bm B_s\bigg|^{ p }\,\bigg]\leq \sqrt{D_{ p }\,T^{ p -1}}\,C_{ p ,d}\,\Delta_{n-1}.
\end{equation*}
With
\begin{equation*}
D_{ p }:=\min\bigg(\frac{1}{6}\cdot\frac{1}{3^{ p -1}\,T^{ p -1}},\,\bigg(\frac{1}{6}\cdot\frac{1}{3^{ p -1}\,C_{ p ,d}}\bigg)^2\frac{1}{(8d^2)^{ p }\,T^{ p -1}}\bigg),
\end{equation*}
we have:
\begin{equation*}
3^{ p -1}\,D_{ p }\,T^{ p -1}\leq \frac{1}{6}, \quad
3^{ p -1}\,C_{ p ,d}\,\sqrt{(8d^2)^{ p }\,D_{ p }\,T^{ p -1}}
\leq \frac{1}{6},\quad 
3^{ p -1}\,C_{ p ,d}\,\sqrt{D_{ p }\,T^{ p -1}}\leq \frac{1}{6}.
\end{equation*}
Thus, $\Delta_n\leq\frac{1}{2}\Delta_{n-1}$. 
	
\medskip	

Similarly, for $\Delta_1$ and upon increasing the value of $C_{p,d}$ if necessary, 
\begin{equation*}
\begin{split}
&\E\bigg[\sup_{t\in[0,T/N]} \bigg|\int_0^{t} \bm b\big(s,\omega,\bm X^{(1)}_s\big) - \bm b\big(s,\omega,\bm X^{(0)}_s\big) \d s\bigg|^{ p }\,\bigg] \\
&\qquad \leq 2^{ p -1}\bigg(\frac{T}{N}\bigg)^{ p -1}\,\E\bigg[\int_0^{T/N}\big\|\bm b(t,\omega,\cdot)\big\|^{ p }_{C_b(\R^d)}\d t\bigg]\leq 2^{ p -1}D_{ p }\,T^{ p -1}, \\
&\E\bigg[\sup_{t\in[0,T/N]}\bigg|\int_0^{t} \bm \sigma\big(s,\omega,\bm X^{(1)}_s\big) - \bm \sigma\big(s,\omega,\bm X^{0)}_s\big) \d\bm W_s\bigg|^{ p }\,\bigg] \\
&\qquad \leq 2^{ p -1}\,C_{ p ,d}\,
\bigg(\frac{T}{N}\bigg)^{\frac{ p -1}{2}}\,\E\bigg[\int_0^{T/N}\big\|\bm \sigma^2(s,\omega,\cdot)\big\|^{ p }_{C_b(\R^d)}\d s\bigg]^{1/2}
\leq 2^{ p -1}\, C_{ p ,d}\,\sqrt{D_{ p }\,T^{ p -1}}, \\
&\E\bigg[\sup_{t\in[0,T/N]}\bigg|\int_0^{t} \bm \gamma\big(s,\omega,\bm X^{(1)}_s\big) - \gamma\big(s,\omega,\bm X^{(0)}_s\big) \d\bm B_s\bigg|^{ p }\,\bigg] \\
&\qquad \leq 2^{ p -1}\,C_{ p ,d}\,\bigg(\frac{T}{N}\bigg)^{\frac{ p -1}{2}}\,\E\bigg[\int_0^{T/N}\big\|\bm \gamma(s,\omega,\cdot)\big\|^{2 p }_{C_b(\R^d)}\d s\bigg]^{1/2}
\leq 2^{ p -1}\, C_{ p ,d}\,\sqrt{D_{ p }\,T^{ p -1}}.
\end{split}
\end{equation*}
Hence, $\Delta_1\leq 2^{ p -1}$ and, by induction, $\Delta_n\leq 2^{ p -n}$. This and Markov's inequality give 
\begin{equation*}
\sum_{n=1}^\infty \, \PP\Big(\sup_{[0,T/N]}\big|\bm X^{(n+1)}-\bm X^{(n)}\big|^{ p }\ge 2^{-n/2}\Big)\leq\sum_{n=1}^\infty 2^{ p -n/2}<\infty.
\end{equation*}
We now use the Borel-Cantelli lemma to extract an a.s.~uniform limit of $(\bm X^{(n)})_{n\in\N}$ on $[0,T/N]$,
and we easily verify that the latter is a strong solution of \eqref{SDE_gen} on $[0,T/N]$. Iterating the construction we build a strong solution of \eqref{SDE_gen} on $[0,T]$.

\bigskip\bigskip

\bibliographystyle{abbrv}
\bibliography{Superposition}

\bigskip\bigskip\bigskip\bigskip
\end{document}